\documentclass[11pt,a4paper]{article}
\usepackage[margin=1in]{geometry}
\usepackage{tabularx, enumerate, authblk,soul}
\usepackage[dvipsnames]{xcolor}
\usepackage{amsmath,amssymb,dsfont,subcaption}
\usepackage{mathtools, graphicx,color,dsfont,amsthm,bbm, mathrsfs}
\usepackage{lipsum,bm,comment}
\usepackage{natbib}
\usepackage{longtable}
\usepackage{array}
\usepackage{caption}
\usepackage[colorlinks=true,linkcolor=blue,anchorcolor=blue,citecolor=blue,filecolor=black,menucolor=black,runcolor=black,urlcolor=black]{hyperref}
\usepackage{cleveref}

\theoremstyle{plain}
\newtheorem{Theor}{Theorem}[section]
\newtheorem{Lem}[Theor]{Lemma}

\newtheorem{Prop}[Theor]{Proposition}
\newtheorem{Corol}[Theor]{Corollary}

\theoremstyle{definition}
\newtheorem{Defin}[Theor]{Definition}

\newtheorem{Examples}[Theor]{Examples}

\newcommand{\C}{\mathbb{C}}

\newcommand{\N}{\mathbb{N}}

\newcommand{\R}{\mathbb{R}}
\renewcommand{\SS}{\mathbb{S}}
\newcommand{\BB}{\mathbb{B}}

\newcommand{\CC}{\mathcal{C}}
\newcommand{\D}{\mathcal{D}}

\newcommand{\SC}{\mathcal{S}}

\newcommand{\F}{\mathcal{F}}
\newcommand{\FF}{\mathcal F}
\newcommand{\FT}{\mathcal{F}}

\newcommand{\HH}{\mathcal{H}}
\newcommand{\LL}{\mathcal{L}}

\newcommand{\RR}{\mathcal{R}}

\newcommand{\I}{\mathbb{I}}

\newcommand{\ie}{\emph{i.e.\ }}

\newcommand{\g}{{\rm g}}
\newcommand{\Fpg}{F_P^\g}
\newcommand{\Qpg}{Q_P^\g}
\newcommand{\braS}{{\SC', \SC}}
\newcommand{\braD}{{\D', \D}}
\newcommand{\PV}{{\rm PV}}

\newcommand{\E}{\mathbb{E}}

\newcommand{\Ind}[1]{\mathbb{I}\left[#1\right]}

\newcommand{\ps}[2]{\left\langle #1,#2 \right\rangle}
\newcommand{\sm}{\setminus}

\DeclareMathOperator{\supp}{supp}

\DeclareMathOperator{\Reel}{Re}

\newcommand{\ve}{\varepsilon}
\newcommand{\Om}{\Omega}

\newcommand{\loc}{{\rm loc}}
\newcommand{\spa}{\quad\quad}

\usepackage{scalerel,stackengine}
\stackMath
\newcommand\wwidehat[1]{%
	\savestack{\tmpbox}{\stretchto{%
			\scaleto{%
				\scalerel*[\widthof{\ensuremath{#1}}]{\kern-.6pt\bigwedge\kern-.6pt}%
				{\rule[-\textheight/2]{1ex}{\textheight}}%WIDTH-LIMITED BIG WEDGE
			}{\textheight}% 
		}{0.5ex}}%
	\stackon[1pt]{#1}{\tmpbox}%
}

\newcommand{\RN}[1]{%
	(\textup{\uppercase\expandafter{\romannumeral#1}})%
}
\newcolumntype{H}{>{\setbox0=\hbox\bgroup}c<{\egroup}@{}}

\newenvironment{Proof}[1]
{
	\begin{proof}[Proof of #1]
}
{
	\end{proof}
}

\allowdisplaybreaks

\title{PDE characterisation of geometric distribution functions and quantiles}
\author{Dimitri Konen}
\date{\today\\ First version: 24 August, 2022}

\begin{document}

\maketitle

\begin{abstract}
	We show that in any Euclidean space, an arbitrary probability measure can be reconstructed explicitly by its geometric (or spatial) distribution function. The reconstruction takes the form of a (potentially fractional) linear PDE, where the differential operator is given in closed form. This result implies that, contrary to a common belief in the statistical depth community, geometric cdf's in principle provide exact control over the probability content of all depth regions. We present a comprehensive study of the regularity of the geometric cdf, and show that a continuous density in general does not give rise to a geometric cdf with enough regularity to reconstruct the density pointwise. Surprisingly, we prove that the reconstruction displays different behaviours in odd and even dimension: it is local in odd dimension and completely nonlocal in even dimension. We investigate this issue and provide a partial counterpart for even dimensions, and establish a general representation formula of the geometric cdf of spherically symmetric probability laws in odd dimension. We provide explicit examples of reconstruction of a density from its geometric cdf in dimension 2 and 3.
\end{abstract}
%
%	\setcounter{tocdepth}{1}
%
%	{
%			\hypersetup{linkcolor=black}
%			\tableofcontents
%		}

	\section{Introduction} \label{sec:Introduction}
	
	%	\subsection{Context of this work}
	
	The cdf and quantile maps of univariate distributions play a vital role in statistics and probability. For instance, statistical procedures that combine broad validity (no need for moment assumption, resistance to possible outlying observations) and efficiency are typically based on ranks, computed by evaluating the cdf at observed data points. That is, each data point is used only through the fact it is the smallest one, the second smallest, etc. These procedures found many applications in hypothesis testing, outliers detection, and extreme value theory, for instance. In the multivariate setting $\R^d$ with $d\geq 2$, however, no canonical ordering is available, so that there is no natural concept of cdf that can be used to define cdf-based statistical procedures. In this context, \emph{statistical depth} is a general device that allows one to define a center-outward ordering of data points in $\R^d$, thus providing the basis for the definition of suitable multivariate rank-based procedures. Since the introduction of the celebrated \emph{halfspace depth} in \cite{Tuk1975}, and the subsequent work of Regina Liu that made statistical depth a field in itself (\cite{Liu1990}, \cite{Liuetal1999}), 
	statisticians have proposed many concepts, called \emph{multivariate quantiles, cdf's, and depth}, to extend these ideas to a multivariate framework. The most celebrated depths are the aforementioned halfspace depth, the \emph{simplicial depth} \cite{Liu1990}, the \emph{geometric (or spatial) depth} \cite{VarZha2000}, and the \emph{projection depth} \cite{Zuo2003}.
	%In the univariate case $d=1$, cdfs are provided by the cdf through evaluation at observed data points. Therefore, the notions of cdfs, cdf, and quantiles are closely related.
	% The concept of \emph{halfspace depth}, introduced in \cite{Tuk1975}, is arguably the most famous attempt to establish a multivariate analogue of quantiles. Many other depth notions followed, among which the \emph{simplicial depth} \cite{Liu1990}, the \emph{spatial (or geometric) depth} \cite{VarZha2000}, and the \emph{projection depth} \cite{Zuo2003} to cite only a few. 
	%
	%In the univariate case $d=1$, cdfs are obtained by evaluating the empirical cdf of the data at the observed data points. Therefore,  For this reason, a lot of effort has been made over the past decades to extend the concepts of cdf and quantiles to a multivariate setting. It easy to notice that, for a probability measure $P$ on $\R^d$ with $d\geq 2$, considering quantiles of $P$ component-wise leads to quantile regions which are not equivariant under orthogonal transformations. For obvious geometric reasons, it automatically discarded this approach. In this context, \emph{statistical depth} is a general device that allows one to define a center-outward ordering of data points in $\R^d$, thus providing the basis for the definition of suitable multivariate cdf-based procedures . 
	%To any point of $\R^d$ it associates a non-negative number, its \emph{depth}. The regions of points the depth of which does not exceed a given treshold value are interpreted as analogues of quantile regions in $\R$. 
	Other approaches have been adopted, attempting to define a proper notion of multivariate cdf and quantiles; the most notable are based on regression quantiles \cite{Haletal2010}, or optimal transport \cite{HallEspagne21}. We also refer the reader to \cite{Ser2002C} for a review on the topic.

	Among the concepts extending cdf's and quantiles to a multivariate setting, a popular approach is that of \textit{geometric} multivariate quantiles and cdf's, introduced in \cite{Cha1996}, on which this paper focuses. They enjoy important advantages over other competing approaches. Among them, let us stress that the geometric cdf is available in closed form, which leads to trivial evaluation in the empirical case, unlike most competing concepts. As a consequence, explicit Bahadur-type representations and asymptotic normality results are provided in \cite{Cha1996} and \cite{Kol1997}, when competing approaches offer at best consistency results only; see also \cite{PasReid2022}. More recently, \cite{GirStu2015} and \cite{GirStu2017} derived a fine asymptotic behaviour of extreme geometric quantiles. Finally, these concepts can be extended to infinite-dimensional Banach spaces; see, e.g., \cite{Kem1987}, \cite{ChaCha2014}, \cite{Chow19}, and \cite{Romon2022}. We also refer the reader to \cite{Oja2010} for some applications of geometric quantiles and cdf's.
	
	Before addressing the main question of the paper, we recall basic definitions and standard results about geometric quantiles and cdf's. We start by introducing a closely related map that will play an important role in our proofs.
	
	%	We denote the Euclidean inner product between two vectors $u=(u_1,\ldots, u_d)\in\R^d$ and $v=(v_1,\ldots, v_d)\in\R^d$ by $(u,v):=\sum_{i=1}^d u_i v_i$, $\|u\|:=\sqrt{(u,u)}$ for the Euclidean norm. We denote the unit sphere of~$\R^d$ by $\SS^{d-1}=\{x\in\R^d : (x,x)=1\}$. 

	\begin{Defin}\label{DefinObjectiveFunction}
		Let $d\geq 1$ and $P$ be a Borel probability measure on $\R^d$. Let $Z$ be a random $d$-vector with law $P$. We define the map $h_P : \R^d\to\R$ by letting 
		$$
		h_P(x)
		=
		\E\big[\|x-Z\|-\|Z\|\big] 
		%		\int_{\R^d} \big( \|z-x\|-\|z\| \big)\, dP(z)
		,\spa \forall\ x\in\R^d
		.
		$$
		%	for any $x\in\R^d$.
	\end{Defin}
	
	The triangle inequality entails that $h_P$ is well-defined and continuous over $\R^d$, irrespective of the probability measure $P$, without any moment assumption.
	%Theorem 5.1 in \cite{KonPai1} implies that $h_P$ is continuously differentiable over an open subset $U\subset\R^d$ if and only if $P$ has no atoms over $U$. In that case, we have $\Fpg(x)=\nabla h_P(x)$ for any $x\in U$, where $\Fpg$ is given in the next definition. 
	Theorem 5.1 in \cite{KonPai1} further entails that $h_P$ is continuously differentiable over an open subset~$\Om\subset\R^d$ if and only if $P$ has no atoms in $\Om$; in this case, we have
	$$
	\nabla h_P(x)
	=
	\E\bigg[\frac{x-Z}{\|x-Z\|}\I[Z\neq x]\bigg]
	%	\int_{\R^d\sm\{x\}} \frac{x-z}{\|x-z\|}\, dP(z)
	,\spa \forall\ x\in\R^d
	.
	$$
	%for any $x\in\R^d$.\\
	
	We now turn to the definition of multivariate geometric quantiles, which are directional in nature: they are indexed by an order $\alpha\in [0,1)$ and a direction $u\in\SS^{d-1}:=\{x\in\R^d : \|x\| = 1\}$. 
	
	\begin{Defin}\label{DefinQuantiles}
		Let $d\geq 1$ and $P$ be a probability measure on $\R^d$. A geometric quantile of order $\alpha\in [0,1)$ in direction $u\in \SS^{d-1}$ for $P$ is an arbitrary minimizer of the objective function 
		$$
		O_{\alpha, u}^P : \R^d\to \R,\
		x\mapsto 
		O_{\alpha,u}^P(x)
		:=
		h_P(x)-\ps{\alpha u}{x}
		,
		%	\int_{\R^d} \Big\{|z-x|-\|z\|- (\alpha u,x)\Big\}\, dP(z)
		$$
		where $\ps{\alpha u}{x}$ denotes the inner product of $\R^d$ between $\alpha u$ and $x$.
	\end{Defin} 
	
	For $d=1$, any geometric quantile of order $\alpha\in[0,1)$ in direction~$u\in\{-1,+1\}$ reduces to a usual univariate quantile of order~$(\alpha u + 1)/2\in (0,1)$.
	
	As a consequence of the remarks preceeding Definition \ref{DefinQuantiles}, we have
	$
	\nabla O_{\alpha,u}^P(x)
	=
	\nabla h_P(x)-\alpha u 
	$
	for any $x\in\R^d$ as soon as $P$ has no atoms. Further requiring that $P$ is not supported on a single line of $\R^d$, 
	%	\textcolor{red}{REF Theorem 2.17 in Kemperman (1987)}
	Theorem 1 in \cite{PaiVir2021} entails that $O_{\alpha, u}^P$ is strictly convex over $\R^d$ and, therefore, that the geometric quantile of order $\alpha$ in direction $u$ for $P$ is unique for any $\alpha\in [0,1)$ and $u\in \SS^{d-1}$; we write such a quantile~$\Qpg(\alpha u)$. In particular, $x:=\Qpg(\alpha u)$ is the unique solution to the equation $\nabla O_{\alpha,u}^P(x)=0$, \ie 
	$$
	%	\nabla h_P(x)
	%	=
	\E\bigg[\frac{x-Z}{\|x-Z\|}\I[Z\neq x]\bigg]
	=
	\alpha u
	.
	$$
	Under these assumptions ($P$ is non-atomic and not supported on a line),  Theorem 6.2 in \cite{KonPai1} entails that the geometric quantile map $\alpha u\mapsto \Qpg(\alpha u)$ is invertible with inverse~$(\Qpg)^{-1}=\nabla h_P$. In anology with the univariate case, the equality $(\Qpg)^{-1} = \nabla h_P$ provides a good motivation to call $\nabla h_P$ the \emph{multivariate geometric cdf of $P$}.
	
	\begin{Defin}\label{DefinRank}
		Let $d\geq 1$ and $P$ be a Borel probability measure on $\R^d$. The geometric cdf $\Fpg$ of $P$ is the map $\Fpg:\R^d\to\R^d$ defined by letting
		$$
		\Fpg(x)
		=
		\E\bigg[\frac{x-Z}{\|x-Z\|}\I[Z\neq x]\bigg]
		,\spa 
		\forall\ x\in\R^d
		.
		%	=
		%	\int_{\R^d\sm\{z\}} \frac{x-z}{\|x-z\|}\, dP(z)
		%	=
		%	(K*P)(x)
		$$
		%	for any $x\in\R^d$.
	\end{Defin}
	
	Observe that this definition still makes sense when $h_P$ is not differentiable, \ie when $P$ has atoms; only the equality $\Fpg=\nabla h_P$ requires $P$ to be non-atomic.	When $d=1$, denoting by $F_P$ the usual univariate cdf of $P$, we have
	$$
	\Fpg(x)
	=
	\E\big[{\rm sign}(x-Z)\I[Z\neq x]\big]
	%	\int_{\R\sm\{z\}} \textrm{sign}(x-z)\, dP(z)
	=
	2F_P(x)-1
	,\spa \forall\ x\in\R
	.
	$$
	Therefore, $\Qpg$ and $\Fpg$ qualify as multivariate extensions of the univariate notions of quantile map and cdf, respectively. Like its univariate counterpart, the geometric cdf characterises probability measures in any dimension: if $P$ and $Q$ are Borel probability measures on $\R^d$, and if $\Fpg(x) = F_Q^\g(x)$ for all $x\in\R^d$, then $P=Q$; see Theorem 2.5 and Corollary 2.9 in \cite{Kol1997}. This very desirable property is also shared by the so-called \emph{center-outward cdf} introduced in \cite{HallEspagne21} and based on optimal transport; in addition, when $P$ admits a sufficiently smooth density, the density is linked to its center-outward cdf through a non-linear partial differential equation. However, the characterisation property is not shared by the Tukey halfspace depth; see \cite{Nag2021}. Nevertheless, halfspace depth possesses this property within some classes of probability measures; see, e.g., \cite{StrRou1999} who gave the first positive result for empirical probability measures by algorithmically reconstructing the measure. We refer to \cite{Nagy2020} for a recent review on the characterisation property of this depth concept. These results provided a strong motivation to explore this question for the geometric cdf: since the characterization property established in \cite{Kol1997}, the existence of an operational procedure to reconstruct $P$ from $\Fpg$ has been left unanswered. The present work completely solves this question by providing the explicit way to recover $P$ from $\Fpg$.
	
	%	\subsection{Our contribution}
	%	While computing $\Fpg$ from $P$ is straightforward from the definition, it is not easy to understand what information about $P$ can be extracted from $\Fpg$. 
	%	This question remained unanswered for 25 years, and was completely solved by this work.\\
	Through a Fourier transform approach inspired by \cite{Kol1997}, we prove in Section \ref{sec:PDEGeneral} that the link between an \emph{arbitrary} probability measure $P$ and its geometric cdf $\Fpg$ is given by the (potentially fractional) linear partial differential equation $P=\gamma_d\ (-\Delta)^{(d-1)/2}(\nabla\cdot \Fpg)$, where $\gamma_d$ is an explicit dimensional constant; we provide a brief introduction to fractional Laplacians, to which the previous PDE reduces when $d$ is even, in Section \ref{sec:IntroLaplace}. While the PDE relating $P$ and $\Fpg$ holds in the sense of distributions for any probability measure $P$, we prove in Section \ref{sec:Regularity} that when $P$ admits a density $f_P$ the equality holds pointwise when $P$ is replaced by $f_P$. Surprisingly, the nature of the reconstruction procedure depends on whether $d$ is odd or even; we show in Section \ref{sec:Localisation} that when $d$ is odd the PDE is local in nature: if $\Fpg=F_Q^\g$ on an open set $\Om\subset\R^d$, then $P$ and $Q$ coincide on $\Om$ as well, whereas if $d$ is even and $P$ and $Q$ coincide on $\Om$, then the equality $\Fpg=F_Q^g$ on $\Om$ automatically implies that $P=Q$ on $\R^d$. When $d$ is odd and $P$ is a spherically symmetric probability measure, we establish an explicit formula for $\Fpg$ that only involves cumulative moments of $P$ of order comprised between~$-(d-2)$ and~$(d-1)$, from which we deduce that there indeed exist distinct probability measures $P$ and $Q$ such that their geometric cdf coincide on an open set. We then provide a counterpart to the nonlocality when $d$ is even by essentially embedding probability measures in higher dimensional spaces; we further show that this procedure is in fact equivalent to well-known results of localisation of fractional Laplacians used by PDE analysts. Our results allow us to determine in Section \ref{sec:DepthRegions} the regularity of geometric depth regions associated with a probability measure: we show that they are smooth manifolds of dimension $d-1$ and that their regularity is controlled by the regularity of the density. Proofs and additional comments are collected in the Appendix.  

	\subsection*{Notation}
	%	\label{sec:Notation}
	Let $\N=\{0,1,2,\ldots,\}$ stand for the natural numbers. We denote the Euclidean inner product on $\R^d$ by $\ps{\cdot}{\cdot}$, and the induced norm by $\|\cdot \|$. We write $\SS^{d-1}=\{x\in\R^d : \ps{x}{x}=1\}$ for the unit sphere of $\R^d$, and $\BB^d$ the open unit ball of $\R^d$; for any $x\in\R^d$ and $r>0$, we denote by $\BB^d_r(x)$ the open ball of radius $r$ centered at $x$---when $x=0$ we will just write $\BB^d_r$. For any subset $A\subset\R^d$, we write $\overline{A}$ the closure of $A$ with respect to the usual topology of $\R^d$. We let $\Ind{A}$ denote the indicator function of the condition $A$. 
	
	For any open subset $\Om\subset\R^d$, with $d\geq 1$, and any $k$-times differentiable function $u:\Om\to\C$ we let 
	$$
	(\partial^\alpha u)(x)
	:=
	\frac{\partial^{|\alpha|} u}{\partial x_1^{\alpha_1} \ldots \partial x_d^{\alpha_d}}(x)
	,\spa 
	\forall\ x\in \Om,
	$$
	for any multi-index $\alpha=(\alpha_1,\ldots, \alpha_d)\in\N^d$ such that $|\alpha|:=\sum_{j=1}^d \alpha_j \leq k$. By convention, we let $\partial^\alpha u := u$ if $\alpha=(0,\ldots, 0)$. In addition, we write
	\begin{itemize}
		%			\item $u\in\CC^k(U)$ if $u$ is $k$-times differentiable and such that $\partial^\alpha u$ is continuous over $U$ for any $\alpha\in\N^d$ with $|\alpha|\leq k$;
		\item $u\in \CC^k_b(\Om)$ when $u\in\CC^k(\Om)$ and $\partial^\alpha u$ is bounded over $\Om$ for any $\alpha\in\N^d$ with $|\alpha|\leq k$;
		
		\item $u\in \CC^k_c(\Om)$ when $u\in\CC^k(\Om)$ has a compact support contained in $\Om$---the set $\CC^\infty_c(\Om)$ of infinitely differentiable maps with compact support in $\Om$ is also denoted $\D(\Om)$;
		
		\item $u\in \CC^{k,\alpha}(\Om)$ for some $\alpha \in (0,1]$ when $u\in\CC^k(\Om)$, $\partial^\beta u$ is bounded over $\Om$ for any $\beta\in\N^d$ with $|\beta|\leq k$, and  $\partial^\beta u$ is $\alpha$-H{\"o}lder continuous over $\Om$ when $|\beta|=k$;
		
		\item $u\in\CC_0(\R^d)$ when $u$ is continuous and converges to $0$ at infinity.
	\end{itemize}
	When $V$ is a collection of functions $u:\mathcal{T}\to \C^d$ defined over a topological space $\mathcal{T}$, we let $V_{\textrm{loc}}$ denote the collection of functions $u:\mathcal{T}\to\C^d$ such that the restriction $u_{|K}$ of $u$ to any compact set $K\subset \mathcal{T}$ belongs to $V$; for instance, we will consider the spaces $L^p_\loc(\Om)$.
	
	We denote by $\D(\Om)$---or, equivalently, $\CC^\infty_c(\Om)$---the set of infinitely differentiable maps whose support is compact and included in $\Om$, and by $\SC(\R^d)$ the Schwartz class. We denote by $\D(\Om)'$ the set of distributions over $\Om$, and by $\SC(\R^d)'$ the set of tempered distributions; see Section \ref{sec:Distributions} of the Appendix for a short review on the theory of distributions. We will denote by $\ps{\cdot}{\cdot}_\braD$---or $\ps{\cdot}{\cdot}_\braS$, depending on the context---the duality bracket of distributions, \ie we write the value of a distribution $T$ computed at the test function $\varphi$ as $\ps{T}{\varphi}_\braD$.
	
	When $T\in\D(\Om)'$ is a $\C$-valued distribution---or a function, which is just a particular case---and for an integer $\ell\geq 0$, we write $\partial_i^\ell T$ for the partial derivative of $T$ taken $\ell$-times with respect to the $i$th coordinate, $i=1,2,\ldots, d$. We refer the reader to Section \ref{sec:Distributions} of the Appendix for the precise meaning of $\partial_i T$ when $T$ is a general distribution. We then denote by $\nabla T = (\partial_1 T, \partial_2 T,\ldots, \partial_d T)$ the gradient of $T$ and by $\Delta T = \sum_{i=1}^d \partial_i^2 T$ the Laplacian of $T$. If $T$ is a $\C^d$-valued distribution on $\R^d$---or a vector-field on $\R^d$, which is a particular case---we let $\nabla\cdot T=\sum_{i=1}^d \partial_i T_i$ stand for the divergence of $T$, where $T_i$ denotes the $i$th component of $T$. When a map $u\in L^1_\loc(\R^d\sm\{0\})$ is such that $u\notin L^1_\loc(\BB^d)$ and the limit 
	$$
	\ps{\PV(u)}{\varphi}_{\D',\D}
	:=
	\lim_{\eta\downarrow 0} \int_{\R^d\sm\BB^d_\eta} u(x)\varphi(x)\, dx 
	$$
	exists for all $\varphi\in\CC^\infty_c(\R^d)$, we denote by $\PV(u)$ the distribution that, to any $\varphi\in\CC^\infty_c(\R^d)$, associates the above integral; we call it \emph{the principal value of $u$}.
	
	For any $u\in L^1(\R^d)$, we define the Fourier transform $\FT u$ of $u$ by letting
	$$
	(\FT u)(\xi)
	=
	\int_{\R^d} u(x)e^{-2i\pi\ps{x}{\xi}}\, dx 
	,\spa 
	\forall\ \xi\in\R^d
	.
	$$
	We also denote by $\FT$ the Fourier transform on $L^2(\R^d)$, defined as the unique continuous extension to $L^2(\R^d$) of the the Fourier transform restricted to the Schwartz class. The Fourier transform on $\SC(\R^d)'$ will be denoted by $\FT$ as well. We further let $\FT$ act componentwise on $\SC^k(\R^d)':=(\SC(\R^d)')^k$ for any integer $k\geq 1$; see Section \ref{sec:Distributions} of the Appendix.

	\section{Basics of fractional Laplacians}\label{sec:IntroLaplace}
	
	Because our main result involves a fractional Laplacian, we first present a brief construction of fractional Laplacians. This section assumes basic knowledge of distribution theory and Sobolev spaces, for which we provide a concise review in Section \ref{sec:Review} of the Appendix; see also the Notation subsection of Section \ref{sec:Introduction} for notations we will use in the present section. Different definitions of fractional Laplacians exist in the literature; some rely on the Fourier transform, others on singular integrals, or Sobolev spaces. They all coincide for smooth enough functions---such as the Schwartz class---but may differ in general, or at least have different domains. In this section, we provide a self-contained introduction to fractional Laplacians based on Fourier transforms because this approach is the one that naturally appears in our proofs. We refer the reader to \cite{Silvestre2007}, \cite{HitchhikerLaplace}, \cite{Stinga2019}, \cite{Kwa2019}, and \cite{Lischke2020} for further details.
	%Observe that 
	%$$
	%\FT(\Delta u)
	%=
	%(2i\pi \|\xi\|)^2 \FT u
	%,
	%$$
	%which we rewrite
	%$$
	%\FT((-\Delta) u)
	%=
	%(2\pi \|\xi\|)^2 \FT u
	%.
	%$$
	
	Fix $u\in\SC(\R^d)$, a Schwartz function on $\R^d$. Throughout the paper, we always implicitly denote by $x$ the space-variable---belonging to the domain of $u$---and by $\xi$ the frequency variable---belonging to the domain of $\FT u$. Recalling that the Fourier transform $\FT((-\Delta)^\ell u)(\xi)$ of the map $(-\Delta)^\ell u$ evaluated at $\xi$ satisfies, for any integer $\ell\geq 0$, 
	$$
	\FT((-\Delta)^\ell u)(\xi)
	=
	(2\pi \|\xi\|)^{2\ell} (\FT u)(\xi)
	,\spa \forall\ \xi\in\R^d
	,
	$$
	we define, for any positive real number $s>0$,
	$$
	((-\Delta)^s u)(x)
	:=
	(2\pi)^{2s} \FT^{-1}\Big(\|\xi\|^{2s} (\FT u)(\xi)\Big)(x)
	,
	\spa \forall\ x\in\R^d
	.
	$$
	%for any real $s>0$.
	%This motivates the definition of a generic operator $(-\Delta)^s$ for a real $s>0$ as the operator associated to the Fourier symbol $(2\pi\|\xi\|)^{2s}$. We therefore let 
	%$$
	%(-\Delta)^s u(x)
	%:=
	%(2\pi)^{2s} \FT^{-1}(\|\xi\|^{2s} \FT u(\xi))(x)
	%$$
	%for any $x\in\R^d$. 
	Before going further, let us make a comment about the factor $(2\pi)^{2s}$ in the definition $(-\Delta)^s u$. Obivously, this is a consequence of our choice of normalization in the definition of the Fourier transform.
	% since it derives from the fact that 
	%$$
	%\FT((-\Delta)u)
	%=
	%(2\pi\|\xi\|)^2 \FT u
	%.
	%$$
	If one chooses to work with a different normalization of the Fourier transform, $\FT_{a,b}$ say, defined by
	$$
	(\FT_{a,b} u)(\xi)
	:=
	\frac1b\int_{\R^d} u(x)e^{-ia \ps{x}{\xi}}\, dx 
	,\spa
	\forall\ \xi\in\R^d,
	$$
	for some $a>0$ and $b>0$, we instead define $(-\Delta)^s u$ as 
	$$
	((-\Delta)^s u)(x)
	:=
	a^{2s} \FT_{a,b}^{-1}\Big(\|\xi\|^{2s} (\FT_{a,b} u)(\xi)\Big)(x)
	,\spa \forall\ x\in\R^d 
	.
	$$
	It is, however, easy to see that 
	$$
	a^{2s} \FT_{a,b}^{-1}\Big(\|\xi\|^{2s} (\FT_{a,b} u)(\xi)\Big)(x)
	=
	\FT_{1,1}^{-1}\Big(\|\xi\|^{2s} (\FT_{1,1} u)(\xi)\Big)(x)
	,\spa
	\forall\ x\in\R^d
	.
	$$
	This entails that the definition of $(-\Delta)^s u$ is, in fact, independent of the normalization coefficients $a$ and $b$ of the Fourier transform. Throughout the paper, our choice corresponds to $a=2\pi$ and $b=1$.
	%	 we have
	%	$$
	%	(\FT_{a,b}^{-1}u)(\xi)
	%	=
	%	b\Big(\frac{a}{2\pi}\Big)^d \int_{\R^d} u(x)e^{ia(x,\xi)}\, dx 
	%	,
	%	\spa \forall\ \xi\in\R^d
	%	.
	%	$$
	%	It is easy to show that 
	%	$$
	%	a^{2s} \FT_{a,b}^{-1}(\|\xi\|^{2s} \FT_{a,b} u(\xi))(x)
	%	=
	%	\FT_{1,1}^{-1}(\|\xi\|^{2s} \FT_{1,1} u(\xi))(x)
	%	,\spa
	%	\forall\ x\in\R^d
	%	.
	%	$$
	%	It follows that any choice of $a$ and $b$ leads to the same value of $(-\Delta)^s u$ if we let
	%	$$
	%	((-\Delta)^s u)(x)
	%	=
	%	a^{2s} \FT_{a,b}^{-1}(\|\xi\|^{2s} \FT_{a,b} u(\xi))(x)
	%	,\spa \forall\ x\in\R^d
	%	.
	%	$$
	%	In the sequel, we will be working with $a=2\pi$ and $b=1$.\\
	
	When $s=n+\sigma$, with $n\in\N$ and $\sigma\in (0,1)$, taking the Fourier transform readily implies that 
	$$
	(-\Delta)^s u 
	=
	(-\Delta)^n((-\Delta)^\sigma u)
	=
	(-\Delta)^\sigma((-\Delta)^n u)
	,
	$$
	where $(-\Delta)^n$ is the usual differential operator $-\Delta$ taken $n$ times. When $s\in (0,1)$, Proposition 3.3 in \cite{HitchhikerLaplace} gives
	\begin{equation}\label{Eq:PointwiseFracLaplace}
		((-\Delta)^s u)(x)
		=
		c_{d,s}\ \lim_{\eta \downarrow 0}  \int_{\R^d\sm \BB^d_\eta(x)} \frac{u(x)-u(z)}{\|x-z\|^{d+2s}}\, dz
		,\spa \forall\ x\in\R^d,
	\end{equation}
	for some constant $c_{d,s}$ that only depends on $d$ and $s$. Note that the normalization of the Fourier transform used in \cite{HitchhikerLaplace} corresponds to $a=1$ and $b=(2\pi)^{d/2}$ in our previous discussion. The value of the constant $c_{d,s}$ can be found in \cite{Stinga2019} (see Theorem 1), and is given by
	\begin{equation}\label{Eq:ConstantCds}
		c_{d,s}
		=
		\frac{s(1-s)4^s \Gamma(d/2 + s)}{\Gamma(2-s)\pi^{d/2}}
		.
	\end{equation}
	
	In the paper, we will need to apply a fractional Laplacian $(-\Delta)^s$ to a function for which the previous definition of $(-\Delta)^s$ through Fourier transforms is not well-defined. Thus, let us explain how we can extend the domain of $(-\Delta)^s$ when $s\in (0,\infty)$. First, observe that Plancherel's formula, \ie permuting the Fourier transform on the functions under the integral, gives
	$$
	\int_{\R^d} ((-\Delta)^s u)(x)  v(x)\, dx 
	=
	\int_{\R^d} u(x) ((-\Delta)^s v)(x)  \, dx
	,\spa \forall\ u,v\in\SC(\R^d)
	.
	$$
	It is thus tempting to define the fractional Laplacian $(-\Delta)^s T$ of an arbitrary tempered distribution $T\in\SC(\R^d)'$
	%	---recall that Fourier transforms can only be defined for tempered distributions and that Fourier transforms are involved in the previous definition of $(-\Delta)^s$---
	by letting 
	\begin{equation}\label{eq:DualityFracLaplace}
		\ps{(-\Delta)^sT}{\psi}_\braS
		:=
		\ps{T}{(-\Delta)^s\psi}_\braS
		,\spa \forall\ \psi\in\SC(\R^d)
		.
	\end{equation}
	%	where $\ps{\Lambda}{\varphi}$ simply stands for $\Lambda(\varphi)$ for any distribution $\Lambda$ and test function $\varphi$. 
	This approach, however, must be discarded because $(-\Delta)^s \psi$ does not belong to $\SC(\R^d)$ in general so that $\ps{T}{(-\Delta)^s \psi}_\braS$ might not be defined. Consequently, we need to understand how the operation $(-\Delta)^s$ affects the regularity of test functions. This is the content of the next proposition, which is partially stated in \cite{Silvestre2007}---see after Definition 2.3---but not proved. For the sake of completness, we provide a proof of the next Proposition in Section \ref{sec:Sobolev} of the Appendix.
	
	\begin{Prop}\label{PropRegFracLaplace}
		Fix $d\geq 1$ and $s\in (0,\infty)$. Let $n=\lfloor s\rfloor\in\N$ and $\sigma=s-n\in [0,1)$. For all $u\in\SC(\R^d)$ we have $(-\Delta)^s u\in\CC^\infty(\R^d)$. Furthermore, there exists a positive constant $C_{d,\sigma}$ such that
		%		\begin{equation}\label{eq:PropRegFracLaplace}
			%			\sup_{x\in\R^d} |(1+\|x\|^{d+2s}) \partial^\alpha ( (-\Delta)^s u)(x)|
			%			<
			%			\infty 
			%%			,\spa \forall\ \alpha\in\N^d
			%%			.
			%		\end{equation}
		%		and
		%	\begin{equation}\label{eq:PropRegFracLaplace}
			%		\sup_{x\in\R^d} |(1+\|x\|^{d+2s}) \partial^\alpha ((-\Delta)^s u)(x)|
			%		\lesssim 
			%		|\partial^\alpha u|_{L^1(\R^d)}
			%		+
			%		\sup_{z\in\R^d} \Big( (1+\|z\|)^{d+2} |\nabla^2 (\partial^\alpha u)(z)| \Big)
			%		,
			%	\end{equation}
		%		\begin{equation}
			\begin{eqnarray}
				\label{eq:PropRegFracLaplace}
				\lefteqn{
					\hspace{-5mm}
					\sup_{x\in\R^d} |(1+\|x\|^{d+2\sigma}) \partial^\alpha ((-\Delta)^s u)(x)|
				}
				\\
				&&
				\hspace{5mm} 
				\leq 
				C_{d,\sigma}
				\bigg\{
				\|\partial^\alpha (-\Delta)^n u\|_{L^1(\R^d)}
				+
				\sup_{z\in\R^d} \Big( (1+\|z\|)^{d+2} \| \nabla^2 (\partial^\alpha (-\Delta)^n u)(z) \|_{\rm op} \Big)
				\bigg\}
				\nonumber
				,
			\end{eqnarray}
			%		\end{equation}
		%	\begin{eqnarray}\label{eq:PropRegFracLaplace}
			%		\lefteqn{ 
				%			%			\hspace{-10mm}
				%			\sup_{x\in\R^d} |(1+\|x\|^{d+2s}) \partial^\alpha ((-\Delta)^s u)(x)|
				%%			\lesssim 
				%%			+
				%%			|\nabla^2 (\partial^\alpha u)|_{L^\infty(\R^d)}
				%		}
			%		%		\\[2mm]
			%		\\
			%		&&
			%%		\hspace{60mm}
			%		\hspace{10mm}
			%		\lesssim
			%		|\partial^\alpha u|_{L^1(\R^d)}
			%		+
			%		\sup_{z\in\R^d} \Big( (1+\|z\|)^{d+2} |\nabla^2 (\partial^\alpha u)(z)| \Big)\nonumber
			%		,
			%	\end{eqnarray}
		holds for all $u\in\SC(\R^d)$ and $\alpha\in\N^d$, where $\|\nabla^2 (\cdot)(z)\|_{\rm op}$ stands for the operator norm of the Hessian matrix at $z$.
	\end{Prop}
	
	%We therefore need $T$ to be more regular than general tempered distributions. 
	According to Proposition \ref{PropRegFracLaplace}, defining the fractional Laplacian $(-\Delta)^s T$ of a tempered distribution $T\in\SC(\R^d)'$ by duality through (\ref{eq:DualityFracLaplace}) when $s\in (0,\infty)\sm \N$ requires $T$ to be defined not only on $\SC(\R^d)$, but on the larger space of test functions
	$$
	\SC_s(\R^d)
	:=
	\Big\{\psi \in\CC^\infty(\R^d) : \sup_{x\in\R^d} |(1+\|x\|^{d+2(s-\lfloor s\rfloor)})\partial^\alpha \psi(x)| < \infty,\ \forall \alpha\in\N^d \Big\}
	.
	$$
	We endow $\SC_s(\R^d)$ with the topology associated with the following convergence: a sequence $(\psi_k)\subset\SC_s(\R^d)$ converges to $\psi\in\SC_s(\R^d)$ in the space $\SC_s(\R^d)$ if 
	$$
	\lim_{k\to\infty} 
	\sup_{x\in\R^d} |(1+\|x\|^{d+2(s-\lfloor s\rfloor)}) \partial^\alpha (\psi_k-\psi)(x)|
	= 
	0
	,\spa \forall\ \alpha\in\N^d
	.
	$$
	
	\begin{Defin}
		Fix $s \in (0,\infty)$. We let $\SC_s(\R^d)'$ be the set of maps 
		$
		T:\SC_s(\R^d)\to\C,\ \psi\mapsto \ps{T}{\psi}_{\SC_s', \SC_s}
		$
		which are linear and continuous with respect to the topology of $\SC_s(\R^d)$. 
		%		We write the corresponding duality bracket as $\ps{T}{\psi}_{\SC_s', \SC_s}$ for all $\psi\in \SC_s(\R^d)$.
	\end{Defin}
	
	%	Note that we safely discard the case $s\in\N$ as $\SS_s(\R^d)$ is, in this case, not the appropriate space
	Proposition \ref{PropRegFracLaplace} entails that if a sequence $(\psi_k)\subset \SC(\R^d)$ converges to $0$ in the space $\SC(\R^d)$, then the sequence $((-\Delta)^s \psi_k)$ converges to $0$ in the space $\SC_s(\R^d)$ as $k\to \infty$. 
	%	Therefore, $\SC(\R^d)$ continuously embeds into $\SC_s(\R^d)$. 
	In particular, if $T\in\SC_s(\R^d)'$ and $(-\Delta)^s T$ is defined according to (\ref{eq:DualityFracLaplace}), then $(-\Delta)^s T$ is a tempered distribution. This provides the motivation for the next definition.
	
	%	is more appropriate to define fractional Laplacians by duality, as we will explain below. 
	%as $k\to \infty$ for any $\alpha\in\N^d$. 
	
	%	Because $\SC(\R^d)$ continuously embeds into $\SC_s(\R^d)$, we have $\SC_s(\R^d)'\subset\SC(\R^d)'$. In particular, distributions in $\SC_s(\R^d)'$ are more regular than those from $\SC(\R^d)'$.
	
	%In order to define the fractional Laplacian $(-\Delta)^s T$ of a tempered distribution $T$ according to \eqref{eq:DualityFracLaplace} and for $(-\Delta)^s T$ to be a tempered distribution as well, we will then require $T\in \SC_s(\R^d)'$.
	%The right domain of $(-\Delta)^s$ is therefore $\SC_s(\R^d)'$, the space of linear maps $\SC_s(\R^d)\to\C$ that are continuous with respect to the topology induced by $(p_\alpha)_{\alpha\in\N^d}$ on $\SC_s(\R^d)$.
	
	\begin{Defin}\label{DefinFracLaplace}
		Fix $s\in (0,\infty)$. For any $T\in \SC_s(\R^d)'$, we let $(-\Delta)^s T : \SC(\R^d)\to \C$ be the tempered distribution defined by letting
		$$
		\ps{(-\Delta)^s T}{\psi}_\braS
		:=
		\ps{T}{(-\Delta)^s\psi}_{\SC_s', \SC_s}
		,\spa 
		\forall\ \psi\in\SC(\R^d)
		.
		$$ 
		%		The map $(-\Delta)^s T$ is a tempered distribution on $\R^d$.
	\end{Defin}
	
	Note that when $s=n\in\N$, Definition \ref{DefinFracLaplace} coincides with the usual definition of $(-\Delta)^n$ and simply extends its domain.
	For any $k\geq 1$, we let $\SC^k_s(\R^d)$ denote the \emph{cartesian} product of $\SC_s(\R^d)$, \ie the collection of vector fields $\Psi=(\psi_1,\ldots, \psi_k)$ such that $\psi_i\in\SC_s(\R^d)$ for each $i=1,2,\ldots, k$. Similarly to tempered distributions, we let $\SC^k_s(\R^d)'$ stand for the \emph{tensor} product of $\SC_s(\R^d)'$ acting multilinearly on $\SC^k_s(\R^d)$; see Definition \ref{DefinTemperedDistributions} of the Appendix and the comments below. The class $\SC_s(\R^d)$ is obviously closed under differentiation but not under Fourier transform, so that the same holds for $\SC_s(\R^d)'$. In addition, because
	$$
	\partial^\alpha (-\Delta)^s u 
	= 
	(-\Delta)^s \partial^\alpha u
	,\spa \forall\ \alpha\in\N^d,\ \forall\ u\in\SC(\R^d),
	$$ 
	we have
	$$
	\partial^\alpha (-\Delta)^s T
	= 
	(-\Delta)^s \partial^\alpha T
	,\spa \forall\ \alpha\in\N^d,\ \forall\ T\in\SC_s(\R^d)'
	.
	$$
	%for any $T\in\SC_s(\R^d)'$ and $\alpha\in\N^d$.\\
	
	Let us now give examples of tempered distributions that also belong to $\SC_s(\R^d)'$. It is straightforward to see that any (non-negative) measure $\mu$ belongs to $\SC_s(\R^d)'$ as soon as 
	$$
	\int_{\R^d} \frac{1}{1+\|x\|^{d+2(s-\lfloor s\rfloor)}}\, d\mu(x)
	<
	\infty 
	.
	$$
	In particular, any Borel probability measure on $\R^d$ belongs to $\SC_s(\R^d)'$. In addition, it is not difficult to see that 
	\begin{equation}\label{eq:L1locFrac}
		L^1_\loc(\R^d)\cap \SC_s(\R^d)'
		=
		\bigg\{ 
		u:\R^d\to \C\ \text{measurable},
		\int_{\R^d} \frac{|u(x)|}{1+\|x\|^{d+2(s-\lfloor s\rfloor)}}\, dx 
		< 
		\infty
		\bigg\}
		.
	\end{equation}
	In particular, we have $L^p(\R^d)\subset \SC_s(\R^d)'$ for any $p\in [1,+\infty]$.
	
	The next proposition provides different ways to compute $(-\Delta)^s u$ depending on the regularity of $u$. Among others, we consider the situation where $u$ belongs to a Sobolev space $H^{2s}$; see Section \ref{sec:Sobolev} of the Appendix for a definition of that space. Once again, we will tacitly write the argument of $u$ as $x$ and the argument of $\FT u$ as $\xi$. For the sake of completness, we prove the next result in Section \ref{Appendix:Laplacian} of the Appendix.
	% in terms of $\FT^{-1}(\|\xi\|^{2s}\FT u)$.
	
	%\begin{Prop}\label{PropLaplaceH2}
	%	Let $s\in (0,1)$ and $u\in H^{2s}(\R^d)$. Then, $(-\Delta)^s u \in L^2(\R^d)$ and we have
	%	$$
	%	(-\Delta)^s u 
	%	=
	%	(2\pi)^{2s} \FT^{-1}(\|\xi\|^{2s}\FT u)
	%	$$
	%	in $L^2(\R^d)$.
	%\end{Prop}

	\begin{Prop}\label{PropLaplaceExplicit}
		Fix $s>0$ and $u\in \SC_s(\R^d)'$.	
		\begin{enumerate}
			\item[(i)] If $s\in (0,1)$ and either $u\in \CC^{0,2s+\ve}(\Om)$ or $u\in\CC^{1,2s+\ve - 1}(\Om)$ for some $\ve>0$ and an open subset $\Om\subset\R^d$, then $(-\Delta)^s u$ is continuous over $\Om$ and we have
			$$
			((-\Delta)^s u)(x)
			=
			c_{d,s} \lim_{\eta\downarrow 0} \int_{\R^d\sm \BB^d_\eta(x)} \frac{u(x)-u(z)}{\|x-z\|^{d+2s}}\, dz
			,\spa \forall\ x\in\Om
			.
			$$
			
			\item[(ii)] If $u\in H^{2s}(\R^d)$, then $(-\Delta)^s u \in L^2(\R^d)$ and the equality
			$$
			(-\Delta)^s u 
			=
			(2\pi)^{2s} \FT^{-1}\Big(\|\xi\|^{2s}(\FT u)(\xi)\Big)
			$$
			holds in $L^2(\R^d)$, \ie almost everywhere.
			
			\item[(iii)] If $u\in L^1_\loc(\R^d)\cap \SC_s(\R^d)'$ and $\FF u\in L^1_\loc(\R^d)\cap \SC_s(\R^d)'$, then $\|\xi\|^{2s} \FF u\in \SC(\R^d)'$ and
			%			\FT u\in L^1_\loc(\R^d)$ and $\|\xi\|^{2s} (\FT u)(\xi)\in\SC(\R^d)'$, 
			we have
			$$
			(-\Delta)^s u
			=
			(2\pi)^{2s} \FT^{-1}\Big(\|\xi\|^{2s}(\FT u)(\xi)\Big)
			$$
			in $\SC(\R^d)'$. If, in addition, we have $\|\xi\|^{2s} (\FF u)(\xi)\in L^1(\R^d)$, then $(-\Delta)^s u\in \CC_0(\R^d)$ and the previous equality holds pointwise on $\R^d$.
			%		for any $x\in\R^d$ ;
			
		\end{enumerate}
	\end{Prop}

	Let us stress the fact that Proposition \ref{PropLaplaceExplicit} allows $s$ to be an integer. Lemma \ref{LemFourierTransformRiesz} of the Appendix combined with Proposition \ref{PropLaplaceExplicit} allows, in particular, to compute the fractional Laplacian of the map $1/\|x\|^{\alpha}$, for $\alpha\in (0,d-2s]$, which will be of interest in the next section.

	\begin{Corol}\label{CorolFTRiesz}
		Let $d\geq 2$ and $s\in (0,d/2)$. Then, for all $\alpha\in (0,d-2s]$ we have
		$$
		(-\Delta)^s\Big(\frac{1}{\|x\|^\alpha}\Big)	
		=
		\begin{cases}
			4^s\ 
			\frac{
				\Gamma(\frac{d-\alpha}{2}) \Gamma(\frac{\alpha+2s}{2})
			}{
				\Gamma(\frac{\alpha}{2})\Gamma(\frac{d-\alpha-2s}{2})
			}\
			\frac{1}{\|x\|^{\alpha+2s}} & \quad \text{if} \quad \alpha<d-2s,\\[3mm]
			4^{s}\pi^{d/2} 
			\frac{\Gamma(s)}{\Gamma(\frac{d-2s}{2})}\ \delta & \quad \text{if}\quad \alpha=d-2s,
		\end{cases}
		$$
		where $\delta$ stands for Dirac's distribution at $0$.
		%		, and
		%		$$
		%		C
		%		=
		%		\begin{cases}
			%			4^s\ 
			%			\frac{
				%			\Gamma(\frac{d-\alpha}{2}) \Gamma(\frac{\alpha+2s}{2})
				%			}{
				%				\Gamma(\frac{\alpha}{2})\Gamma(\frac{d-\alpha-2s}{2})
				%			} 
			%			& \quad \text{if} \quad \alpha<d-2s,\\[5mm]
			%			4^{s}\pi^{d/2} 
			%			\frac{\Gamma(s)}{\Gamma(\frac{d-2s}{2})} & \quad \text{if}\quad \alpha=d-2s.
			%		\end{cases}
		%		$$
	\end{Corol}
	
	In terms of regularity, we see from this example that $(-\Delta)^s$ acts by reducing the order of $\|x\|^{-\alpha}$ by $2s$, as is the case when $s$ is an integer. This heuristics, beyond this particular example, is formally stated in Proposition 2.7 of \cite{Silvestre2007}. We further refer to Section 3.4 of \cite{Kwa2019} for other explicit computations of fractional Laplacians.
	
	It is well-know that $1/\|x\|^{d-2}$ is the fundamental solution of $-\Delta$ when $d\geq 3$, \ie we have $(-\Delta) (1/\|x\|^{d-2}) = \delta$ up to a dimensional constant. Corollary \ref{CorolFTRiesz} entails that $1/\|x\|^{d-2s}$ is the fundamental solution of $(-\Delta)^s$ for any $s\in (0,d/2)$ as soon as $d\geq 2$. Notice that the value $s=1$ is indeed excluded when $d=2$ since $d-2s$ needs to be positive.

	\section{General PDE reconstruction}\label{sec:PDEGeneral}
	
	While the computation of the geometric cdf $\Fpg$ of a probability measure $P$ over $\R^d$ is straightforward (see Definition \ref{DefinRank}) and given by an expectation, it is not easy to understand what information about $P$ can be extracted from $\Fpg$. In this section, we address the problem of recovering $P$ from $\Fpg$. Our approach to tackle this question is to treat it as an \emph{inverse potential problem}. Indeed, our starting observation is that $\Fpg$ writes as the convolution
	$$
	\Fpg 
	=
	K_d*P
	$$
	between $P$ and a fixed kernel $K_d:\R^d\to\R^d$, defined as
	$$
	K_d(x)
	=
	\frac{x}{\|x\|} \I[x\neq 0]
	.
	%	\begin{cases}
		%		\frac{x}{\|x\|} & \textrm{ if } x\neq 0,\\
		%		0 & \textrm{ if } x = 0.
		%	\end{cases}
	$$	
	When no confusion about the dimension is possible we will omit the dependence of $K_d$ on $d$ and simply write $K$. The heuristics to solve the problem then consists in taking the Fourier transform $\FT(\Fpg)$ of $\Fpg$ which, at least formally, gives
	\vspace{-1mm}
	\begin{equation}\label{eq:HeuristicsEqFT}
		\FT(\Fpg) 
		=
		\FT(K_d) \FT(P)
		.
	\end{equation}
	This strategy was already used in \cite{Kol1997} to prove that if $P$ and $Q$ are probability measures on $\R^d$ such that $\Fpg(x)=F_Q^\g(x)$ for all $x\in\R^d$, then necessarily $P=Q$; see Theorem 2.5 and Corollary 2.9 in \cite{Kol1997}. For this purpose, they consider the map $h_P$, introduced in the present paper in Definition \ref{DefinObjectiveFunction}, and notice that $h_P$ is the convolution between $P$ and the kernel $x\mapsto \|x\|$. They provide the (distributional) Fourier transform of $x\mapsto \|x\|$ and, from their Theorem 2.5, deduce in Corollary 2.9 that the equality $\Fpg = F_Q^\g$ necessarily implies that $P=Q$. Through the same Fourier transform approach, the aim of the present section is to go a step further by establishing the explicit way of reconstructing an \emph{arbitrary} $P$ from $\Fpg$: we will show that the solution to this inverse problem takes the form $P=\LL_d(\Fpg)$, where $\LL_d$ is the operator
	$$
	\LL_d
	:=
	\gamma_d\ (-\Delta)^{\frac{d-1}{2}}(\nabla\ \cdot\ )
	.
	$$
	The operator $\LL_d$ involves the divergence $\nabla\ \cdot$ and the Laplacian~$(-\Delta)^{(d-1)/2}$ defined in Section \ref{sec:IntroLaplace}, and a positive constant $\gamma_d$ given by $\gamma_d^{-1}
	=
	2^d \pi^{(d-1)/2} \Gamma((d+1)/2)
	.
	$
	%	The operator $\LL_d$ then writes
	%	$$
	%	\LL_d
	%	:=
	%	\gamma_d\ (-\Delta)^{\frac{d-1}{2}}(\nabla\ \cdot\ )
	%	.
	%	$$
	The dual operator $\LL_d^*$ of $\LL_d$ is given by
	$$ 
	\LL_d^*
	:=
	\gamma_d\ \nabla (-\Delta)^{\frac{d-1}{2}}
	.
	$$
	While the precise definitions of the domains of $\LL_d$ and $\LL_d^*$ are given in Section \ref{Appendix:PDEGeneralOperators} of the Appendix, let us mention that $(L^\infty(\R^d))^d$ is always contained in the domain of $\LL_d$, so that $\LL_d(\Fpg)$ is defined for an \emph{arbitrary} probability measure $P$ on $\R^d$.
	%	Recalling the definition of $\SC_{1/2}(\R^d)'$ and $\SC_{1/2}^d(\R^d)'$, introduced in Section \ref{sec:IntroLaplace}, we see that the domain of $\LL_d$ always contains $\SC_{1/2}^d(\R^d)'$ and that the domain of $\LL_d^*$ always contains $\SC_{1/2}(\R^d)'$; the exact definition of their domains depends on $d$ and is discussed in Appendix \ref{Appendix:PDEGeneralOperators}. In particular, since $\Fpg\in (L^\infty(\R^d))^d\subset\SC_{1/2}^d(\R^d)$ for any probability measure $P$ on $\R^d$, we see that $\Fpg$ always belongs to the domain of $\LL_d$, so that the equality $P=\LL_d(\Fpg)$ makes sense for an \emph{arbitrary} probability measure $P$ on $\R^d$.
	We say that $\LL_d^*$ is dual to $\LL_d$ because the equality 
	$$
	\ps{\LL_d T}{\psi}_{\SC', \SC}
	=
	\ps{T}{\LL_d^* \psi}_{\SC', \SC}
	%	,\spa 
	%	\forall\ \psi\in\SC(\R^d)
	,
	$$
	holds for all $\psi\in\SC(\R^d)$ and $T\in\SC_{1/2}^d(\R^d)'$---see Section \ref{sec:IntroLaplace} for the definition of $\SC^d_{1/2}(\R^d)'$.
	%	which will provide the solution to (\ref{eq:HeuristicsEqFT}). 
	%	Notice that $\LL_d$ acts on vector fields and $\LL_d^*$ on real- or complex-valued maps. Defining these operators properly on the approriate space of distributions requires some care and depends on the dimension $d$. We discuss this in more details in Appendix \ref{Appendix:PDEGeneralOperators}. The domain of $\LL_d^*$ contains $\SC_{1/2}(\R^d)'$, and the domain of $\LL_d$ contains the space $\SC_{1/2}^d(\R^d)'$ of $\C^d$-valued distributions in $\SC_{1/2}(\R^d)'$, \ie ``vector fields'' in the sense of distributions; see Section \ref{sec:IntroLaplace} and Appendix \ref{sec:Distributions}. 

	%	Before stating the first theorem of this section, let us emphasize that, in $\R$, $\LL_1$ and $\LL_1^*$ reduce to 
	%	$
	%	\LL_1
	%	=
	%	\frac{1}{2}\frac{d}{dx}
	%	=
	%	\LL_1^*
	%	$
	%	.\\
	
	%We will show that the Fourier transform $\FT(K_d)$ of $K_d$ is given, up to a multiplicative constant $C_d$ that only depends on $d$, by 
	%$$
	%\pv\Big(\frac{\xi}{\|\xi\|^{d+1}}\Big)
	%,
	%$$
	%where $\pv$ stands for \emph{principal value}. In other words, we have
	%$$
	%\int_{\R^d} K_d(x) \wwidehat{\psi}(x)\, dx 
	%=
	%C_d \lim_{\eta\to 0} \int_{\R^d\sm B_\eta} \frac{\xi}{\|\xi\|^{d+1}} \psi(\xi)\, d\xi 
	%$$
	%for any $\psi\in\SC(\R^d)$, where $B_\eta$ stands for the ball of radius $\eta>0$ centered at the origin, and $\SC(\R^d)$ (or, equivalently, $\SC(\R^d,\C)$) stands for the complex-valued Schwartz class over $\R^d$.\\
	
	%	\subsection{Distributional recovery}\label{sec:DistributionalRecovery}
	
	To make the heuristics (\ref{eq:HeuristicsEqFT}) rigorous, one has to make sure that the product $\F(K_d)\F(P)$ is well-defined. Because $K_d$ is bounded, it is a tempered distribution, hence so is $\F(K_d)$. Proposition \ref{PropConvolution} of the Appendix then entails that $\F(K_d)\F(P)$ is a tempered distribution provided $P$ belongs to $\SC(\R^d)$, \ie $P$ admits a density in the Schwartz class. We thus start our investigation in Theorem \ref{TheorEDPSchwartz} by assuming that $P$ has such a density before extending the result in Theorem \ref{TheorEDPDistributions}.
	%	
	%	In Theorem \ref{TheorEDPSchwartz}, we assume that $P$ admits a well-behaved density and explicitly recover the density from $\Fpg$. Then, we extend this result to an arbitrary probability measure $P$ in Theorem \ref{TheorEDPDistributions}, and recover $P$ from $\Fpg$ in the sense of distributions. 
	
	\begin{Theor}
		\label{TheorEDPSchwartz}
		Let $d\geq 1$ and $P$ be a Borel probability measure on $\R^d$. Assume that $P$ admits a density $f_P\in\SC(\R^d)$ in the Schwartz class with respect to the Lebesgue measure. Then $\Fpg\in\CC^\infty(\R^d)$ and we have
		$
		f_P(x)
		=
		(\LL_d \Fpg)(x)
		$
		for all $x\in\R^d$.
		%	(see Definition \ref{DefinDiffOperator} for the definition of $\LL_d$ and $D(\LL_d)$).
		%	 where $\LL_d$ is the differential operator introduced in Definition \ref{DefinDiffOperator}.
	\end{Theor}
	
	The proof of Theorem \ref{TheorEDPSchwartz} is constructive in the sense that it \emph{finds out} the operator $\LL_d$ by which $P$ is to be recovered from $\Fpg$. Once one believes that $\LL_d(\Fpg)$ is the correct way to reconstruct $P$, a direct proof relying on PDE arguments can be given as follows. Proposition \ref{PropRankIntermediateRegularity} below entails that
	$$
	\nabla\cdot \Fpg 
	=
	(\nabla\cdot K)*f_P
	,
	$$
	where $(\nabla\cdot K)(x)=(d-1)/\|x\|$ for all $x\in\R^d\sm\{0\}$, and that $\nabla\cdot \Fpg$ is bounded on $\R^d$. As we already observed after Corollary \ref{CorolFTRiesz}, the map $1/\|x\|$ is, up to a constant, the fundamental solution of $(-\Delta)^{(d-1)/2}$, which will provide the conclusion. In what follows, we justify the use of Corollary \ref{CorolFTRiesz} and show that $\gamma_d$ is the right constant. Since $\nabla\cdot \Fpg$ is bounded, it belongs to $L^1_\loc(\R^d)$ and $\SC_{(d-1)/2}(\R^d)'$. Also observe that Proposition \ref{PropConvolution} of the Appendix entails that
	$$
	\FF(\nabla\cdot \Fpg)
	=
	(d-1)\ \FF\Big(\frac{1}{\|x\|}\Big)\ \FF(f_P)
	,
	$$
	which, given the explicit form of $\FF(1/\|x\|)$ provided in Lemma \ref{LemFourierTransformRiesz} of the Appendix, entails that $\FF(\nabla\cdot \Fpg)$ belongs to $L^1_\loc(\R^d)$ and $\SC_{(d-1)/2}(\R^d)'$ as well. Consequently, 	Proposition \ref{PropLaplaceExplicit} (iii) implies that 
	$$
	(-\Delta)^{\frac{d-1}{2}}(\nabla\cdot \Fpg) 
	=
	(d-1)\ \Big((-\Delta)^{\frac{d-1}{2}} \frac{1}{\|x\|} \Big)*f_P
	,
	$$
	which, by Corollary \ref{CorolFTRiesz} with $s=(d-1)/2$ and $\alpha=1$, yields 
	$$
	(-\Delta)^{\frac{d-1}{2}}(\nabla\cdot \Fpg) 
	=
	(d-1)\ 4^{\frac{d-1}{2}} \pi^{\frac{d}{2}} \frac{\Gamma(\frac{d-1}{2})}{\Gamma(\frac{1}{2})}\
	(\delta*f_P)
	=
	2^d \pi^{\frac{d-1}{2}} \Gamma\Big(\frac{d+1}{2}\Big)\ f_P 
	=
	\gamma_d^{-1}\ f_P
	.
	$$

	As we mentioned before Theorem \ref{TheorEDPSchwartz}, in general the heuristics (\ref{eq:HeuristicsEqFT}) only works provided that $P$ has a density in the Schwartz class. It cannot be extended to general probability measures; notice that the straightforward approach presented right after Theorem \ref{TheorEDPSchwartz} cannot be largely extended in general either. Nonetheless, Theorem \ref{TheorEDPSchwartz} is enough to reconstruct an \emph{arbitrary} probability measure from its geometric cdf in the sense of distributions by a density argument. This is the content of the next theorem.
	%Because $\LL_1=\frac{1}{2}\frac{d}{dx}$ and $\Fpg=2F_P-1$ over $\R$ when $d=1$, where $F_P$ is the usual univariate cdf of $P$, the equality
	%$
	%\frac{d}{dx} F_P
	%=
	%f_P
	%$
	%is a particular case of Theorem \ref{TheorEDPSchwartz} when $f_P$ belongs to the Schwartz class on $\R$.

	%
	%\vspace{3mm}
	%
	%We are now able to give the proof of Theorem \ref{TheorEDPSchwartz}.
	
	%	\vspace{3mm}

	\begin{Theor}\label{TheorEDPDistributions}
		Let $d\geq 1$ and $P$ be a Borel probability measure on $\R^d$. The equality 
		$
		P
		=
		\LL_d(\Fpg)
		$
		holds in $\SC(\R^d)'$, \emph{i.e.} we have $\ps{P}{\psi}_{\SC', \SC}=\ps{\LL_d (\Fpg)}{\psi}_{\SC',\SC}$ for any $\psi\in \SC(\R^d)$. In other words, letting $Z$ be a random vector with law $P$, we have
		$$
		%		\int_{\R^d} \psi(x)\, dP(x)
		\E\big[\psi(Z)\big]
		=
		\int_{\R^d} \ps{\Fpg(z)}{(\LL_d^* \psi)(z)}\, dz
		,\spa \forall\ \psi\in\SC(\R^d)
		.
		$$
		%	for any $\psi\in\SC(\R^d)$.
	\end{Theor}

	%	\subsection{Pointwise recovery}\label{sec:PointwiseRecovery}
	
	%	In Section \ref{sec:DistributionalRecovery}, 
	%	We established that any probability measure $P$ on $\R^d$ writes $P=\LL_d(\Fpg)$ in the sense of tempered distributions, where 
	%	$$
	%	\LL_d
	%	=
	%	\gamma_d\ (-\Delta)^{\frac{d-1}{2}}\nabla\ \cdot
	%	$$
	%	is the operator introduced at the beggining of Section \ref{sec:PDEGeneral}. We further showed that if $P$ admits a smooth and fast-decreasing density $f_P\in \SC(\R^d)$, the equality $f_P=\LL_d (\Fpg)$ actually holds pointwise over $\R^d$. 
	
	%	In this section, 
	
	\section{Regularity and pointwise reconstruction}\label{sec:Regularity}
	
	Let us now investigate conditions less restrictive than Theorem \ref{TheorEDPSchwartz} that will ensure that, when $P$ admits a density $f_P$, the equality $f_P=\LL_d(\Fpg)$ holds pointwise; see also Section \ref{sec:Examples} of the Appendix for examples where we explicitly compute $\LL_d(\Fpg)$ and show that this indeed coincides with $f_P$. When this is the case, one can compute the value of $(\LL_d \Fpg)(x)$ by successively applying the differential operators involved in the definition of $\LL_d$ to $\Fpg$ pointwise at $x$. In the univariate case $d=1$, minimal assumptions are already well-known. Indeed, let us observe that $\LL_1=2^{-1}d/dx$ and $\Fpg=2F_P-1$, where 
	$
	F_P
	%	(x)
	%	=
	%	\int_{-\infty}^x f_P(s)\, ds 
	$
	is the usual cdf of $P$; see the comments after Definition \ref{DefinRank}. Provided $f_P$ is continuous in a neighbourhood of $x_0\in\R$, the fundamental theorem of calculus then yields 
	$$
	f_P(x_0)=\frac{d F_P}{dx}(x_0)= (\LL_1 \Fpg)(x_0)
	.
	$$
	In dimension $d\geq 2$, since computing $\LL_d(\Fpg)$ pointwise requires taking $d$ derivatives, we expect $\Fpg$ to be of class $\CC^d$; in fact, we need to compute $d$ derivatives when $d$ is odd, but $d-1$ derivatives and one pseudo-derivative when $d$ is even, which, on account of Proposition \ref{PropLaplaceExplicit} (i) can be performed provided $\Fpg\in\CC^d$. Since differentiating $\Fpg=K*f_P$ amounts to differentiating $K$, and $\partial^\alpha K$ behaves like $1/\|x\|^{|\alpha|}$, then the derivatives $\partial^\alpha \Fpg$ are well-defined on $\R^d$ under mild assumptions as long as $|\alpha|<d$. In particular, we show in Proposition \ref{PropRankIntermediateRegularity} that $\Fpg$ is of class $\CC^{d-1}$.	
	
	As a preliminary result, let us mention that a straightforward application of Lebesgue's dominated convergence theorem yields that $\Fpg$ is continuous at a point $x\in\R^d$ if and only if $P[\{x\}]=0$. In particular, $\Fpg$ is continuous over $\R^d$ when $P$ admits a density with respect to the Lebesgue measure. We now turn to the differentiability of $\Fpg$.
	
	%	In this section, we will always let $Z$ denote a random $d$-vector with law $P$.
	
	\begin{Prop} \label{PropRankIntermediateRegularity}
		Fix $d\geq 2$, an open subset $\Om\subset\R^d$, and an integer $\ell\in [1,d-1]$. Let $P$ be a Borel probability measure on $\R^d$, and assume that $P$ admits a density $f_\Om\in L^1(\Om)$ over $\Om$ with respect to the Lebesgue measure. (i)
		%		\begin{enumerate}
			%			\item 
			If $f_\Om\in L^p_\loc(\Om)$ for some $p\in (d/(d-\ell),\infty]$, then $\Fpg\in\CC^{\ell}(\Om)$ and
			$$
			(\partial^\alpha \Fpg)(x)
			=
			\E\big[(\partial^\alpha K)(x-Z)\big]
			,\spa \forall\ x\in \Om 
			,
			%			\int_{\R^d} (\partial^\alpha K)(x-z)\, dP(z)
			$$
			for all $\alpha\in\N^d$ with $|\alpha|\leq \ell$, where $Z$ is a random vector with law $P$. (ii)
			%			\item 
			If $\Om=\R^d$ and the density 
			%		$f_P:=f_{\R^d}$ 
			belongs to $L^p(\R^d)$ for some $p\in (d/(d-\ell),\infty]$, then $\partial^\alpha \Fpg$ converges to $0$ at infinity for any $\alpha\in\N^d$ such that $1\leq |\alpha|\leq \ell$.
			%		\item If $\Om=\R^d$ and the density is $\beta$-H{\"o}lder continuous over $\R^d$ for some $\beta\in (0,1]$, then $\partial^\alpha \Fpg$ is $(\frac{|\alpha|}{d} \beta)$-H{\"o}lder continuous over $\R^d$ for any $\alpha\in\N^d$ with $|\alpha|\leq \ell$.
			%		\end{enumerate}  
	\end{Prop}
	
	Let us stress that Part (i) of Proposition \ref{PropRankIntermediateRegularity} requires $P$ to have a density only on $\Om$. In particular, the result holds even if $P$ admits atoms outside $\Om$. Also note that if $\Om$ lies outside the support of $P$, in the sense that $P[\Om]=0$, then the result also applies since $P$ then admits the density $f_\Om\equiv 0$ on $\Om$. 
	%	When $\Om=\R^d$, one can see from the proof that the condition $f_\Om\in L^p(\R^d)$ in Part (ii) can be relaxed into $f_\Om\in L^p(\R^d\sm K)$ for some compact set $K\subset \R^d$.
	
	%		Let us make a comment on the requirement $p>d/(d-\ell)$.
	%%		The requirement that the density $f_\Om$ be in $L^p_\loc$ for some $p>\frac{d}{d-\ell}$ might feel unnatural in the first part of Proposition \ref{PropRankIntermediateRegularity}. 
	%%		If $f_\Om$ is bounded in a neighbourhood of $x$, the result is a straightforward application of Lebesgue's dominated convergence theorem. 
	%		This integrability condition allows $f_\Om$ to have singularities but prevents them from growing too fast. For instance, if $f_\Om(z)$ behaves like $\|z\|^{-\beta}$ around $z=0$ for some $\beta>0$, then $f_\Om$ will be integrable to the $p$th power around $0$ provided $\beta < d/p$. Because $p$ can be taken arbitrarily close to $d/(d-\ell)$, this is equivalent to $\beta < d-\ell$ without loss of generality. The latter condition has a clear interpretation: the more derivatives we want, \ie the higher we take $\ell$, the more derivatives of $K$ must be integrable, the smaller the spikes of $f_\Om$ should be.

	When $d\geq 2$, reaching a regularity up to order $d$ is more challenging than in the univariate case $d=1$; when $|\alpha|=d$, then $\partial^\alpha K$ scales like $1/\|x\|^d$ and is not integrable near the origin. Lemma \ref{LemDerivOK} of the Appendix entails that $(\partial^\alpha \Fpg)(x)$
	%	 the convolution $(\partial^\alpha K_d)*f_P$ 
	is \emph{formally} given by the singular integral
	\begin{equation}\label{eq:SingInt}
		\lim_{\eta\downarrow 0} 
		\int_{\R^d\sm \BB^d_\eta(x)}
		%	\E\Big[(\partial^\alpha K)(x-Z)\ \I\big[\|Z-x\|>\eta\big]\Big]
		(\partial^\alpha K)(x-z)\  f_P(z)\, dz
		+
		\kappa_\alpha f_P(x)
		,
	\end{equation}
	for some constant $\kappa_\alpha\in\R$. Still with $|\alpha|=d$, one easily proves, by induction for instance, that there exists a polynomial in $d$ variables $p_\alpha:\R^d\to\R$ such that $(\partial^\alpha K)(z)=\|z\|^{-d} p_\alpha(z/\|z\|)$ for all $z\neq 0$. It follows from Lemma \ref{GreenFormula} of the Appendix that the integral of $z\mapsto (\partial^\alpha K)(x-z)$ over $\BB^d_1(x)\sm \BB^d_\eta(x)$ vanishes. In particular, when $f_P$ is $\beta$-H{\"o}lder continuous for some $\beta\in (0,1]$, we can write the singular integral in (\ref{eq:SingInt}) as the pointwise well-defined map
	\begin{equation}\label{eq:SingIntBis}
		\int_{\R^d\sm \BB^d_1(x)}
		%	\E\Big[(\partial^\alpha K)(x-Z)\ \I\big[\|Z-x\|>\eta\big]\Big]
		(\partial^\alpha K)(x-z)\  f_P(z)\, dz
		+
		\int_{\BB^d_1(x)}
		%	\E\Big[(\partial^\alpha K)(x-Z)\ \I\big[\|Z-x\|>\eta\big]\Big]
		(\partial^\alpha K)(x-z)\  (f_P(z)-f_P(x))\, dz
		.
	\end{equation}
	The representation (\ref{eq:SingIntBis}) remains valid at $x$ if we only require that there exist positive constants $C_x$, $\delta_x$, and $\beta_x\in (0,1]$ such that $|f_P(y)-f_P(x)|\leq C_x \|y-x\|^{\beta_x}$ for all $y\in\R^d$ with $\|y-x\|<\delta_x$, if one replaces $\BB^d_1(x)$ by $\BB^d_{\delta_x}(x)$ in (\ref{eq:SingIntBis}). In fact, writing the second term of (\ref{eq:SingIntBis}) in spherical coordinates centered at $x$ it is easy to see that this expression is well-defined at $x$ as soon as the modulus of continuity $r\mapsto \omega_x(r):= \sup\{|f_P(y)-f_P(x)| : \|y-x\|\leq r\}$ of $f_P$ at $x$ is integrable on $(0,1)$ with respect to the measure $(1/r)\times dr$, known as a Dini-type condition. Furthermore, when $f_P$ does not satisfy the Dini condition, the formal expression (\ref{eq:SingInt}) remains almost everywhere well-defined. Indeed, using the fact that $(\FF K)(\xi)\propto \PV(\xi/\|\xi\|^{d+1})$, established in the proof of Theorem \ref{TheorEDPSchwartz}, we see that $(\FF(\partial^\alpha K))(\xi) \propto\xi^\alpha (\partial^\alpha K)(\xi)$ is a bounded map on $\R^d$ when $|\alpha|=d$.
	%	; let $B>0$ be such that $|\FF(\partial^\alpha K)|\leq B$ on $\R^d$.
	Therefore, Plancherel's identity and Proposition \ref{PropConvolution} of the Appendix entail that 
	\begin{equation}\label{eq:BoundL2}
		%	\|\partial^\alpha \Fpg\|_{L^2(\R^d)}
		%	=
		%	\|\partial^\alpha (K*g)\|_{L^2(\R^d)}
		%	=
		\|(\partial^\alpha K)*g\|_{L^2(\R^d)}
		%	=
		%	\|\FF(\partial^\alpha K)\FF(g)\|_{L^2(\R^d)}
		\lesssim
		%	\|\FF(g)\|_{L^2(\R^d)}
		%	=
		\|g\|_{L^2(\R^d)}
		,\spa 
		\forall\ g\in\SC(\R^d)
		;
	\end{equation}
	throughout, the missing constants do not depend on $g$. It follows from (\ref{eq:BoundL2}) and the representation of $\partial^\alpha (K*g)$ from (\ref{eq:SingInt}) that $\|\partial^\alpha (K*g)\|_{L^2(\R^d)}\lesssim \|g\|_{L^2(\R^d)}$ for all $g\in\SC(\R^d)$; this inequality thus extends to any $g\in L^2(\R^d)$ by density of $\SC(\R^d)$ in $L^2(\R^d)$. Theorem 5.3.3 in \cite{Gra2014} then entails that for any $p\in (1,\infty)$ we have $\|\partial^\alpha (K*g)\|_{L^p(\R^d)}\lesssim \|g\|_{L^p(\R^d)}$ for all $g\in L^p(\R^d)$. In addition, denoting by $\lambda$ the Lebesgue measure in $\R^d$, we have
	$$
	\lambda\Big(\big\{x\in\R^d : |\partial^\alpha (K*g) (x)|>t\big\}\Big)
	\lesssim 
	\frac{1}{t}
	,\spa 
	\forall\ t>0,\ 
	\forall\ g\in L^1(\R^d)
	.
	$$
	Consequently, we see that the distributional derivatives $\partial^\alpha \Fpg$ of order $d$ of $\Fpg$ are, in fact, almost everywhere well-defined maps as soon as $P$ admits a density $f_P\in L^1(\R^d)$, and that these derivatives further belong to $L^p(\R^d)$, $1<p<\infty$, when $f_P$ does. The map $\partial^\alpha \Fpg$ is, however, not continuous for an arbitrary continuous density $f_P$ without additional requirements as we show in the next Proposition. In fact, the quantity (\ref{eq:SingInt}) cannot even map any continuous density $f_P$ to a locally bounded function. 
	%	In particular, when a probability measure $P$ admits 
	%	
	%	We deduce that the linear map $f\mapsto \partial^\alpha (K*f)$ extends to a bounded operator on $L^2(\R^d)$ onto itself. In particular, there exists $c_\alpha > 0$ such that $\|\partial^\alpha(K*f)$
	%	This last integral, however, does not define a continuous function for an arbitrary continuous density $f_P$, which is in strong contrast with the univariate case. In fact, the quantity (\ref{eq:SingInt}) cannot even map any continuous density $f_P$ to a locally bounded function. 
	
	\begin{Prop}\label{PropContinuousCounterExample}
		Fix $d\geq 2$. For any open set $\Om\subset\R^d$ and $\alpha\in \N^d$ with $|\alpha|=d$, there exists a Borel probability measure $P$ on $\R^d$ with a continuous density $f_P\in\CC^0_c(\Om)$ compactly supported in $\Om$ such that the distributional derivative $\partial^\alpha \Fpg$ of $\Fpg$ is not a bounded map on $\Om$. In particular, we have $\Fpg\in \CC^{d-1}(\R^d)\sm \CC^d(\Om)$.
	\end{Prop}
	
	Let us stress again how Proposition \ref{PropContinuousCounterExample} contrasts with the univariate case $d=1$, where $\Fpg\in \CC^1(\Om)$ as soon as $P$ admits a continuous density on $\Om$. As we already briefly outlined, the negative result of Proposition \ref{PropContinuousCounterExample} can be overcome when $f_P\in \CC^{0,\beta}_\loc(\Om)$ for some $\beta\in (0,1]$. This is the content of the next Theorem.

	\begin{Theor}
		\label{TheorRankRegularityOdd}
		Fix $d\geq 1$. Let $\Om\subset\R^d$ be an open subset and $P$ a Borel probability measure on $\R^d$. Assume that $P$ admits a density $f_\Om\in L^1(\Om)$ over $\Om$ with respect to the Lebesgue measure. If $f_{\Om}\in \CC^{k,\beta}_\loc(\Om)$ for some $k\in\N$ and $\beta\in (0,1)$, then $\Fpg\in \CC^{d+k,\beta}_\loc(\Om)$ and
		\begin{equation}\label{eq:PointwiserRecovery}
			%		\label{EqFDimGreater2}
			f_{\Om}(x)
			=
			(\LL_d \Fpg)(x)
			,\spa \forall\ x\in\Om
			.
		\end{equation}
		%		(ii) If $f_P\in\CC^\infty(\Om)$, then $\Fpg\in\CC^\infty(\Om)$ and (\ref{eq:PointwiserRecovery}) holds. 
	\end{Theor}
	
	%	\vspace{3mm}
	When $\Om=\R^d$, we provide in Section \ref{sec:Examples} of the Appendix a few examples where we explicitly compute $(\LL_d \Fpg)(x)$ and show that this indeed coincides with $f_P(x)$. 
	
	We already showed that Theorem \ref{TheorRankRegularityOdd} fails for $\beta=0$. The value $\beta=1$ must also be discarded as the spaces $\CC^{k,1}$ are not well-suited to study ellitic regularity; see, e.g., the proof of Theorem 2.14 in \cite{RosOtonRegulairty} and the H{\"o}lder estimates preceeding Theorem 1.7 in the same reference to see why $\beta<1$ is an important assumption. This, however, is not a major issue: if $f_\Om\in \CC^{k,1}_\loc(\Om)$ then $f_\Om\in\CC^{k,1-\ve}_\loc(\Om)$ so that $\Fpg\in\CC^{d+k,1-\ve}_\loc(\Om)$ for all $\ve\in (0,1)$. Therefore, if $f_\Om\in \CC^k(\Om)$ then $f_\Om\in \CC^{k-1,1}_\loc(\Om)$ so that $\Fpg\in \CC^{d+k-1}(\Om)$. In particular, if $f_\Om\in \CC^\infty(\Om)$ then $\Fpg\in\CC^\infty(\Om)$. 
	
	When~$d$ is even, Proposition \ref{PropLaplaceExplicit} (i) and Theorem \ref{TheorRankRegularityOdd} entail that, with $c_{d, 1/2}$ be the constant defined in (\ref{Eq:ConstantCds}) and letting
	%	If $f_\Om\in C^1(\Om)$, then $f\in \CC^{0,\alpha}_\loc(\Om)$ for any $\alpha\in [0,1]$ so that $\Fpg\in \CC^{d}(\Om)$ with $f_\Om = \LL_d(\Fpg)$ pointwise in $\Om$. Similarly, if $f_\Om\in \CC^k(\Om)$ then $\Fpg\in \CC^{d+k-1}(\Om)$. In particular, if $f_\Om\in\CC^\infty(\Om)$ then $\Fpg\in \CC^\infty(\Om)$. 
	%	
	%		If $\Om=\R^d$ and $f:=f_\Om$ is smooth except on a closed set $S$, i.e. $f\in \CC^\infty(\R^d\sm S)$
	%	
	%	A few comments are in order. (i) If $\beta=1$ then $\CC^{d+k,1-\ve}(\Om)$ for any $\ve\in (0,1)$. (ii) 
	%	In addition, letting 
	$
	u
	:=
	\gamma_d\ (-\Delta)^{(d-2)/2} (\nabla\cdot \Fpg)
	$, we can compute $f_\Om(x)=(\LL_d \Fpg)(x)=((-\Delta)^{1/2}u)(x)$ pointwise as
	$$
	f_\Om(x)
	%		=
	%		((-\Delta)^{1/2}\Fpg^{(d-1)})(x)
	%=
	%\big((-\Delta)^{1/2}u\big)(x)
	=
	c_{d,1/2} 
	\lim_{\eta\downarrow 0} 
	\int_{\R^d\sm B_\eta(x)} \frac{u(x)-u(z)}{\|x-z\|^{d+1}}\, dz
	,\spa \forall\ x\in\Om
	.
	$$
	Let us observe that this implies that, for $d$ even, the reconstruction of $P$ from $\Fpg$ is nonlocal in the sense that the value of $f_\Om(x)$ depends on the values of $\Fpg$ on all of $\R^d$ through $(-\Delta)^{1/2} u$; this is not the case when $d$ is odd as $\LL_d(\Fpg)$ only involves local differential operators, so that the value of $f_\Om(x)$ depends only on the values of $\Fpg$ in a neighbourhood of $x$. Thus, it could be surprising that Theorem~\ref{TheorRankRegularityOdd} establishes the regularity of $\Fpg$ locally from that of the density of $P$. In fact, when we forget about the \emph{values} and restrict to \emph{regularity} only, the interaction between $u$ and $f_\Om=(-\Delta)^{1/2}u$ somehow becomes local; see the proofs of Proposition \ref{PropLaplaceExplicit} (i) and Proposition \ref{PropSchauderFrac} of the Appendix.

\section{Local characterisation}\label{sec:Localisation}

In Section \ref{sec:PDEGeneral} and Section \ref{sec:Regularity} we explained how $P$ and $\Fpg$ interact globally; $P$ is recovered from $\Fpg$ through the \emph{global} identity $P=(\LL_d \Fpg)$ where $\LL_d = \gamma_d\ (-\Delta)^{(d-1)/2}(\nabla\cdot\ )$. The present section is devoted to studying the \emph{local} interplay between $P$ and $\Fpg$: what information about $P$ can one recover through the knowledge of $\Fpg$ on a subset $\Om$ of $\R^d$? As we already noticed at the end of Section \ref{sec:Regularity}, we expect the answer to depend on whether $d$ is odd or even due to the presence of the \emph{nonlocal} fractional Laplacian $(-\Delta)^{1/2}$ in the operator $\LL_d$ when $d$ is even. When $d$ is odd, however, the operator $\LL_d$ is purely local so that the knowledge of $\Fpg$ on an open subset $\Om\subset\R^d$ immediately provides the corresponding information about $P$ through the restriction of the identity $P=\LL_d(\Fpg)$ to $\Om$. We have the following result.

\begin{Prop}\label{PropLocalOdd}
Fix an odd integer $d\geq 1$. Let $P$ and $Q$ be a Borel probability measures on $\R^d$. Let $\Om\subset\R^d$ be an open subset, and assume that $\Fpg(x)=F_Q^\g(x)$ for all $x\in \Om$. Then $P$ and $Q$ coincide over $\Om$, \ie $P(E)=Q(E)$ for all Borel subsets $E\subset\Om$.
\end{Prop}

% in section we investigate the local properties of the operator $\LL_d=(-\Delta)^{\frac{d-1}{2}}(\nabla\ \cdot\ )$. We already mentioned in Section \ref{sec:PDEGeneral} it displays substantially different behaviours in odd and even dimensions. This is due the nature of $(-\Delta)^{\frac{d-1}{2}}$, which depends on whether $\frac{d-1}{2}$ is an integer or not. When $\frac{d-1}{2}\in\N$, then $(-\Delta)^{\frac{d-1}{2}}$ is the classical differential operator that consists in applying the Laplacian $-\Delta$ successively $(d-1)/2$ times. This operator is local in nature: if smooth functions $u_1$ and $u_2$ coincide over an open subset $\Om\subset\R^d$, then $(-\Delta)^{\frac{d-1}{2}}u_1$ and $(-\Delta)^{\frac{d-1}{2}}u_2$ also coincide over $\Om$. When $d$ is even, then $\frac{d-1}{2}\in\R\sm\N$; in this case, we write 
%$$
%(-\Delta)^{\frac{d-1}{2}}
%=
%(-\Delta)^{\frac12}(-\Delta)^{\frac{d-2}{2}}
%.
%$$
%Although $(-\Delta)^{1/2}$ acts like a derivative in terms of regularity (see Proposition 2.6 in \cite{Silvestre2007}), it is a non-local operator: the value of $(-\Delta)^{1/2}u$ at a point $x\in\R^d$ depends on the values of $u$ over all of $\R^d$.

It has been long known that the geometric cdf characterises probability measures in arbitrary dimension $d$: if $P$ and $Q$ are Borel probability measures on $\R^d$ and if $\Fpg(x)=F_Q^\g(x)$ for all $x\in\R^d$, then $P=Q$; see Theorem 2.5 and Corollary 2.9 in \cite{Kol1997}. When $d$ is odd, Proposition \ref{PropLocalOdd} is thus a refinement that leverages the explicit reconstruction of $P$ established in Theorem \ref{TheorEDPDistributions}. When $d$ is even, however, the proof does not apply, but one could still expect or hope that the local characterisation of Proposition \ref{PropLocalOdd} still holds in some particular situations. The next result establishes that this is \emph{never} the case.

\begin{Prop}\label{PropLocalEven}
Fix an even integer $d\geq 2$. Let $P$ and $Q$ be Borel probability measures on $\R^d$, and $\Om\subset\R^d$ an open subset. If $P=Q$ on $\Om$ and $\Fpg=F_Q^\g$ on $\Om$, then $P=Q$ on $\R^d$.
\end{Prop}

At this point, one could legitimately wonder if Proposition \ref{PropLocalOdd} is in fact not an empty statement when $d>1$, in the sense that the conclusion of Proposition \ref{PropLocalEven} still holds when $d$ is odd, perhaps for different reasons. We will see, however, that this is not the case: when $d\geq 3$ is odd, there exist distinct probability measures $P$ and $Q$ on $\R^d$ and an open set $\Om\subset\R^d$ such that $\Fpg$ and $F_Q^\g$ coincide over $\Om$. This will be a consequence of the next Proposition. 

We say that a Borel probability measure $P$ on $\R^d$ is spherically symmetric if $P[O A]=P[A]$ for all Borel subset $A\subset \R^d$ and $d\times d$ orthogonal matrix $O$.

\begin{Prop}\label{ExplicitSpherOdd}
Fix an odd integer $d\geq 3$. There exist real coefficients $w_0, w_1\ldots, w_{(d-3)/2}$ and a continuous map $p:(0,1]\to\R$ such that for any Borel probability measure $P$ on $\R^d$ admitting a bounded density in a neighbourhood of the origin, and for $Z$ a random $d$-vector with law $P$, we have
$$
\Fpg(su)
=
\bigg( 
\sum_{i=0}^{\frac{d-3}{2}} \Big( \E\Big[\frac{1}{\|Z\|^{2i+1}}\Big] w_i s^{2i+1} \Big) 
+ 
\E\Big[ p\Big(\frac{\|Z\|}{s}\Big) \I\big[\|Z\|\leq s\big]\Big] 
\bigg) 
u
,\spa \forall\ s>0,\ \forall\ u\in\SS^{d-1}
.
$$
%	
%	g(s)
%	=
%	\sum_{i=0}^{\frac{d-3}{2}} \Big( \E\Big[\frac{1}{\|Z\|^{2i+1}}\Big] w_i s^{2i+1} \Big) 
%	+ 
%	\E\Big[ p\Big(\frac{\|Z\|}{s}\Big) \I\big[\|Z\|\leq s\big]\Big] 
%	,\spa \forall\ s>0
%	$$
\end{Prop}

It follows from the proof that $p$ is in fact a (Laurent) polynomial of the form
$$
p(y)
=
\sum_{i=0}^{\frac{d-1}{2}} \alpha_i y^{d-1-2i} 
+
\sum_{i=1}^{\frac{d-1}{2}} \beta_i \frac{1}{y^{2i-1}}
,\spa 
\forall\ y\in (0,1]
,
$$
with $\alpha_0\neq 0$ and $\beta_i\neq 0$ for all $i\in\{1,2,\ldots, (d-1)/2\}$.
%
%\begin{Prop}\label{ExplicitSpherOdd}
%	Fix $d=3$. Let $P$ be a spherically symmetric Borel probability measure on $\R^3$, i.e. $P[O A]=P[A]$ for all Borel subset $A\subset \R^3$ and $3\times 3$ orthogonal matrix $O$, and let $Z$ be a random vector with law $P$. If $\E[\|Z\|^2]<\infty$, then 
%%	 and $p(y):= (1 - (5/3)y^2 + (2/3)y^{-1})\I[0\leq y\leq 1]$ for all $y$, we have 
%	$$
%	\Fpg(su)
%	=
%	\Bigg( 
%	1
%	-
%	\frac{5}{3} \frac{1}{s^2}  \E\big[\|Z\|^2\big]
%	-
%	\E\bigg[ \bigg(1 - \frac{5}{3} \frac{\|Z\|^2}{s^2} + \frac{2}{3} \frac{s}{\|Z\|} \bigg) \I\big[\|Z\|>s\big]\bigg] 
%	\Bigg) 
%	u
%%	\E\Big[p\Big(\frac{\|Z\|}{\|x\|}\Big)\Big]\  
%%	,\spa \forall\ x\in\R^d
%	$$
%	for all $s>0$ and $u\in\SS^{2}$.
%%	$$
%%	\E\bigg[ \bigg(1 - \frac{5}{3} \frac{\|Z\|^2}{\|x\|^2} + \frac{2}{3} \frac{\|x\|}{\|Z\|} \bigg) \I\big[\|Z\|<\|x\|\big]\bigg] 
%%	$$
%%	P[\BB^d_s] - \frac{5}{3} \frac{1}{s^2} \E\big[\|Z\|^2\ \I[\|Z\|<s]\big] + \frac{2}{3} u \E\Big[\frac{1}{\|Z\|}\ \I[\|Z\|<S]\Big]
%%	$$
%%	,\spa \forall\ s>0,\ \forall\ u\in\SS^{d-1}
%%	.
%%	$$
%\end{Prop}
We deduce from Proposition \ref{ExplicitSpherOdd} that if $P$ and $Q$ are spherically symmetric probability measures on $\R^d$ with bounded densities that coincide on a ball of radius $R>0$ and have equal moments of order $-(2i+1)$, $i=0,1,\ldots, (d-3)/2$, then $F_P^\g(x)=F_Q^\g(x)$ for all $\|x\|\leq R$. In particular, since such probability measures obviously exist, we deduce that Proposition~\ref{PropLocalOdd} is not an empty statement. When $d=3$, a straightforward computation following the reasoning of the proof of Proposition \ref{ExplicitSpherOdd} gives
$$
\Fpg(su)
=
\Bigg( 
\E\bigg[ \bigg(1 - \frac{5}{3} \frac{\|Z\|^2}{s^2} + \frac{2}{3} \frac{s}{\|Z\|} \bigg) \I\big[\|Z\|<s\big]\bigg] 
- \frac{2}{3} s\ \E\Big[\frac{1}{\|Z\|}\Big]
\Bigg)
u
.
$$

%Observe the similarity with the univariate case $d=1$, where one directly verifies that the cdf's of $P$ and $Q$ coincide on $\R\sm [-R,R]$ provided $P$ and $Q$ do. It is also worth noting that requiring second moments in Proposition \ref{ExplicitSpherOdd} is not essential. Indeed, provided $P$ has a bounded density near the origin so that $\E[\|Z\|^{-1}]<\infty$, one has 
%$$
%\Fpg(su)
%=
%\Bigg( 
%\E\bigg[ \bigg(1 - \frac{5}{3} \frac{\|Z\|^2}{s^2} + \frac{2}{3} \frac{s}{\|Z\|} \bigg) \I\big[\|Z\|<s\big]\bigg] 
%- \frac{2}{3} s\ \E\Big[\frac{1}{\|Z\|}\Big]
%\Bigg)
%u
%.
%$$
%We deduce that if $P$ and $Q$ admit bounded densities near the origin, are spherically symmetric, coincide on a ball $\BB^3_R$, and have equal moments of order $-1$, then $\Fpg(x)=F_Q^\g(x)$ for all $\|x\|\leq R$. This statement, in turn, does not extend to the univariate case $d=1$. Indeed, there exist two probability measures on $\R$ that coincide on an interval $[-R,R]$ such that their cdf's do not. 
%In particular, having a bounded density near the origin is not enough when $d=1$ since univariate densities do not possess moments of order $-1$ in general; for this, they must vanish fast enough at the origin.

While it is relatively straightforward to show, as in \cite{Kol1997}, that the geometric cdf characterises the underlying probability measure in \emph{arbitrary} dimension, the difference in behaviour when $d$ is odd or even, described in Proposition \ref{PropLocalOdd} and Proposition \ref{PropLocalEven}, is quite unexpected. So far we don't have an intuitive or geometric way to explain this phenomenon from the geometric cdf itself; it is rather an analytical consequence of the nature of the reconstruction operator $\LL_d$ that is local when $d$ is odd and nonlocal when $d$ is even (through the fractional Laplacian $(-\Delta)^{1/2}$). 

%displayed by the geometric that the geometric cdf allows for a local characterisation

%It is natural to embed $P$ into $\R^{d+1}$, when $d$ is even, to try to recover a better understood and computable cdf in odd dimension. From a PDE perspective, 
Letting $u:=\gamma_d\ (-\Delta)^{(d-2)/2}(\nabla\cdot \Fpg)$ when $d$ is even, \cite{CaffarelliSilvestre2008} and \cite{Silvestre2007} introduced a procedure to compute $f_P = (-\Delta)^{1/2}u$ from another quantity, $U$ say, by applying only local differential operators to $U$; in particular, this allows one to recover a notion of locality with respect to $U$ (but not with respect to $u$). This approach consists in
%can be performed this amounts to computing $(-\Delta)^{1/2} u$ through an auxiliary function $U$ the equation $(-\Delta)^{1/2} u = f$, defined on $\R^d$, by embedding $f$ in $\R^{d+1}$. This approach, introduced in \cite{CaffarelliSilvestre2008} and \cite{Silvestre2007}, is performed by  in computing $(-\Delta)^{1/2} u$ by 
% We then present two approaches attempting to recover a localisation result similar to Proposition \ref{PropLocalOdd}.\\
%
%Fix $d$ even, and consider a probability measure $P$ on $\R^d$.  The first idea that naturally comes to mind is to embed $P$ into $\R^{d+1}$ (with $d+1$ odd); this gives rise to a probability measure $P^*$ supported on the hyperplane $x_{d+1}=0$ of $\R^{d+1}$. Proposition \ref{PropLocalOdd} then applies to $P^*$. The other approach consists in localising the operator $(-\Delta)^{1/2}$. For a smooth function $g$ on $\R^d$, computing $(-\Delta)^{1/2} g$ can be achieved by 
first solving $-\Delta U = 0$ on $\R^{d+1}_+:=\R^d\times (0,\infty)$, subject to the boundary condition~$U(x,0) = u(x)$ for all $x\in\R^d$, and then recovering $f_P=(-\Delta)^{1/2} u$ as 
$$
f_P(x)
=
((-\Delta)^{1/2} u)(x)
=
-
\lim_{\ve\downarrow 0} (\partial_{d+1} U)(x,\ve)
,\spa \forall\ x\in\R^d
.
$$
Because the values of $\partial_{d+1} U$ in some open subset $\Om\subset \R^{d+1}_+$ depend on the values of $U$ on $\Om$ only, this formulation is now local with respect to $U$. It is important to note, however, that this does not yield a local relation between $u$ and $f$; the values of $U$ on an open subset $\Om$ of $\R^{d+1}_+$, for instance $\Om=\tilde{\Om}\times (a_i, b_i)$ for some open subset $\tilde{\Om}\subset\R^d$, are not determined solely by the values of $u$ on $\tilde{\Om}$. Another very natural idea is to embed $P$ in $\R^{d+1}$, still when $d$ is even: defining $P^*$ as the probability measure on $\R^{d+1}$ supported on the hyperplane $x_{d+1}=0$ with density $f_P$ with respect to the $d$-dimensional Hausdorff measure $\HH_d$ on $\R^{d+1}$, the geometric cdf $F_{P^*}^\g$ of $P^*$ now enjoys the local characterisation property of Proposition \ref{PropLocalOdd}. In the next Proposition we show that both approaches are in fact equivalent, with $U$ essentially given by~$(-\Delta)^{(d-2)/2}(\nabla\cdot F_{P^*}^\g)$. We recall that $\HH_d$ stands for the $d$-dimensional Hausdorff measure on $\R^{d+1}$.

\begin{Prop}\label{PropEquivLocalMethods}
Fix an even integer $d\geq 2$. Let $\Om\subset\R^d$ be an open subset and $P$ a Borel probability measure on $\R^d$. Assume that $P$ admits a bounded density $f_P$ on $\R^d$ with respect to the Lebesgue measure and that $f_P\in \CC^{0,\beta}_\loc(\Om)$ for some $\beta\in (0,1)$. Let $P^*$ be the probability measure on $\R^{d+1}$ supported on the hyperplane $x_{d+1}=0$ with density $f_P$ with respect to $\HH_d$. 
%Let $Z$ be a random $d$-vector with law $P$, and $Z^*$ a random $(d+1)$-vector with law $P^*$. 
Define 
%$$
%G(x)
%=
%2\gamma_{d+1} \E\Big[((-\Delta)^{\frac{d-2}{2}}(\nabla\cdot K_{d+1}))(x-Z^*)\Big]
%,\spa 
%\forall\ x\in\R^{d+1}_+:=\R^d\times (0,\infty)
%,
%$$
$
u(x)
=
\gamma_d\ (-\Delta)^{(d-2)/2}(\nabla\cdot \Fpg)(x)
$
for all $x\in\Om$, and 
$
U(x^*)
=
2\gamma_{d+1}\ (-\Delta)^{(d-2)/2}(\nabla\cdot F_{P^*}^\g)(x^*)
$
for all $x^*\in\R^{d+1}_+:=\R^d\times (0,\infty)$. Then, (i) $U\in \CC^\infty(\R^{d+1}_+)$ and $-\Delta U=0$ on $\R^{d+1}_+$, (ii)~$U$ extends continuously on $\Om\times [0,\infty)$ with $U(x,0)=u(x)$ for all $x\in\Om$, and (iii) we have
$$
f_P(x)
=
((-\Delta)^{1/2} u)(x)
=
\lim_{\ve \downarrow 0}
-(\partial_{d+1} U)(x,\ve)
,\spa 
\forall\ x\in\Om
.
$$
%In addition, the following holds:
%\begin{enumerate}
%	\item $G(x)=2\gamma_{d+1}\ (-\Delta)^{\frac{d-2}{2}}(\nabla\cdot R_{P^*})(x)$ and $-\Delta G(x)=0$, for any $x\in\R^{d+1}_+$ ;\vspace{2mm}
%	\item for any $\tilde{x}\in\R^d$, $G(\tilde{x},0)=u(\tilde{x})$ and
%	$$
%	f_P(\tilde{x})
%	=
%	((-\Delta)^{1/2} g)(\tilde{x})
%	=
%	\lim_{x_{d+1} \downarrow 0}
%	-(\partial_{d+1} G)(\tilde{x},x_{d+1})
%	.
%	$$
%\end{enumerate}
\end{Prop}
%Let us stress that the map $U$ defined in Proposition \ref{PropEquivLocalMethods} is pointwise well-defined on $\R^{d+1}_+$. Indeed, since $P^*$ admits the null density on $\R^{d+1}_+$, then Theorem \ref{TheorRankRegularityOdd} entails that $F_{P^*}^\g\in \CC^\infty(\R^{d+1}_+)$.
%In addition, 
The requirement $f_P\in \CC^{0,\beta}_\loc(\R^d)$ is only used to ensure that $f_P(x)=((-\Delta)^{1/2} u)(x)$. However, the fact that $f_P(x)=\lim_{\ve \downarrow 0} -(\partial_{d+1} U)(x,\ve)$, which is the main result of Proposition \ref{PropEquivLocalMethods}, only requires $f_P$ to be continuous over $\Om$. In this case, it is  straightforward to see from the proof that the convergence in Part (iii) of Proposition \ref{PropEquivLocalMethods} is uniform over the compact subsets of $\Om$. In addition, a quantiative result can easily be derived if we assume, for instance, that $f_P\in \CC^{0,\beta}(\R^d)$ for some $\beta\in (0,1)$; in this case the limit as $\ve\downarrow 0$, which is uniform in $x\in\R^d$, converges at speed $\ve^\beta$. More generally, a quantitative result can be obtained in terms of the modulus of continuity of $f_P$. 

The fact that we assumed that $P$ has a density is not essential to prove that $P$ can be recovered from $-\partial_{d+1} U$. Indeed, define $\rho(y)=C_d\ (1+\|y\|^2)^{-(d+1)/2}$ for all $y\in\R^d$, where $C_d$ is a normalizing constant such that $\int_{\R^d} \rho(y)\, dy=1$. Letting $\rho_\ve(y):=\ve^{-d} \rho(y/\ve)$ for all $\ve>0$, then the proof of Part (iii) of Proposition \ref{PropEquivLocalMethods}, still with $d$ even, reveals that for an \emph{arbitrary} probability measure $P$, and $P^*$ the corresponding embedding of $P$, we have $F_{P^*}^\g\in\CC^\infty(\R^{d+1}_+)$ and
$$
%-2\gamma_{d+1}\ \partial_{d+1} (-\Delta)^{\frac{d-2}{2}} (\nabla\cdot F_{P^*}^\g)(x,\ve)
-(\partial_{d+1} U)(x,\ve)
=
\int_{\R^d} \rho_\ve(x-z)\, dP(z)
,\spa 
\forall\ x\in \R^d,\ \forall\ \ve>0 
.
$$
Considering each map $x\mapsto (\partial_{d+1} U)(x,\ve)$, for a fixed $\ve>0$, as a distribution gives rise to a collection $\{(\partial_{d+1} U)(\cdot, \ve) : \ve > 0\}$ in $\D(\R^d)'$. Then, it is easy to show that the equality
$
P 
=
\lim_{\ve\downarrow 0} -(\partial_{d+1} U)(\cdot, \ve)
$
holds in 
$\D(\R^d)'$; the limit of distributions is to be understood as a pointwise limit by evaluation on test functions from $\CC^\infty_c(\R^d)$.

One can in principle recover $P^*$, hence $f_P$, from $F_{P^*}^\g$ through the equality $P^*=\gamma_{d+1}(-\Delta)^{d/2}(\nabla\cdot F_{P^*}^\g)$ established in Theorem \ref{TheorEDPDistributions}. This identity, however, will be useful in practice only when $P^*$ has a density with respect to the ambiant space Lebesgue measure. But $P^*$ has a density with respect to a lower dimensional Hausdorff measure. In particular, because $P^*$ has a vanishing density on the complement of the hyperplane $H$ of equation $x_{d+1}=0$, Theorem \ref{TheorRankRegularityOdd} yields 
$$
\gamma_{d+1}\ (-\Delta)^{d/2}(\nabla\cdot F_{P^*}^\g)(x^*)
=
\frac12 (-\Delta U)(x^*)
=
0,\spa 
\forall\ x^*\in \R^d\sm H 
.
$$
Proposition \ref{PropEquivLocalMethods} then entails that in order to actually see the density $f_P$ pointwise, which lives in a lower dimensional subspace, less derivatives of $F_{P^*}^\g$ must be computed. In fact, this procedure explains how one can retrieve the density of a probability measure $P$ supported on a lower dimensional hyperplane (of co-dimension $1$). We extend it to hyperplanes of arbitrary \emph{even} dimension in the next Proposition.

\begin{Prop}\label{PropLowerScaleEven}
Fix $d\geq 3$ and a probability measure $P$ on $\R^d$ supported on a hyperplane $H$ of even dimension $\kappa\in [2,d)$ with density $f_P\in \CC^0_b(H)$ with respect to the $\kappa$-dimensional Hausdorff measure. Then $\Fpg\in \CC^\infty(\R^d\sm H)$ and for any unit vector $\nu$ orthogonal to $H$ we have
$$
f_P(x)
=
-\frac{\Gamma(\frac{d-\kappa+1}{2})}{(4\pi)^{\frac{\kappa}{2}} \Gamma(\frac{d+1}{2})}-
\lim_{\ve\downarrow 0} \big(\partial_{\nu} (-\Delta)^{\frac{\kappa}{2}-1}(\nabla\cdot F_P^\g)\big)(x+\ve \nu)
,\spa 
\forall\ x\in H 
,
$$
where $\partial_\nu =\sum_{i=1}^d \nu_i \partial_i$ denotes the directional derivative in direction $\nu$.
\end{Prop}

The same argument as the one following Proposition \ref{PropEquivLocalMethods} entails that, more generally, the equality 
$$
P
=
-\frac{\Gamma(\frac{d-\kappa+1}{2})}{(4\pi)^{\frac{\kappa}{2}} \Gamma(\frac{d+1}{2})} 
\lim_{\ve\downarrow 0} \big(\partial_{\nu} (-\Delta)^{\frac{\kappa}{2}-1}(\nabla\cdot F_P^\g)\big)(\cdot+\ve \nu)
$$
holds in $\D(\R^d)'$ even when $P$ admits no density. This reconstruction procedure has the advantage to involve local differential operators only. It is worth noticing that, while a usual $d$-dimensional Lebesgue density in $\R^d$ is recovered from the $d$th derivatives of $\Fpg$ (see Theorem \ref{TheorRankRegularityOdd}), it follows from Proposition~\ref{PropLowerScaleEven} that a density living at scale $\kappa<d$, with $\kappa$ even, is reconstructed from the derivatives of $\Fpg$ of order $\kappa$. For $\kappa=0$, \ie when $P$ is atomic, a straightforward computation yields
$$
%\|F_P^\g(x)-F_P^\g(x_0)\| 
%=
%P[\{x_0\}] + O(\|x-x_0\|)
P[\{x\}]
=
\lim_{\ve\downarrow 0} \ps{F_P^\g(x+\ve \nu)-F_P^\g(x)}{\nu}
%\lesssim 
%\ve
,\spa 
\forall\ x\in\R^d,\ \forall\ \nu\in\SS^{d-1}
%\quad
%\text{as }\ x\to x_0
.
$$
Here again we see that the ``density of $P$ at scale $\kappa=0$'', \ie the weight function $x\mapsto P[\{x\}]$, is recovered from $\Fpg$ through differences of order $0$. While we believe that an analogue of Proposition~\ref{PropLowerScaleEven} holds when $\kappa$ is odd, with another differential operator of order $\kappa$, our proof strategy does not apply and the existence of such a result is still open.

\section{Depth regions}\label{sec:DepthRegions}

Fix $d\geq 2$ and a probability measure $P$ on $\R^d$. For any $\beta\in [0,1)$ and $u\in \SS^{d-1}$ recall that a geometric quantile of order $\beta$ in direction $u$ for $P$ is an arbitrary minimizer of the objective function $O_{\beta,u}^P$ introduced in Section \ref{sec:Introduction}. Further assuming that $P$ is not supported on a single line of $\R^d$, Theorem 1 in \cite{PaiVir2021} implies that such a quantile is unique; let us denote it by $\Qpg(\beta u)$. This gives rise to the \emph{geometric quantile map} $Q_P^\g : \mathbb{B}^d\to \R^d$ which, as we explained in Section \ref{sec:Introduction}, is a natural multivariate analogue of the usual univariate quantile function. As such, it allows one to compute centrality regions associated with $P$, \ie an increasing family of nested sets that each contain the center of $P$ (the geometric median). They provide invaluable tools to perform various statistical inference tasks such as depth-based classification, outliers detection, analysis of extremes, or testing some symmetry. In a broader context, such regions are often called \emph{depth regions} in reference to the more general theory of statistical depth, and in our context they lead to \emph{geometric depth regions}.

%define the geometric quantile regions $\mathcal{D}_P^\beta$ and contours $\mathcal{C}_P^\beta$ of arbitrary order $\beta\in [0,1)$ in the next definition. 

\begin{Defin}
Fix $d\geq 2$ and $P$ a probability measure on $\R^d$. Assume that $P$ is not supported on a single line of $\R^d$. For any $\beta\in [0,1)$, the geometric depth region $\mathcal{R}_P^\g(\beta)$ and geometric depth contour $\mathcal{C}_P^\g(\beta)$ of order $\beta$ for $P$ are defined as
$$
\mathcal{R}_P^\g(\beta)
=
\big\{ \Qpg(\alpha u) : \alpha \in [0,\beta],\ u\in \SS^{d-1} \big\}
\quad\quad 
\textrm{and} 
\quad\quad
\mathcal{C}_P^\g(\beta)
=
\big\{\Qpg(\beta u) : u\in \SS^{d-1}\big\}
.
$$

\end{Defin}

When $P$ is not supported on a line of $\R^d$ and non-atomic, then Theorem 6.1 in \cite{KonPai1} entails that the map $Q_P$ is a homeomorphism with inverse $\Fpg$. In particular, as the continous images of the closed ball and the sphere of radius $\beta$, respectively, the geometric depth regions and contours are compact and connected. In addition, the depth regions $\{\RR_P^\g(\beta) : \beta\in [0,1)\}$ are increasing, \ie $\RR_P^\g(\beta_1)\subset \RR_P^\g(\beta_2)$ if $\beta_1\leq \beta_2$, and the depth contours $\{ \mathcal{C}_P^\g(\beta) : \beta\in [0,1)\}$ are disjoint. Although depth regions are convex in most cases, they may fail to be convex in general; see Figure 2, Figure 5, and Figure 7 in Section 4 of \cite{GirStu2017}. We also refer to \cite{Nag2017} for a quantified discussion on the shape of depth regions. 

In the rest of this section, we investigate a few properties of geometric regions and contours. First, we show an explicit construction that allows one to make sure that the probability content of $\RR_P^\g(\beta)$ is equal to $\beta$. Then, we establish the regularity of geometric contours thanks to the results we proved in Section \ref{sec:Regularity}. 

%We conclude by showing that, quite surprisingly, the geometric depth function, which is essentially given by the map $x\mapsto \|\Fpg(x)\|$, fully characterises $P$.

\subsection{Relabeling}\label{sec:Relabeling}

Unlike center-outward quantiles based on optimal transport (see, e.g., \cite{HallEspagne21}), geometric quantile regions are not indexed by their probability content, \ie we do not have $P[\RR_P^\g(\beta)]=\beta$ in general. To see this, consider a probability measure $P$ which is spherically symmetric, \ie $P[O A] = P[A]$ for any $d\times d$ orthogonal matrix $O$ and Borel set $A$; further assuming that $P$ is not a Dirac at $0$, this implies that $P$ has no atoms and is not supported on a single of $\R^d$ so that $Q_P$ is a homeomorphism. Proposition 2.2 in \cite{GirStu2017} then entails that there exists a continuous and increasing map $g : [0,1)\to [0,\infty)$ such that $\Qpg(\beta u)=g(\beta)u$ for all $\beta\in [0,1)$ and $u\in\SS^{d-1}$. Since $Q_P$ is a homeomorphism, we further have $g(\beta)\uparrow \infty$ as $\beta\uparrow 1$. In particular, each depth region $\RR_P^\g(\beta)$ is a ball centered at the origin with radius $g(\beta)$. Now, let us further assume that the support of $P$ is compact and contained in $\BB^d_r$ for some $r>0$. It follows that for all $\beta\in [0,1)$ large enough so that $g(\beta)>r$, we have $\BB^d_r\subset \RR_P^\g(\beta)$; in particular, for any such $\beta$ we have $1=P[\RR_P^\g(\beta)]>\beta$.  Although this issue is a conceptual flaw of geometric quantiles, it has no practical impact since one can always re-label the depth regions so that they match their probability content. For this purpose let us define, for an arbitrary non-atomic probability measure $P$ that is not supported on a line, the map $w_P : [0,1)\to[0,1)$ by letting $w_P(\beta)=P[\RR_P^\g(\beta)]$ for all $\beta\in [0,1)$. Recalling that $\Qpg$ is thus invertible with inverse $\Fpg$ yields
\begin{equation}\label{eq:DefWeight}
w_P(\beta)
=
P\big[\|\Fpg(Z)\|\leq \beta\big]
,\spa \forall\ \beta\in [0,1)
,
\end{equation}
where $Z$ is a random $d$-vector with law $P$. It follows that the map $w_P$ is, in fact, the cdf of the univariate random variable $\|\Fpg(Z)\|$. This implies that $w_P(\|\Fpg(Z)\|)$ is uniformly distributed on $[0,1)$, \ie 
$$
P\big[w_P(\|\Fpg(Z)\|)\leq \beta\big]
=
\beta
,\spa \forall\ \beta\in [0,1)
.
$$
This motivates the following definition, that allows $P$ to have atoms or be supported on a single line.

\begin{Defin}\label{DefRelabel}
Fix $d\geq 2$ and $P$ a Borel probability measure on $\R^d$. Let $w_P:[0,1)\to [0,1)$ be the map defined by (\ref{eq:DefWeight}). The relabeled geometric cdf of $P$ is the map $\widetilde{F}_P^\g:\R^d\to \overline{\BB^d}$ defined as
\begin{equation}\label{eq:ReIndexRank}
	\widetilde{F}_P^\g(x)
	=
	w_P(\|\Fpg(x)\|) \frac{\Fpg(x)}{\|\Fpg(x)\|} \I\big[\|\Fpg(x)\|\neq 0\big]
	,\spa 
	\forall\ x\in\R^d
	.
\end{equation}
\end{Defin}
%Consequently, we can define the \emph{relabeled geometric cdf} of $P$ as 

The relabeled geometric cdf naturally gives rise to the corresponding relabeled regions and contours 
$$
\widetilde{\RR}_P^\g(\beta)
=
\big\{z\in\R^d : \|\widetilde{F}_P^\g(z)\|\leq \beta\big\}
\quad\quad 
\textrm{and} 
\quad\quad
\widetilde{\CC}_P^\g(\beta)
=
\big\{z\in\R^d : \|\widetilde{F}_P^\g(z)\|= \beta\big\}
.
$$
When $P$ has no atoms and is not supported on a line, we then have $P[\widetilde{\RR}_P^\g(\beta)]=\beta$ for all~$\beta\in [0,1)$. It is important to note that the relabeled regions and contours are the same, as collections of sets, as the original ones; only their enumeration changed. 
%In particular, this re-indexed geometric cdf has quantile regions that match their probability content. 
%	Let us make some comments about the behaviour of $\widetilde{R}_P$ when $\|\Fpg(x)\|\to 0$. For this purpose, we need to understand $\theta_P(\beta)$ when $\beta\to 0$, \ie the cdf of $\|\Fpg(Z)\|$ around $0$. 
%	\begin{Prop}
%		Let $d\geq 1$ and $P$ a Borel probability measure on $\R^d$. Let $Z$ be a random $d$-vector with law $P$, and assume that $P$ admits a bounded density $f_P$ with respect to the Lebesgue measure. Then the law of $\|\Fpg(Z)\|$ admits a density with respect to the Lebesgue measure over $[0,1]$.
%	\end{Prop}
The relabeled geometric cdf can easily be estimated in practice by considering its empirical counterpart: denoting by $P_N$ the empirical probability measure associated with a random sample $Z_1,\ldots, Z_N$ of size $N$ drawn from $P$, consider 
\begin{equation}\label{eq:ReIndexRankEmpirical}
\widetilde{F}_{P_N}^\g(x)
:=
w_{P_N}(\|F_{P_N}^\g(x)\|) \frac{F_{P_N}^\g(x)}{\|F_{P_N}^\g(x)\|} \I\big[\|F_{P_N}^\g(x)\|\neq 0\big]
,\spa \forall\ x\in\R^d
,
\end{equation}
where 
$$
F_{P_N}^\g(x)
=
\frac{1}{N}
\sum_{j=1}^N \frac{x-Z_j}{\|x-Z_j\|}\I[Z_j\neq x]
%	,\spa \forall\ x\in\R^d
\quad 
\text{and} 
\quad
w_{P_N}(\|F_{P_N}^\g(x)\|)
=
\frac{1}{N} \sum_{i=1}^N \I\big[\|F_{P_N}^\g(Z_i)\|\leq \|F_{P_N}^\g(x)\|\big]
.
$$
Since the points $\|\widetilde{F}_{P_N}(Z_1)\|,\ldots, \|\widetilde{F}_{P_N}^\g(Z_N)\|$ are uniformly distributed on $\{1/N, 2/N,\ldots, 1\}$, they have a natural interpretation as multivariate analogues of the usual univariate ranks. As such, they provide the basis to design multivariate (distribution-free) rank-based statistical procedures. To the best of our knowledge, the behaviour of $F_{P_N}^\g(Z_i)/\|F_{P_N}^\g(Z_i)\|$ has not been studied yet in the literature and remains unknown.

In order to be useful for inference purposes, a good notion of cdf should satisfy a Glivenko-Cantelli result, \ie the empirical cdf converges uniformly to the population one as the sample size becomes arbitrarily large. A Glivenko-Cantelli result has been proved for the geometric cdf when $P$ admits a bounded density; see Lemma 2 in \cite{Mottonen97}. The next Proposition establishes that this extends to its relabeled counterpart. 
%The next proposition is thus a consequence of this fact. For the sake of completness we provide a proof of Proposition \ref{PropReIndexedRanksConv} in Appendix \ref{Appendix:ReIndexedUnifConv}.

\begin{Prop}\label{PropReIndexedRanksConv}
Fix $d\geq 1$ and $P$ a probability measure on $\R^d$. Assume that $P$ admits a bounded density with respect to the Lebesgue measure. Let $P_N$ denote the empirical measure associated with a random sample of size $N$ drawn from $P$. Then, we have
$$
\sup_{x\in\R^d} \|\widetilde{F}_{P_N}^\g(x)-\widetilde{F}_P^\g(x)\|
\to 0 
$$
$P$-almost surely as $N\to\infty$.
\end{Prop}

In particular, when $Z_1,\ldots, Z_N$ is a random sample drawn from $P$, Proposition \ref{PropReIndexedRanksConv} implies that
$$
\max_{i=1,\ldots, N} 
\|\widetilde{F}_{P_N}^\g(Z_i)-\widetilde{F}_P^\g(Z_i)\|
\to 
0 
$$
$P$-almost surely as $N\to\infty$. While the original proof of \cite{Mottonen97} requires $P$ to have a bounded density, we believe that this assumption can be dropped and that the Glivenko-Cantelli result holds under much weaker assumptions on $P$.

\subsection{Regularity of quantile maps and depth contours}\label{sec:Depth}

%In this section, we prove new results about geometric depths regions of an arbitrary probability measure on $\R^d$. We use our results from Section \ref{sec:PDEGeneral} to characterise the regularity of depth regions.\\

We now use the results of Section \ref{sec:Regularity} to establish the regularity of the quantile map $Q_P^\g$ and, subsequently, to determine the regularity of the geometric depth contours $\{\CC_P^\g(\beta) : \beta\in [0,1)\}$. In the next Proposition, we use the Inverse Function Theorem to derive the regularity of the geometric quantile map $Q_P^\g$ from that of $\Fpg$ established in Theorem \ref{TheorRankRegularityOdd}.

\begin{Prop}\label{PropRegQuantiles}
Fix $d\geq 2$ and $P$ a probability measure on $\R^d$. Assume that $P$ admits a density $f_P$ with respect to the Lebesgue measure. (i) If $f_P\in L^p_\loc(\R^d)$ for some $p\in (d/(d-\ell),\infty]$ and an integer $\ell\in [1,d-1]$, then $\Qpg:\mathbb{B}^d\to\R^d$ is a diffeomorphism of class $\CC^\ell$ with inverse $\Fpg$. (ii) If $f_P\in \CC^{k,\beta}_\loc(\R^d)$ for some $k\in\N$ and $\beta\in (0,1)$, then $\Qpg$ is a diffeomorphism of class $\CC^{k+d}$ with inverse $\Fpg$. 
%		(iii) If $f_P\in\CC^\infty(\R^d)$
%		\end{enumerate} 
\end{Prop}

%	To state regularity properties of depth contours, let us first rewrite depth contours in terms of the cdf map $\Fpg$. Theorem 6.1 in \cite{KonPai1} entails that $x=Q_P(\alpha u)$ if and only if $\Fpg(x)=\alpha u$. Consequently, we have
%	$$
%	\mathcal{D}_P^\beta 
%	=
%	\Big\{ x\in\R^d : \|\Fpg(x)\|\leq \beta \Big\}
%	\quad\textrm{and}\quad
%	\mathcal{C}_P^\beta
%	=
%	\Big\{x\in\R^d : \|\Fpg(x)\| = \beta\Big\}
%	.
%	$$

Each geometric depth contour $\CC_P^\g(\beta)$ being the image under $\Qpg$ of the smooth manifold $\beta \SS^{d-1}$, Proposition \ref{PropRegQuantiles} can be used to determine the minimal expected regularity of all geometric contours.

\begin{Corol}\label{CorolRegularContour}
Fix $d\geq 2$ and $P$ a probability measure on $\R^d$. Assume that $P$ admits a density $f_P$ with respect to the Lebesgue measure. Fix $\beta\in [0,1)$. (i) If $f_P\in L^p_\loc(\R^d)$ for some $p\in (d/(d-\ell),\infty]$ and an integer $\ell\in [1,d-1]$, then $\CC_P^\g(\beta)$ is a $(d-1)$-dimensional manifold of class $\CC^\ell$. (ii) If $f_P\in \CC^{k,\alpha}_\loc(\R^d)$ for some $k\in\N$ and $\alpha\in (0,1)$, then $\CC_P^\g(\beta)$ is a $(d-1)$-dimensional manifold of class $\CC^{k+d}$.
%		\end{enumerate} 
\end{Corol}

Although these results may seem natural in the context of statistical depth---regular probability measures give rise to regular depth contours ---, they are shared by very few concepts of depth in general. Not only are these types of results hard to establish for other existing depth notions, but depth contours sometimes display singularities. For instance, Tukey's halfspace depth contours may have infinitely many singular points even for mixtures of Gaussian measures; see \cite{GijNag2016}. The regularity of halfspace depth regions can be difficult to establish in practice, and is often a prerequiste to estalish theoretical properties of the halfspace depth. More recently, \cite{NagEtAl2019} showed that the regularity of a halfspace depth region for $P$ implies the existence of a certain convex floating body for $P$, which is an important geometric constraint on the depth region; see Proposition 31 in \cite{NagEtAl2019}.

To conclude this section, let us mention that whether the depth function of a probability measure $P$ itself charachterises $P$ or not is still an open question. More precisely, if $P$ and $Q$ are two probability measures such that $\|\Fpg(x)\| = \|F_Q^\g(x)\|$ for all $x\in\R^d$, do we necessarily have $P=Q$? Note that this requirement is stronger than simply asking that the depth regions of $P$ and $Q$ be identical; it further requires that the labeling of these regions coincide for both measures as well. Equality of the family of depth regions alone is not enough to guarantee $P=Q$, in general; for instance, if $P$ and $Q$ are distinct spherically symmetric probability measures, then the family of contours is the collection of all spheres $\{r\SS^{d-1} : r > 0\}\cup\{0\}$ centered at the origin. Although no proof or counter-example exists in the literature, we believe the characterisation property of depth functions to be true. 

A related conjecture, which is specific to spherically symmetric probability measures, can be stated as follows: if all geometric depth regions of $P$ are hyperspheres, then $P$ is spherically symmetric. We believe that solving these open questions would be valuable for at least two reasons. First, the techniques used to prove either of these conjectures would arguably provide further insights on the way geometric quantiles and cdf's work. Second, this would legitimate using depth-based methods in statistical inference as universally consistent procedures.

%\subsection{Characterization property of the geometric depth}
%
%
%+ conjecture sphérique 
%
%Pour tout $\alpha\in [0,1]$, we must have $Q_P(\alpha \SS^{d-1}) = O Q_P(\alpha \SS^{d-1})$, \ie 
%$$
%\alpha \SS^{d-1} 
%=
%F_P ( O Q_P(\alpha \SS^{d-1}))
%.
%$$
%$
%\|F_P(O x)\| 
%=
%\alpha
%$
%pour tout $x\in \RR_P(\alpha)$, \ie $\|F_P(O x)\| = \|F_P(x)\|$ pour tout $x$.
%
%$$
%\int_0^1 \int_{\RR_P^\g(\alpha)} \big|\|F_P(O x)\| - \alpha\big|\, d\sigma_\alpha(x)\, d\alpha 
%$$

%\section{Examples in dimension 2 and 3} \label{sec:Examples}
%
%blabla

\section{Final comments}\label{sec:Conclusion}

While we provided a complete solution to the inverse problem of recovering a probability measure $P$ from its geometric cdf $\Fpg$, our results raise other questions and call upon various potential extensions. 

We established in Proposition \ref{ExplicitSpherOdd} an explicit expression for $\Fpg$ when $P$ is a spherically symmetric probability measure and $d$ is odd. This sheds some light on results obtained in \cite{GirStu2017} (see Theorem 2.2) and \cite{SinEtAl2023} (see Theorem 2.5) that show that the behaviour of extreme geometric quantiles is fully determined by a fixed number of moments---although our result only applies to spherical laws in odd dimension. An analoguous result cannot hold when $d$ is even since this would contradict the nonlocality property established in Proposition \ref{PropLocalEven}. We know, however, that there is a map $g:[0,1)\to[0,\infty)$ such that $\Fpg(su)=g(s)u$ for all $s>0$ and $u\in\SS^{d-1}$, and that the asymptotic behaviour of $g$ is the same as in odd dimension; a general formula for $g$ in this setting would fill the gap that exists between odd and even dimensions in the spherical case.

%We showed in Proposition \ref{PropLowerScaleEven} that 

Whether the geometric (or so-called \emph{spatial}) depth function $x\mapsto \|\Fpg(x)\|$ is enough to characterise $P$ or not is arguably one of the most important conceptual questions about geometric quantiles that remains unanswered. Proving this result or providing a counterexample would be a decisive insight to help moving the field forward: if it always characterises $P$ then inference can be performed through the depth function $\|\Fpg(x)\|$ only, which is a univariate quantity that leads to universally consistent procedures---in particular, the 'directional information' contained in $\Fpg(x)/\|\Fpg(x)\|$ is superfluous---, whereas if it is false then Serfling's DOQR paradigm---which states that depth and quantile functions correspond to each other (see, e.g., \cite{Ser2010}, \cite{SerZuo2010}, and \cite{Ser2019a})---is wrong and a quantile function can be strictly more informative than its depth counterpart.

While we focused on finite-dimensional Euclidean spaces, geometric quantiles of a probability measure $P$ supported on an infinite-dimensional Banach space can still be defined along with the corresponding geometric cdf. To the best of our knowledge, there is no result in the literature on the characterisation property of the geometric cdf in infinite dimension. More specifically, one can wonder how the operator $\LL_d$ behaves as $d\uparrow \infty$, and whether, similarly to the finite dimensional case, there even exists a linear operator through which $P$ is to be reconstructed from $\Fpg$. 

Finally, the most puzzling result of the present paper is arguably the two distinct behaviours displayed by the geometric cdf in odd and even dimension. It is interesting to note, however, that the nonlocality in even dimension is a direct consequence of the form of the Fourier transform of the kernel $K_d$ introduced in Section \ref{sec:PDEGeneral}. In particular, other choices of kernels would give rise to different cdf concepts that might not exhibit dimension-dependent features. As long as they satisfy minimal equivariance requirements and essentially coincide with the univariate cdf when $d=1$, these concepts would be legitimate multivariate cdf and quantile functions.

\section*{Acknowledgements}
We thank Prof. Stanislav Nagy for the many valuable comments and suggestions he made on the first preprint version of the manuscript. We would also like to thank Prof. Davy Paindaveine for his insightful comments and the help he provided during the preparation of this work, and Prof. Yvik Swan for his suggestions and remarks. We thank Prof. Bruno Premoselli for pointing out to us many references on fractional Laplacians. 
%We finally thank the Editor-In-Chief, Davy Paindaveine, the Associate Editor, and two anonymous referees for their insightful comments that led us to investigate some aspects further and resulted in a significant improvement of the manuscript.
%We would also like to thank Alberto Gonsalez-Sanz for pointing out 
%\textcolor{red}{Alberto, Bruno, two anonymous referees}

%%%%%%%%%%%%%%%%%%%%%%%%%%%%%%%%%%%%%%%%%%%%%%
%% Funding information, if any,             %%
%% should be provided in the                %%
%% funding section.                         %%
%%%%%%%%%%%%%%%%%%%%%%%%%%%%%%%%%%%%%%%%%%%%%%
%\begin{funding}
%The author was supported by an Aspirant fellowship from the FNRS (Fonds National pour la Recherche Scientifique), Communauté Française de Belgique. 
%\end{funding}

%\section*{Acknowledgements}
%The author thanks Davy Paindaveine for the comments, suggestions, and help he provided during this work. The author also thanks Stanislav Nagy and Yvik Swan for their precious comments and remarks on the manuscript, as well as the F.R.S.-FNRS for financial support.

\bibliographystyle{agsm}
\bibliography{mybib}

\appendix

\section*{Appendices}

	\section{Review of background material}\label{sec:Review}

	In this section, we review some tools of analysis we need to state and prove the main results of this paper: distribution theory and Sobolev spaces.
	
	\subsection{Distributions}\label{sec:Distributions}
	For proofs and further details on the material presented in this section, we refer the reader to \cite{RudinFunctional}, on which this section is based. 
	
	Fix be an open subset $\Om\subset\R^d$. The set of infinitely differentiable functions whose which a compact support included in $\Om$ is denoted by $\CC^\infty_c(\Om)$ or $\D(\Om)$. We endow $\CC^\infty_c(\Om)$ with the following topology: a sequence $(\varphi_k)\subset\CC^\infty_c(\Om)$ converges to $\varphi\in\CC^\infty_c(\Om)$ in the space $\CC^\infty_c(\Om)$ if there exists a compact subset $K\subset \Om$ with $\supp(\varphi_k)\subset K$ for all $k$, and such that
	$$
	\lim_{k\to\infty} 
	\sup_{x\in K} |\partial^\alpha(\varphi_k-\varphi)(x)|
	= 
	0
	,\spa \forall\ \alpha\in\N^d
	.
	$$
	%as $k\to\infty$ for any $\alpha\in\N^d$. 

	\begin{Defin}\label{DefinDistributions}
		A distribution on $\Om$ is a linear
		$
		T : \CC^\infty_c(\Om)\to \C,\ \varphi\mapsto \ps{T}{\varphi}_{\D',\D}
		$
		which is continuous with respect to the topology of $\CC^\infty_c(\Om)$. The set of all distributions on $\Om$ is denoted $\D(\Om)'$.
	\end{Defin}
	%. The space $\D(U)'$ is equipped with the pointwise convergence : we say that a sequence $(T_k)\subset\D(U)'$ converges to $T\in\D(U)'$ in $\D(U)'$ if 
	%$$
	%\ps{T_k}{\varphi}
	%\to 
	%\ps{T}{\varphi}
	%$$
	%as $k\to\infty$ for any $\varphi\in\CC^\infty_c(U)$. 
	%The examples belowWe give some important examples of distributions.
	
	%An important number of distributions can be obtained, as showed in the following examples.
	
	\begin{Examples}
		We list here typical examples of distributions and a few usual ways to obtain distributions from other distributions.
		
		\begin{enumerate}
			\item Any function $f\in L^1_\loc(\Om)$ gives rise to a distribution on $\Om$, also denoted $f$, by letting 
			$$
			\ps{f}{\varphi}_{\D',\D}
			:=
			\int_\Om f(x)\varphi(x)\, dx 
			,\spa 
			\forall\ \varphi\in\CC^\infty_c(\Om)
			.
			$$
			%	for any .
			
			\item Similarly, any Borel measure $\mu$ on $\Om$ that is finite over compact subsets of $\Om$ leads to a distribution on $\Om$ by letting 
			$$
			\ps{\mu}{\varphi}_{\D',\D}
			:=
			\int_\Om \varphi(x)\, d\mu(x)
			,\spa 
			\forall\ \varphi\in\CC^\infty_c(\Om)
			.
			$$
			In particular, any Borel probability measure is a distribution on any open subset of $\R^d$.
			
			\item If $T\in \D(\Om)'$ is a distribution on $\Om$, we define its distributional derivatives $\partial^\alpha T$, $\alpha\in\N^d$, by letting
			$$
			\ps{\partial^\alpha T}{\varphi}_{\D',\D}
			:=
			(-1)^{|\alpha|} \ps{T}{\partial^\alpha\varphi}_{\D',\D}
			,\spa 
			\forall\ \varphi\in\CC^\infty_c(\Om)
			.
			$$
			%	for any $\varphi\in\CC^\infty_c(U)$. 
			
			\item For any smooth function $f\in\CC^\infty(\Om)$ and distribution $T\in\D(\Om)'$ on $\Om$, we define the distribution $fT$ on $\Om$ by letting 
			$$
			\ps{fT}{\varphi}_{\D',\D}
			:=
			\ps{T}{f\varphi}_{\D',\D}
			,\spa 
			\forall\varphi\in\CC^\infty_c(\Om)
			.
			$$
			%	for any $\varphi\in\CC^\infty_c(U)$.
		\end{enumerate}  
		%\end{Examples}
	\end{Examples}
	
	%	In this paper, we will need to consider convolution between two distributions, which is well-defined when at leats one of them has a compact support: 
	
	We say that $T\in\D(\R^d)'$ vanishes on an open set $\Om\subset\R^d$ if $\ps{T}{\varphi}_{\D',\D}=0$ for all $\varphi\in\CC^\infty_c(\Om)$, and we call the support of $T$ the complement of the union of all open sets on which $T$ vanishes. When $T\in\D(\R^d)'$ and $\varphi\in\CC^\infty_c(\R^d)$, we define their convolution $T*\varphi:\R^d\to\C$ as 
	$$
	(T*\varphi)(x)
	:=
	\ps{T}{\psi(x-\cdot)}_{\D',\D}
	,\spa 
	\forall\ x\in\R^d 
	.
	$$
	In this case, we have $T*\varphi\in\CC^\infty(\R^d)$ and the equalities
	$$
	\partial^\alpha(T*\varphi)
	=
	(\partial^\alpha T)*\varphi 
	=
	T*(\partial^\alpha \varphi)
	$$
	hold in $\D(\R^d)'$; see Theorem 6.30 in \cite{RudinFunctional}. In addition, when $T$ has compact support then $T*\varphi$ has compact support as well, \ie $T*\varphi\in \CC^\infty_c(\R^d)$. This legitimates the next definition.
	
	\begin{Defin}\label{DefinConvolDistributions}
		Let $S,T\in\D(\R^d)'$ and assume that $T$ has compact support. We define the distribution $S*T$ by letting 
		$$
		\ps{S*T}{\varphi}_{\D',\D}
		:=
		\ps{S}{T*\varphi}_{\D', \D}
		,\spa\forall\ \varphi\in\CC^\infty_c(\R^d)
		.
		$$
	\end{Defin}
	
	The next proposition is given in Theorem 6.37 (e) of \cite{RudinFunctional}.
	
	\begin{Prop}\label{PropConvolGeneral}
		Let $S,T\in\D(\R^d)'$ and assume that $T$ is compactly supported. Then 
		$$\partial^\alpha(S*T)=(\partial^\alpha S)*T = S*(\partial^\alpha T)
		,\spa 
		\forall\ \alpha\in\N^d
		.
		$$
	\end{Prop}

	Distributions are closed with respect to multiplication by smooth functions, and by taking derivatives. Other common operations, such as convolution or Fourier transform, do not leave the space $\CC^\infty_c(\Om)$ invariant and, as a consequence, cannot be defined on $\D(\Om)'$. We thus need another class of test functions, and the corresponding new distributions, to apply these operations. This is the role of tempered distributions, that rely on the Schwartz class $\SC(\R^d)$ defined as
	$$
	\SC(\R^d)
	=
	\Big\{f\in \CC^\infty(\R^d) : \sup_{x\in\R^d} (1+\|x\|)^m |\partial^\alpha f(x)| < \infty,\ \forall\ m\in\N,\ \forall\ \alpha \in\N^d\Big\}
	.
	$$
	Following our definition of $\CC^\infty(\R^d)$, functions from $\SC(\R^d)$ are complex-valued (see Section \ref{sec:Introduction}).
	%\begin{Prop}
	%	\label{PropDensityTestInDistributions}
	%	Let $\Om\subset\R^d$ be an open subset and $T\in \D(\Om)'$ be a distribution on $\Om$. There exists a sequence $(\psi_k)\subset\CC^\infty_c(\Om)$ such that $\psi_k\to T$ in $\D(\Om)'$, \ie 
	%	$$
	%	\int_\Om \psi_k(x) \varphi(x)\, dx 
	%	\to 
	%	\ps{T}{\varphi}
	%	$$
	%	as $k\to\infty$ for any $\varphi\in\CC^\infty_c(\Om)$.
	%\end{Prop}
	The Schwartz class is a vector space. It is also closed with respect to scalar multiplication, multiplication by smooth functions all derivatives of which have at most polynomial growth at infinity, convolution, differentiation and Fourier transform. We further have the inclusion 
	$
	\CC^\infty_c(\Om)\subset\SC(\R^d)
	$
	for any open subset $\Om\subset\R^d$.\\
	
	%	\textcolor{red}{
		%		Attention : le produit $(\SC(\R^d)')^k=\SC(\R^d)'\times \ldots \times \SC(\R^d)'$ n'est pas la même chose que le produit tensoriel $\SC(\R^d)'\otimes \ldots \otimes \SC(\R^d)'$. Or il semble, dans la définition de $\SC^k(\R^d)'$ ci-dessous, ce soit le produit tensoriel qu'on considère (et dont on a besoin), tandis que pour la classe de Schwartz ce soit le produit cartésien $\SC(\R^d)^k$ dont a besoin. Différence fondamentale : 
		%		$$
		%		u_1\otimes \ldots \otimes u_k : X^k \to \C^k, 
		%		(x_1,\ldots, x_k)\mapsto (u_1(x_1), \ldots, u_k(x_k))
		%		$$
		%		tandis que 
		%		$$
		%		u_1\times \ldots \times u_k : X\to \C^k, 
		%		x \mapsto (u_1(x),\ldots, u_k(x))
		%		$$
		%	}
	
	Similarly to the space $\CC^\infty_c(\R^d)$, we endow $\SC(\R^d)$ with a topology: a sequence $(\psi_k)\subset\SC(\R^d)$ converges to $\psi\in\SC(\R^d)$ in the space $\SC(\R^d)$ if 
	$$
	\lim_{k\to\infty} 
	\sup_{x\in \R^d} |(1+\|x\|)^m \partial^\alpha(\psi_k-\psi)(x)|
	= 
	0
	,\spa \forall\ m\in\N,\ \forall\ \alpha\in\N^d
	.
	$$
	%as $k\to\infty$ for any $m\in\N$ and $\alpha\in\N^d$.
	
	\begin{Defin}\label{DefinTemperedDistributions}
		A tempered distribution is a linear map $
		T:\SC(\R^d)\to\C,\ \psi\mapsto \ps{T}{\psi}_{\SC',\SC}
		$
		which is continuous with respect to the topology of $\SC(\R^d)$. The set of tempered distributions is denoted $\SC(\R^d)'$.
	\end{Defin}
	
	For any integer $k\geq 1$, we denote by $\SC^k(\R^d)$ the \emph{cartesian} product of $\SC(\R^d)$, \ie the collection of vector fields $\Psi=(\psi_1,\ldots, \psi_k):\R^d\to\C^k$ such that $\psi_i\in\SC(\R^d)$ for each $i=1,2\ldots, k$. We denote by
	$\SC^k(\R^d)'$ the \emph{tensor} product of $\SC(\R^d)'$, \ie the collection of maps $T=(T_1,\ldots, T_k):\SC^k(\R^d)\to\C^k$ such that $T_i\in\SC(\R^d)'$ for each $i = 1,2,\ldots, k$ and defined, for all $\Psi=(\psi_1,\ldots,\psi_k)\in\SC^k(\R^d)$, by 
	$$
	\ps{T}{\Psi}_{(\SC^k)', \SC^k}
	:=
	(\ps{T_1}{\psi_1}_{\SC',\SC},\ldots, \ps{T_k}{\psi_k}_{\SC', \SC})
	.
	$$
	We let the operations described previously act on $\SC^k(\R^d)'$ componentwise.
	%	In particular, the identities we stated remain valid on $\SC^k(\R^d)'$.\\
	
	It is easy to see that if a sequence $(\varphi_k)\subset\CC^\infty_c(\R^d)$ converges to some $\varphi\in\CC^\infty_c(\R^d)$ in the topology $\CC^\infty_c(\R^d)$, then the convergence also holds in the topology of $\SC(\R^d)'$. In particular, any tempered distribution is a distribution on $\R^d$. 
	
	\begin{Examples}
		We list here typical examples of tempered distributions and a few usual ways to obtain tempered distributions from other tempered distributions.
		
		\begin{enumerate}
			\item If $f:\R^d\to\C$ is a measurable map such that the map
			$
			x\mapsto
			f(x) (1+\|x\|)^{-m}
			%	\frac{f(x)}{(1+\|x\|)^m}
			%	\in L^p(\R^d)
			$
			belongs to $L^p(\R^d)$, for some $m\in\N$ and $p\in [1,\infty)$, then $f\in\SC(\R^d)'$ by letting
			$$
			\ps{f}{\psi}_{\SC',\SC}
			:=
			\int_{\R^d} f(x)\psi(x)\, dx 
			,\spa 
			\forall\ \psi\in\SC(\R^d)
			.
			$$
			%	for any $\psi\in\SC(\R^d)$.
			
			\item Similarly, if $\mu$ is a Borel measure on $\R^d$ such that 
			$$
			\int_{\R^d}
			\frac{1}{(1+\|x\|)^m}\, d\mu(x)
			<
			\infty 
			$$
			for some $m\in\N$, then $\mu\in\SC(\R^d)'$.
			
			\item When $T\in\SC(\R^d)'$ and $f\in\CC^\infty(\R^d)$, we have already seen that the product $fT$ is a distribution on $\R^d$. For $fT$ to be tempered, it is further required that $f\psi$ be in the Schwartz class for all $\psi\in\SC(\R^d)$. This will be the case, for instance, if $f$ and all its derivatives have at most polynomial growth at infinity. If we impose no condition on the growth of $f$ and its derivatives, the product $fT$ might not be tempered; consider, e.g.,  $T\equiv 1$ and $f(x)=e^x$.
			%	\item Let $\psi\in\SC(\R^d)$ and $f\in\CC^\infty(\R^d)$. If $\partial^\alpha f$ has at most polynomial growth at infinity for any $\alpha \in\N^d$, then $f \psi\in\SC(\R^d)$. In particular, if $T\in\SC(\R^d)'$ is a tempered distribution, the product $fT$ defines a tempered distribution of $\SC(\R^d)'$. Notice that the product $f T$ always defines a distribution of $\D(\R^d)'$ if $T\in\D(\R^d)'$ and $f\in\CC^\infty(\R^d)$, but not necessarily a distribution of $\SC(\R^d)'$ if one does not impose restrictions on the growth of the derivatives of $f$ at infinity, e.g. $T\equiv 1$ and $f(x)=e^x$.
			\item If $T\in\SC(\R^d)'$ is a tempered distribution, then  $\partial^\alpha T$ is also a tempered distribution for all $\alpha\in\N^d$.
			
			\item If $T\in\SC(\R^d)'$, we define its Fourier transform $\FT(T)$ by letting 
			$$
			\ps{\FT(T)}{\psi}_{\SC',\SC}
			:=
			\ps{T}{\FT(\psi)}_{\SC',\SC}
			,\spa 
			\forall\ \psi\in\SC(\R^d)
			.
			$$
			Because the Schwartz class is closed under Fourier transforms, then $\FT(T)$ is a tempered distribution. In addition, the identities
			$$
			\FT(\partial^\alpha T) = (2i\pi \xi)^\alpha \FT(T),
			\quad 
			\textrm{ and } 
			\quad 
			\partial^\alpha \FT(T) = \FT((-2i\pi x)^\alpha T)
			$$
			hold in $\SC(\R^d)'$ for all $\alpha\in\N^d$, where $\xi^\alpha \FT(T)$ designates the multiplication of the distribution $\FT(T)$ by the smooth function $\xi^\alpha$ which has polynomial growth.
			%	\item Derivation preserves convergence in $\SC(\R^d)$ : if $(\psi_k)\subset\SC(\R^d)$ and $\psi\in\SC(\R^d)$, and if $(\psi_k)\to \psi$ in $\SC(\R^d)$ (\ie in the sense of the convergence in $\SC(\R^d)$), then $\partial^\alpha \psi_k \to \partial^\alpha \psi$ in $\SC(\R^d)$ for any $\alpha\in\N^d$. It follows that all distributional derivatives $\partial^\alpha T\in \D(\R^d)'$ are in fact tempered distributions of $\SC(\R^d)'$.
			
			%	\item Since the Fourier Transform maps $\SC(\R^d)$ to $\SC(\R^d)$, we can define the Fourier transform $\FT T$ of a tempered distribution $T\in\SC(\R^d)'$ by 
			%	$$
			%	\ps{\FT{T}}{\psi}
			%	:=
			%	\ps{T}{\wwidehat{\psi}}
			%	$$
			%	for any $\psi\in\SC(\R^d)$.
			%	\item Convolution preserves convergence in $\SC(\R^d)$ : if $(\psi_k)\subset\SC(\R^d)$ and $\psi\in\SC(\R^d)$, and if $(\psi_k)\to \psi$ in $\SC(\R^d)$, then $\psi_k*\eta \to \psi*\eta$ in $\SC(\R^d)$ for any $\eta\in\SC(\R^d)$.
		\end{enumerate}
	\end{Examples}
	
	%This allows one to define the convolution of a distribution with a function from $\SC(\R^d)$. 
	%
	%\begin{Defin}\label{DefinConvolutionDistributions}
	%	Let $d\geq 1$ and $T\in\SC(\R^d)'$. For any $\psi\in \SC(\R^d)$, we let $T*\psi:\SC(\R^d)\to\C$ be the tempered distribution defined by
	%	$$
	%	\ps{T*\psi}{\varphi}
	%	=
	%	\ps{T}{\widetilde{\psi}*\varphi}
	%	,
	%	$$
	%	for any $\varphi\in\SC(\R^d)$, where we let $\widetilde{\psi}(x):=\psi(-x)$.
	%\end{Defin}
	
	%	\textcolor{red}{Add a precise reference for this statement (probably Rudin).}
	
	When $T\in\SC(\R^d)'$ is a tempered distribution, we can define its convolution $T*\psi:\R^d\to\C$ with a test function $\psi\in\SC(\R^d)$ from the Schwartz class, by letting
	$$
	(T*\psi)(x)
	=
	\ps{T}{\psi(x-\cdot)}_{\SC',\SC}
	,\spa 
	\forall\ x\in\R^d 
	.
	$$
	The next proposition is given in Theorem 7.19 of \cite{RudinFunctional}.
	
	\begin{Prop}\label{PropConvolution}
		Let $T\in\SC(\R^d)'$ and $\psi\in \SC(\R^d)$. The map $T*\psi$ belongs to $\CC^\infty(\R^d)$, and we have 
		$$
		\partial^\alpha (T*\psi) 
		= 
		(\partial^\alpha T)*\psi 
		= 
		T * (\partial^\alpha \psi)
		,\spa \forall\ \alpha\in\N^d,
		$$
		over $\R^d$. In addition, $T*\psi$ has polynomial growth; in particular, $T*\psi$ is a tempered distribution, and the equality
		$$
		\FT(T*\psi)=\FT(T)\ \FT(\psi)
		$$
		holds in $\SC(\R^d)'$. 	
	\end{Prop}

	%\begin{Prop}\label{PropConvolution}
	%	Let $T\in\SC(\R^d)'$ and $\psi\in\SC(\R^d)$. The following equalities hold in $\SC(\R^d)'$ :
	%	\begin{enumerate}
		%		\item $\FT(T*\psi)=\FT(T)\ \wwidehat{\psi}$ ;
		%		\item $\partial^\alpha (T*\psi) = (\partial^\alpha T)*\psi = T * (\partial^\alpha \psi)$, for any $\alpha\in\N^d$ ;
		%		\item $\FT(\partial^\alpha T) = (2i\pi \xi)^\alpha \FT(T)$, for any $\alpha\in\N^d$ ;
		%		\item $\partial^\alpha \FT(T) = \FT((-2i\pi x)^\alpha T)$, for any $\alpha\in\N^d$.
		%	\end{enumerate}
	%\end{Prop}

	\subsection{Sobolev spaces}\label{sec:Sobolev}
	
	The main reference we used for this section is \cite{Evans1998}. Let $\Om\subset\R^d$ be an open subset. A function $u\in L^1_\loc(\Om)$ has weak derivatives of order $k$ in $\Om$ if, for any $\alpha\in\N^d$ with $|\alpha|\leq k$, the distributional derivative $\partial^\alpha u$ is in fact a function and belongs to $L^1_\loc(\Om)$, \ie there exists $v_\alpha\in L^1_\loc(\Om)$ such that 
	$$
	\int_\Om u(x)\partial^\alpha \varphi(x)\, dx 
	=
	(-1)^{|\alpha|}
	\int_\Om v_\alpha(x)\varphi(x)\, dx
	,\spa \forall\ \varphi\in\CC^\infty_c(\Om)
	.
	$$
	We write $v_\alpha=\partial^\alpha u$. When $u\in\CC^k(\Om)$, then $\partial^\alpha u$ coincides with the usual partial derivative of $u$. In the sequel, when no regularity of $u\in L^1_\loc(\Om)$ is assumed, then $\partial^\alpha u$ will always stand for a distributional derivative.
	
	For any integer $k\geq 1$, we define the Sobolev space
	$$
	H^k(\Om)
	=
	\Big\{u\in L^2(\Om) : \partial^\alpha u \in L^2(\Om), \ \forall \alpha\in\N^d, |\alpha|\leq k \Big\}
	.
	$$
	Since $\FT(\partial^\alpha u)=(2i\pi x)^\alpha \FT u$ holds in $\SC(\R^d)'$, and because a function belongs to $L^2(\R^d)$ exactly when its distributional Fourier transform does, the condition ``$u\in L^2(\R^d)$ and $\partial^\alpha u \in L^2(\R^d)$'' is equivalent to ``$u\in L^2(\R^d)$ and $x^\alpha \FT u\in L^2(\R^d)$''. It follows that we can equivalently define $H^k(\R^d)$ as
	$$
	H^k(\R^d)
	=
	\Big\{u\in L^2(\R^d) : \int_{\R^d} (1+\|\xi\|^{2k}) |(\FT u)(\xi)|^2\, d\xi < \infty \Big\}
	.
	$$
	For any real $s\geq 0$, we finally let
	\begin{align*}
		H^s(\R^d)
		&=
		\Big\{u\in L^2(\R^d) : \int_{\R^d} (1+\|\xi\|^{2s}) |(\FT u)(\xi)|^2\, d\xi < \infty \Big\}
		\\[2mm]
		&=
		\Big\{u\in L^2(\R^d) :  (1+\|\xi\|^2)^{s/2} (\FT u)(\xi)\in L^2(\R^d) \Big\}
		.
	\end{align*}
	It is easy to see that $H^{s'}(\R^d)\subset H^{s}(\R^d)$ if $s<s'$, and $H^0(\R^d)= L^2(\R^d)$ by Plancherel's formula. For any $s>0$, the space $H^s(\R^d)$ is a Hilbert space endowed with the inner product
	$$
	(u,v)_{H^s(\R^d)}
	=
	\int_{\R^d} (1+\|\xi\|^2)^s (\FT u)(\xi) \overline{(\FT v)(\xi)}\, d\xi 
	,
	$$
	where $\overline{(\FT v)(\xi)}$ here denotes the complex conjugate of $(\FT v)(\xi)$.	Sobolev spaces are particularly useful to study the regularity of distributional solutions $u$ to the Laplace equation $-\Delta u = f$, and will be used in Appendix \ref{Appendix:PDEGeneral} to prove Theorem \ref{TheorRankRegularityOdd}.
	
	One can also define Sobolev spaces $H^{-s}(\R^d)$ of negative order, with $s>0$, following the same construction than that of $H^s(\R^d)$. These spaces are not spaces of functions but spaces of distributions. For any real $s\geq0$, we let 
	$$
	H^{-s}(\R^d)
	=
	\Big\{ 
	u\in\SC(\R^d)' :  \xi\mapsto \frac{(\F u)}{(1+\|\xi\|^2)^{s/2}}\in L^2(\R^d)
	\Big\}
	.
	$$
	It follows from the definition that any $u\in H^{-s}(\R^d)$ is such that $\F u\in L^2_\loc(\R^d)$. We see that any $u\in L^2(\R^d)$ belongs to $H^{-s}(\R^d)$ for all $s\geq0$, and that any $u\in\SC(\R^d)'$ such that $\F u\in L^\infty(\R^d)$ belongs to $H^{-s}(\R^d)$ as soon as $s>d/2$. In particular, any probability measure $P$ belongs to $H^{-s}(\R^d)$ for all $s>d/2$ since $\F(P)$ is nothing but a rescaling of the usual characteristic function of $P$ and is thus bounded. Sobolev spaces of negative order will be briefly used in the proof of Proposition \ref{PropLocalEven}.

	\section{Proofs for Section \ref{sec:IntroLaplace}}\label{Appendix:Laplacian}

	\subsection{Proof of Proposition \ref{PropRegFracLaplace}}

	\begin{Proof}{Proposition \ref{PropRegFracLaplace}}
		Because $u\in\SC(\R^d)$, the map $\xi\mapsto (1+\|\xi\|)^m \|\xi\|^{2s} (\FT u)(\xi)$ is integrable on $\R^d$ for any $m\geq 0$. This implies that $\FT^{-1}(\|\xi\|^{2s} \FT u)\in\CC^m(\R^d)$ for all $m\geq 0$, so that $(-\Delta)^s u\in\CC^\infty(\R^d)$. It remains to prove that the inequality stated in Proposition \ref{PropRegFracLaplace} holds for any $\alpha\in\N^d$. We already observed in Section \ref{sec:IntroLaplace} of the main paper that
		$$
		\partial^\alpha (-\Delta)^s u
		=
		(-\Delta)^\sigma \Big((-\Delta)^n (\partial^\alpha u)\Big)
		,
		$$
		%		%	Let us first observe that $(-\Delta)^s u$ is continuous over $\R^d$. Indeed, the map $\xi\mapsto \|\xi\|^{2s} \wwidehat{u}(\xi)$ is integrable over $L^1(\R^d)$ since $u\in\SC(\R^d)$, so that its inverse Fourier transform, \ie $(-\Delta)^s u$, is continuous (and bounded) over $\R^d$. In particular, $(-\Delta)^s u$ is a distribution over $\R^d$.
		%		%	Since we always have
		%		\begin{align*}
			%			%		\lefteqn{
				%				%			\hspace{-10mm}
				%				\partial^\alpha (-\Delta)^s u
				%				&=
				%				\partial^\alpha \FT^{-1} (\|\xi\|^{2s} \FF u)
				%				%		}
			%			\\[2mm]
			%			&
			%			=
			%			\FT^{-1}\Big((2i \pi \xi)^\alpha \|\xi\|^{2s} \FT u\Big)
			%			\\[2mm]
			%			&
			%			=
			%			\FT^{-1}\Big( \|\xi\|^{2s} \FT( \partial^\alpha u) \Big)
			%			\\[2mm]
			%			&
			%			=
			%			(-\Delta)^s(\partial^\alpha u)
			%			=
			%			(-\Delta)^\sigma \Big((-\Delta)^n (\partial^\alpha u)\Big)
			%		\end{align*}
		for any $\alpha\in\N^d$. Since $(-\Delta)^n(\partial^\alpha u)\in\SC(\R^d)$ for all $\alpha\in\N^d$, it is enough to show that 
		\begin{equation}\label{eq:PropRegFracLaplace5}
			\sup_{x\in\R^d} |(1+\|x\|^{d+2\sigma})((-\Delta)^\sigma u)(x)|
			\lesssim 
			\|u\|_{L^1(\R^d)}
			+
			\sup_{z\in\R^d} \Big( (1+\|z\|)^{d+2} \|\nabla^2 u(z)\| \Big)
			.
		\end{equation}
		%	\begin{eqnarray}\label{eq:PropRegFracLaplace5}
			%		\lefteqn{ 
				%%			\hspace{-10mm}
				%			\sup_{x\in\R^d} |(1+\|x\|^{d+2\sigma})(-\Delta)^\sigma u(x)|
				%			\lesssim 
				%			|u|_{L^1(\R^d)}
				%			+
				%			|\nabla^2 u|_{L^\infty(\R^d)}
				%		}
			%%		\\[2mm]
			%		\\[-2mm]
			%		&&
			%		\hspace{60mm}
			%		+
			%		\sup_{z\in\R^d} \Big( (1+\|z\|)^{d+2} |\nabla^2 u(z)| \Big)\nonumber
			%		.
			%	\end{eqnarray}
		%	Since $(-\Delta)^\sigma u$ is the Fourier transform of an integrable function, it is bounded over $\R^d$. Therefore, it remains to show that 
		%	$$
		%	\sup_{x\in\R^d} \Big(\|x\|^{d+2\sigma}|(-\Delta)^\sigma u(x)|\Big)
		%	<
		%	\infty
		%	.
		%	$$
		By simple changes of variable, it is easy to show that 
		$$
		((-\Delta)^\sigma u)(x)
		=
		-\frac12 c_{d,\sigma}\int_{\R^d} \frac{u(x+y)+u(x-y)-2u(x)}{\|y\|^{d+2\sigma}}\, dy
		$$
		for all $x\in\R^d$; see, e.g., Lemma 3.2 in \cite{HitchhikerLaplace}. Thie last integral is not singular at $y=0$ anymore. Indeed, one can easily show that 
		\begin{equation}\label{eq:PropRegFracLaplace2}
			\frac{u(x+y)+u(x-y)-2u(x)}{\|y\|^{d+2\sigma}}
			=
			\frac{1}{\|y\|^{d+2\sigma}} 
			\int_{-1}^1 \ps{y}{(\nabla^2 u)(x+ty) y}\, dt
			.
		\end{equation}
		%	where $\nabla^2 u(z)$ stands for the Hessian matrix of $u$ at $z$. Let $|\nabla^2 u(z)|$ stand for the operator norm of $\nabla^2 u(z)$ and 
		%	$$
		%	|\nabla^2 u|_{L^\infty(\R^d)}
		%	:=
		%	\sup_{z\in\R^d} |\nabla^2 u(z)|
		%	.
		%	$$
		The r.h.s. of \eqref{eq:PropRegFracLaplace2} is then bounded by 
		$
		\|\nabla^2 u\|_{L^\infty(\R^d)}/\|y\|^{d+2\sigma-2}
		$
		,
		which is integrable near the origin since $\sigma<1$. 
		
		Let us first show that 
		\begin{equation}\label{eq:PropRegFracLaplace3}
			\sup_{x\in\R^d} |((-\Delta)^\sigma u)(x)| 
			\lesssim 
			\|u\|_{L^1(\R^d)} + \sup_{z\in\R^d} \Big( (1+\|z\|)^{d+2} \|\nabla^2 u(z)\| \Big)
			.
		\end{equation}
		Fix $x\in\R^d$ and write 
		\begin{eqnarray*}
			\lefteqn{
				\hspace{-10mm}
				-\frac{2}{c_{d,\sigma}} ((-\Delta)^\sigma u)(x)
				=
				\int_{\R^d\sm B_1} \frac{u(x+y)+u(x-y)-2u(x)}{\|y\|^{d+2\sigma}}\, dy
			}
			\\[2mm]
			&&
			\hspace{30mm}
			+
			\int_{B_1} \frac{u(x+y)+u(x-y)-2u(x)}{\|y\|^{d+2\sigma}}\, dy
			\\[2mm]
			&&
			\hspace{17mm}
			=: I_1(x) + I_2(x)
			.
		\end{eqnarray*}
		For $I_1$, we have $|I_1(x)|\leq 4 \|u\|_{L^1(\R^d)}$. For $I_2$, recalling (\ref{eq:PropRegFracLaplace2}), we have
		\begin{align*}
			%		\lefteqn{
				|I_2(x)|
				&\leq 
				\int_{B_1} \frac{1}{\|y\|^{d+2\sigma-2}} \int_{-1}^1 \|\nabla^2 u(x+ty)\|\, dt\, dy
				\\
				%		}
			&\lesssim
			\sup_{z\in\R^d} \|\nabla^2 u(z)\|
			\\
			&\leq 
			\sup_{z\in\R^d} \Big( (1+\|z\|)^{d+2} \|\nabla^2 u(z)\| \Big)
			,
		\end{align*}
		which yields \eqref{eq:PropRegFracLaplace3}. Let us now show that 
		\begin{equation}\label{eq:PropRegFracLaplace4}
			\sup_{x\in\R^d} \Big( \|x\|^{d+2\sigma} |((-\Delta)^\sigma u)(x)| \Big)
			\lesssim 
			\|u\|_{L^1(\R^d)}
			+
			\sup_{z\in\R^d} (1+\|z\|)^{d+2} \|\nabla^2 u(z)\| 
			.
		\end{equation}
		Fix $x\in\R^d$ and write 
		\begin{eqnarray*}
			\lefteqn{
				\hspace{-10mm}
				-\frac{2}{c_{d,\sigma}} \|x\|^{d+2\sigma} ((-\Delta)^\sigma u)(x)
				=
				\|x\|^{d+2\sigma}\int_{\R^d \sm {B_{\frac12 \|x\|}}} \frac{u(x+y)+u(x-y)-2u(x)}{\|y\|^{d+2\sigma}}\, dy
			}
			\\[2mm]
			&&
			\hspace{35mm}
			+
			\|x\|^{d+2\sigma}\int_{{B_{\frac12 \|x\|}}} \frac{u(x+y)+u(x-y)-2u(x)}{\|y\|^{d+2\sigma}}\, dy
			\\[2mm]
			&&
			\hspace{29mm}
			=: J_1(x) + J_2(x)
			.
		\end{eqnarray*}
		For $J_1$, we have $|J_1(x)|\leq 4 \|u\|_{L^1(\R^d)} 2^{d+2\sigma}$.
		%			.
		%	\int_{\|y\|>\frac12\|x\|} \frac{1}{\|y\|^{d+2\sigma}}\, dy
		%	\\[2mm]
		%	&\leq 
		%	4 |u|_{L^1(\R^d)} 2^{d+2\sigma} \int_{\|y\|>\frac12} \frac{1}{\|y\|^{d+2\sigma}}\, dy
		%	\\[2mm]
		%	&
		%	=: C_1 < \infty 
		%		\end{align*}
	For $J_2$, recalling (\ref{eq:PropRegFracLaplace2}), we have
	\begin{align*}
		%		\lefteqn{
			|J_2(x)|
			&\leq 
			\|x\|^{d+2\sigma}\int_{{B_{\frac12 \|x\|}}} \frac{1}{\|y\|^{d+2\sigma-2}} \int_{-1}^1 \|\nabla^2 u(x+ty)\|\, dt\, dy
			%		}
		\\[2mm]
		&=
		\int_{{B_{\frac12 \|x\|}}} 
		\frac{1}{\|y\|^{d+2\sigma-2}}
		\int_{-1}^1 
		\bigg(\frac{\|x\|}{\|x+ty\|}\bigg)^{d+2\sigma} 
		\|x+ty\|^{d+2\sigma} 
		\|\nabla^2 u(x+ty)\|\, 
		dt\, dy
		.
	\end{align*}
	Now observe that, for any $t\in [-1,1]$ and $y$ with $\|y\|\leq \frac12 \|x\|$, we have
	$$
	\frac{\|x\|}{\|x+ty\|}
	\leq 
	\frac{\|x\|}{\|x\|-|t|\|y\|}
	\leq
	2
	.
	$$
	Letting
	$$
	C_k(u)
	:=
	\sup_{x\in\R^d} (1+\|z\|)^k \|\nabla^2 u(z)\|
	<
	\infty
	,
	$$
	for any $k$, we thus have
	$$
	\|x+ty\|^{d+2\sigma} \|\nabla^2 u(x+ty)\|
	\leq 
	\frac{C_{N+d+2\sigma}(u)}{(1+\|x+ty\|)^N}
	\leq 
	\frac{C_{N+d+2\sigma}(u)}{(1+\frac12 \|x\|)^N}
	$$
	for all $N$, $\|y\|\leq \frac12 \|x\|$, and  $t\in [-1,1]$. Then fix $N=2-2\sigma$. We have
	\begin{align*}
		|J_2(x)|
		&\leq  
		2^{d+2\sigma} C_{d+2}(u) \frac{1}{(1+\frac12 \|x\|)^{2-2\sigma}} \int_{{B_{\frac12 \|x\|}}} \frac{1}{\|y\|^{d+2\sigma-2}}\, dy 
		%	\\[2mm]
		%	& 
		%	\frac{\|x\|^{2-2s}}{(1+\frac12 \|x\|)^{2-2s}}
		%	\leq 
		%	2^{2-2s}
		.
	\end{align*}
	In addition, since $\sigma<1$ straightforward computations provide
	$$
	\int_{{B_{\frac12 \|x\|}}} \frac{1}{\|y\|^{d+2\sigma-2}}\, dy 
	\lesssim 
	\|x\|^{2-2\sigma}
	.
	$$
	It follows that
	$$
	\sup_{x\in\R^d} |J_2(x)|
	\lesssim 
	\sup_{x\in\R^d} (1+\|z\|)^{d+2} \|\nabla^2 u(z)\|
	.
	$$
	We deduce that 
	$$
	\sup_{x\in\R^d } \Big( \|x\|^{d+2\sigma} |(-\Delta)^\sigma u(x)| \Big)
	\lesssim 
	\|u\|_{L^1(\R^d)} + \sup_{x\in\R^d} (1+\|z\|)^{d+2} \|\nabla^2 u(z)\| 
	,
	$$
	which establishes \eqref{eq:PropRegFracLaplace4}. Putting \eqref{eq:PropRegFracLaplace3} and \eqref{eq:PropRegFracLaplace4} together yields \eqref{eq:PropRegFracLaplace5}, which concludes the proof.
\end{Proof}

\subsection{Proof of Proposition \ref{PropLaplaceExplicit}}

\begin{Proof}{Proposition \ref{PropLaplaceExplicit}}	

(i) This is Proposition 2.4 in \cite{Silvestre2007}.\\

(ii) Because $u\in L^2(\R^d)$, we have $u\in \SC_s(\R^d)'$. In particular, $(-\Delta)^s u$ is a well-defined tempered distribution. Fix $\psi\in\SC(\R^d)$. We have
$$
\ps{(-\Delta)^s u}{\psi}_{\SC',\SC}
=
\ps{u}{(-\Delta)^s \psi}_{\SC',\SC}
=
(2\pi)^{2s} \int_{\R^d} u(x) \FT^{-1}(\|\xi\|^{2s} \FT \psi)(x)\, dx
.
$$
%	\begin{eqnarray*}
	%		\lefteqn{
		%			\ps{(-\Delta)^s u}{\psi}
		%			=
		%			\ps{u}{(-\Delta)^s \psi}
		%			=
		%			\int_{\R^d} u(x) \FT^{-1}(\|\xi\|^{2s} \FT \psi(\xi))(x)\, dx
		%		}
	%		\\[2mm]
	%		&&
	%	\end{eqnarray*}
Because $u\in L^2(\R^d)$ and $\|\xi\|^{2s} \FT \psi\in L^2(\R^d)$, we have 
$$
\ps{(-\Delta)^s u}{\psi}_{\SC',\SC}
=
(2\pi)^{2s} \int_{\R^d} (\FT^{-1} u)(\xi)\ \|\xi\|^{2s} (\FT \psi)(\xi)\, d\xi
$$
by interchanging the inverse Fourier transform under the integral. Since $u\in H^{2s}(\R^d)$, we have $\|\xi\|^{2s} \FT^{-1}u\in L^2(\R^d)$. Since $\psi\in L^2(\R^d)$, interchanging the Fourier transform again yields
$$
\ps{(-\Delta)^s u}{\psi}_{\SC',\SC}
=
(2\pi)^{2s} \int_{\R^d} \FT\big(\|\xi\|^{2s}(\FT^{-1} u)(\xi)\big)(x)\  \psi(x)\, dx
.
$$
Observing that 
$$
\FT\big(\|\xi\|^{2s}(\FT^{-1} u)(\xi)\big)
=
\FT^{-1}\big(\|\xi\|^{2s}(\FT u)(\xi)\big)
$$
yields
$$
\ps{(-\Delta)^s u}{\psi}_{\SC',\SC}
=
\int_{\R^d} (2\pi)^{2s} \FT^{-1}(\|\xi\|^{2s}(\FT u))(x) \psi(x)\, dx
.
$$
The conclusion follows since $\psi\in\SC(\R^d)$ was arbitrary.\\

(iii) First observe that if $s$ is an integer, then the conclusion is straightforward as there is no differentiability issue at the origin anymore; the map $\|\xi\|^{2s}$ is smooth, so that the product $\|\xi\|^{2s}\FF u$ is a tempered distribution and the product $\|\xi\|^{2s} \psi$ is a Schwartz function for any $\psi\in\SC(\R^d)$. Applying the usual definitions of distributional Fourier transform and multiplication by a smooth function with at most polynomial growth---see Appendix \ref{sec:Distributions}---yields the desired conclusion. Let us, therefore, assume that $s\in (0,\infty)\sm \N$. Since $u\in\SC_s(\R^d)'$, then $(-\Delta)^s u$ is a tempered distribution; see Section \ref{sec:IntroLaplace} in the main paper. In addition, because $\FF u \in L^1_\loc(\R^d)\cap \SC_s(\R^d)'$, we have 
$$
\int_{\R^d} 
\frac{\|\xi\|^{2s} |(\FF u)(\xi)|}{(1+\|\xi\|)^{d+4s}}\, d\xi
\leq 
\int_{\R^d} 
\frac{|(\FF u)(\xi)|}{(1+\|\xi\|)^{d+2s}}\, d\xi
\leq 
\int_{\R^d} 
\frac{|(\FF u)(\xi)|}{(1+\|\xi\|)^{d+2(s-\lfloor s\rfloor)}}\, d\xi
<
\infty
.
$$
This implies that $\|\xi\|^{2s} \FF u$ is a tempered distribution as well; see Appendix \ref{sec:Distributions}. In particular, we will equivalently show that
$$
\FF\big((-\Delta)^s u\big)(\xi) 
=
(2\pi)^{2s} \|\xi\|^{2s} (\FF u)(\xi) 
$$
holds in $\SC(\R^d)'$. In fact, because $\CC^\infty_c(\R^d)$ is dense in $\SC(\R^d)$ with respect to the topology of $\SC(\R^d)$, it is enough to show that the previous equality holds in $\D(\R^d)'$; see Appendix \ref{sec:Distributions}. Fix $\varphi\in\CC^\infty_c(\R^d)$. Because $u$ and $\FF u$ belong to $L^1_\loc(\R^d)$, we need to establish that 
$$
\int_{\R^d} u(x)\ \big((-\Delta)^s(\FF \varphi)\big)(x)\, dx 
=
(2\pi)^s\int_{\R^d} \|\xi\|^{2s} (\FF u)(\xi)\ \varphi(\xi)\, d\xi 
.
$$
Because $\big((-\Delta)^s(\FF \varphi)\big)(x)=(2\pi)^{2s}\FF\big(\|\xi\|^{2s} \varphi(\xi)\big)(x)$, we then need to show that
$$
\int_{\R^d} u(x)\ \FF\big(\|\xi\|^{2s} \varphi(\xi)\big)(x)\, dx 
=
\int_{\R^d} (\FF u)(\xi)\ \|\xi\|^{2s} \varphi(\xi)\, d\xi 
.
$$
This is certainly true when $u$ or $\FF u$ belongs to $L^1(\R^d)$ by "swapping" the Fourier transform under the integral, or if $s=0$ by definition of the distributional Fourier transform. Because we are not assuming that either $u$ or $\FF u$ belongs to $L^1(\R^d)$, we need to show that the distributional identity $\ps{\FF u}{\varphi}_{\SC',\SC}=\ps{u}{\FF\varphi}_{\SC',\SC}$ remains valid when $\varphi(\xi)$ is replaced by $\|\xi\|^{2s}\varphi(\xi)$, which is no longer in $\CC^\infty_c(\R^d)$ or $\SC(\R^d)$ due to the singularity at the origin. Letting $w(\xi):=\|\xi\|^{2s}\varphi(\xi)$, we need to show that 
\begin{equation}\label{eq:FourierToShow}
	\int_{\R^d} u(x)\ (\FF w)(x)\, dx 
	=
	\int_{\R^d} (\FF u)(\xi)\ w(\xi)\, d\xi 
	.
\end{equation}
Notice that the product $u\ \FF w$ is integrable on $\R^d$. Indeed, we have 
$$
\int_{\R^d} |u(x)| |(\FF w)(x)|\, dx 
\leq 
\sup_{x\in\R^d} | (1+\|x\|^{d+2(s-\lfloor s\rfloor)}) (\FF w)(x)|
\times 
\int_{\R^d} \frac{|u(x)|}{1+\|x\|^{d+2(s-\lfloor s\rfloor)}}\, dx
,
$$
which is finite by Proposition \ref{PropRegFracLaplace}, since $\FF w= (-\Delta)^s \psi$ with $\psi = \FF \varphi\in \SC(\R^d)$, and because $u\in L^1_\loc(\R^d)\cap \SC_s(\R^d)'$---see (\ref{eq:L1locFrac}). 

Let $\rho\in\CC^\infty_c(\BB^d)$ be such that $0\leq \rho\leq 1$ and $\int_{\R^d} \rho(x)\, dx=1$. Then define, for all $k\in\N$, $\rho_k(x):=k^d \rho(kx)$ and the convolution
$
w_k(\xi) 
=
(w*\rho_k)(x)
.
$
Because $\rho_k(x)=0$ as soon as $\|x\|\geq 1/k$, and the support of $w_k$ is included in the Minkowski sum $\supp(w) + \BB^d_{1/k}$ of the supports of $w$ and $\rho_k$, we see that $w_k\in\CC^\infty_c(\R^d)$. In particular, (\ref{eq:FourierToShow}) holds for $w_k$ in place of $w$. We have $w_k(x)\to w(w)$ for all $x$ (the convergence is, in fact, uniform over compacts sets since $w$ is continuous) and $\FF w_k = \FF w\ \FF \rho_k \to \FF w$ since $(\FF \rho_k)(x) = (\FF \rho)(x/k)\to (\FF \rho)(0)=1$ as $k\to\infty$. We are going to apply Lebesgue's dominated convergence theorem to (\ref{eq:FourierToShow}) with $w_k$ to show that (\ref{eq:FourierToShow}) holds for $w$. To dominated $u\ \FF w_k$, it suffices to notice that 
$$
|u(x) (\FF w_k)(x)| 
\leq 
|u(x)\ (\FF w)(x)\ (\FF \rho)(x/k)|
\leq 
|u(x)\ (\FF w)(x)|
\in L^1(\R^d)
.
$$
For $(\FF u)\ w_k$, first fix a compact set $K\subset\R^d$ such that $\supp(w_k)\subset \supp(w)+\BB^d_{1/k}\subset K$ for all $k$. Then, notice that 
$$
\int_{\R^d} 
(\FF u)(\xi)\ w_k(\xi)\, d\xi 
=
\int_{K} 
(\FF u)(\xi)\ w_k(\xi)\, d\xi 
.
$$
It remains to observe that 
$
|w_k(x)|
\leq 
\sup_{x\in\R^d} |w(x)| 
$
and to recall that $\FF u$ is integrable over $K$. As a consequence, Lebesgue's dominated convergence applies to (\ref{eq:FourierToShow}) with $w_k$ and yields the conclusion.
\end{Proof}

\subsection{Proof of Corollary \ref{CorolFTRiesz}.}

The proof of Corollary \ref{CorolFTRiesz} requires the following Lemma.

\begin{Lem}\label{LemFourierTransformRiesz}
Let $d\geq 2$ be an integer and $\alpha\in (0,d)$ a real number. Then the Fourier transform of the tempered distribution $1/\|x\|^\alpha$ is given by 
$$
\FT\Big(\frac{1}{\|x\|^{\alpha}}\Big)(\xi)
=
\frac{\Gamma(\frac{d-\alpha}{2})}{\pi^{\frac{d}{2}-\alpha} \Gamma(\frac{\alpha}{2})} \frac{1}{\|\xi\|^{d-\alpha}}
.
$$
\end{Lem}

\begin{Proof}{Lemma \ref{LemFourierTransformRiesz}}
The proof follows the same lines as \cite{Landkof1972} and \cite{Bochner59}---see Equation $(1.1.1)$ in \cite{Landkof1972} and Theorem 56 in \cite{Bochner59}---but we establish the stated distributional equality in a more modern language. 

Let $g\in L^1(\R^d)$ be such that $g(x)=h(\|x\|)$ for all $x\in\R^d$. First, let us prove that 
\begin{equation}\label{Eq:FourierTransformGSpheric}
	\wwidehat{g}(\xi)
	=
	\frac{1}{(2\pi)^{\frac{d}{2}} \|\xi\|^d} \int_0^\infty r^{\frac{d}{2}} h\Big(\frac{r}{2\pi\|\xi\|}\Big) J_{\frac{d-2}{2}}(r)\, dr
	,\spa \forall\ \xi\in\R^d\sm\{0\}
	,
\end{equation}
where $J_\nu$ is the Bessel function of the first kind of order $\nu$. Fix $\xi\in\R^d\sm\{0\}$. Because $g(x)=h(\|x\|)$ for all $x$, we have
$$
(\FT g)(\xi)
=
\int_{\R^d} 
g(x)e^{-2i\pi\ps{x}{\xi}}\, dx 
=
\int_{\R^d} 
g(x)e^{-2i\pi\ps{x}{O\xi}}\, dx 
$$
for any $d\times d$ orthogonal matrix $O$. Then, assume that $\xi = \|\xi\|(1,0,\ldots, 0)$ and compute
\begin{align*}
	%		\lefteqn{
		\int_{\R^d} 
		g(x)e^{-2i\pi\ps{x}{\xi}}\, dx 
		&=
		\int_{\R^d} 
		h(\|x\|)e^{-2i\pi\|\xi\| x_1}\, dx
		%		}
	\\[2mm]
	&=
	\int_{\R} e^{-2i\pi\|\xi\| x_1} \bigg( \int_{\R^{d-1}} h\Big(\sqrt{x_1^2+\|y\|^2}\Big)\, dy\bigg)\, dx_1
	\\[2mm]
	&=
	S_{d-2} \int_{\R} e^{-2i\pi\|\xi\| x_1} \bigg( \int_0^\infty t^{d-2} h\Big(\sqrt{x_1^2+t^2}\Big)\, dt\bigg)\, dx_1
	,
\end{align*}
where $S_{d-2}= 2 \pi^{(d-1)/2}/\Gamma((d-1)/2)$ is the surface area of the $(d-2)$-dimensional unit sphere in $\R^{d-1}$. Write $(x_1,t)$ in spherical coordinates as $(x_1,t)=(r\cos\theta, r\sin\theta)$,	with $r\in [0,\infty)$ and $\theta\in [0,\pi]$, because $t>0$. This provides
\begin{align*}
	(\FT g)(\xi)
	&=
	S_{d-2} \int_0^\infty \bigg( \int_0^\pi e^{-2i\pi \|\xi\| r\cos\theta} (r\sin\theta)^{d-2}  h(r)\, r d\theta\bigg)\, dr
	\\[2mm]
	&=
	S_{d-2} \int_0^\infty r^{d-1} h(r) \bigg( \int_0^\pi e^{-2i\pi \|\xi\| r\cos\theta} (\sin\theta)^{d-2}  d\theta\bigg)\, dr
	.
\end{align*}
Now express $\int_0^\pi e^{-2i\pi \|\xi\| r\cos\theta} (\sin\theta)^{d-2}  d\theta$ in terms of Bessel functions: for all $\nu\in\C$ with $\Reel(\nu)>-1/2$, the Bessel function of the first kind of order $\nu$ satisfies, for all $x\in\R$,
\begin{align*}
	J_\nu(x)
	&=
	\frac{(\frac{x}{2})^\nu}{\sqrt{\pi}\ \Gamma(\nu+\frac{1}{2})} \int_{-1}^1 (1-t^2)^{\nu-\frac{1}{2}} \cos(xt)\, dt
	\\[2mm]
	&=
	\frac{(\frac{x}{2})^\nu}{\sqrt{\pi}\ \Gamma(\nu+\frac{1}{2})} \int_{-1}^1 (1-t^2)^{\nu-\frac{1}{2}} e^{-ixt}\, dt
	;
\end{align*}
see (10.9.4) in \cite{HandBookMathNist}. Substituting $t=\cos\theta$ leads to 
$$
J_\nu(x)
=
\frac{(\frac{x}{2})^\nu}{\sqrt{\pi}\ \Gamma(\nu+\frac{1}{2})}
\int_0^\pi (\sin \theta)^{2\nu} e^{-ix\cos\theta}\, d\theta 
,\spa \forall\ x\in\R
.
$$
It follows that
$$
\int_0^\pi e^{-2i\pi \|\xi\| r\cos\theta} (\sin\theta)^{d-2}  d\theta
=
\frac{\sqrt{\pi}\ \Gamma(\frac{d-1}{2})}{(\pi\|\xi\|r)^\frac{d-2}{2}}
J_{\frac{d-2}{2}}(2\pi\|\xi\|r)
,\spa \forall\ r>0
.
$$
We deduce that 
\begin{align*}
	%	\lefteqn{
		(\FT g)(\xi)
		&=
		\frac{2\pi}{\|\xi\|^{\frac{d-2}{2}}} \int_0^\infty r^{\frac{d}{2}} h(r) J_{\frac{d-2}{2}}(2\pi\|\xi\|r)\, dr
		%	}
	\\[2mm]
	%	\hspace{15mm}
	&=
	\frac{1}{(2\pi)^{\frac{d}{2}} \|\xi\|^d} \int_0^\infty r^{\frac{d}{2}} h\Big(\frac{r}{2\pi\|\xi\|}\Big) J_{\frac{d-2}{2}}(r)\, dr
	,\spa \forall\ \xi\in\R^d\sm\{0\}
	,
\end{align*}
which yields \eqref{Eq:FourierTransformGSpheric}. 

Now fix $\alpha \in (0,(d+1)/2)$. For any $k\in\N$, define 
%		$$
%		g_{k,\alpha}(x)
%		:=
%		\|x\|^{\alpha-n}\Ind{0<\|x\|<k}
%		,\spa \forall\ x\in\R^d
%		$$
%		and 
$$
h_{k,\alpha}(t)
=
t^{\alpha-n}\Ind{0<t<k}
,\spa \forall\ t>0
.
$$
and $g_{k,\alpha}(x):=h_{k,\alpha}(\|x\|)$ for all $x\in\R^d$. For any $k$, we have $g_{k,\alpha}\in L^1(\R^d)$. Then, applying \eqref{Eq:FourierTransformGSpheric} to $g_{k,\alpha}$ provides
$$
(\FT g_{k,\alpha})(\xi)
=
\frac{1}{(2\pi)^{\alpha-\frac{d}{2}} \|\xi\|^{\alpha}} \int_0^{2\pi\|\xi\|k} r^{\alpha-\frac{d}{2}} J_{\frac{d-2}{2}}(r)\, dr
,\spa \forall\ \xi\in\R^d
.
$$
According to (10.22.43) in \cite{HandBookMathNist}, this integral converges to 
$$
\int_0^\infty r^{\alpha-\frac{d}{2}} J_{\frac{d-2}{2}}(r)\, dr 
=
2^{\alpha-\frac{d}{2}} \frac{\Gamma(\frac{\alpha}{2})}{\Gamma(\frac{d-\alpha}{2})}
$$
as $k\to\infty$, since $\alpha \in (0,(d+1)/2)$. It follows that 
$$
\lim_{k\to\infty} (\FT g_{k,\alpha})(\xi)
=
\frac{\pi^{\frac{d}{2}-\alpha} \Gamma(\frac{\alpha}{2})}{\Gamma(\frac{d-\alpha}{2})}\ \frac{1}{\|\xi\|^\alpha}
,\spa \forall\ \xi\in\R^d\sm\{0\}
.
$$
Note that $|(\FT g_{k,\alpha})(\xi)|\lesssim 1/\|\xi\|^\alpha$ for all $\xi\in\R^d\sm\{0\}$, uniformly in $k$. Then, fix $\psi\in\SC(\R^d)$ and observe that we have, uniformly in $k$,
$$
|(\FT g_{k,\alpha})(\xi)\psi(\xi)|
\lesssim
\frac{|\psi(\xi)|}{\|\xi\|^\alpha}
,\spa \forall\ \xi\in\R^d\sm\{0\}
.
$$
Since $\alpha<d$ (recall that $\alpha<(d+1)/2$ and $d\geq 2$) and $\psi\in\SC(\R^d)$, then
$|\psi(\xi)|\|\xi\|^{-\alpha}\in L^1(\R^d)$. Lebesgue's dominated convergence theorem thus yields 
$$
\int_{\R^d} (\FT g_{k,\alpha})(\xi) \psi(\xi)\, d\xi 
\to 
\frac{\pi^{\frac{d}{2}-\alpha} \Gamma(\frac{\alpha}{2})}{\Gamma(\frac{d-\alpha}{2})}\ \int_{\R^d} \frac{1}{\|\xi\|^\alpha}\psi(\xi)\, d\xi
$$
as 
$k\to\infty$. Observe now that, because $g_{k,\alpha}\in L^1(\R^d)$, we have, for all $k$,
$$
\int_{\R^d} (\FT g_{k,\alpha})(\xi) \psi(\xi)\, d\xi 
=
\int_{\R^d} g_{k,\alpha}(\xi) (\FT \psi)(\xi)\, d\xi 
.
$$
Because $|g_{k,\alpha}(x)(\FT \psi)(x)|\leq \|x\|^{\alpha-n} |(\FT \psi)(x)|\in L^1(\R^d)$ holds uniformly in $k$, Lebesgue's dominated convergence theorem yields
$$
\int_{\R^d} g_{k,\alpha}(\xi) (\FT \psi)(\xi)\, d\xi 
\to 
\int_{\R^d} \frac{1}{\|x\|^{d-\alpha}} (\FT \psi)(\xi)\, d\xi 
$$
as $k\to\infty$. It follows that 
$$ 
\int_{\R^d} \frac{1}{\|x\|^{d-\alpha}} (\FT \psi)(\xi)\, d\xi 
=
\frac{\pi^{\frac{d}{2}-\alpha} \Gamma(\frac{\alpha}{2})}{\Gamma(\frac{d-\alpha}{2})}\ \int_{\R^d} \frac{1}{\|\xi\|^\alpha}\psi(\xi)\, d\xi
,\spa \forall\ \psi\in\SC(\R^d)
.
$$
We deduce that
\begin{equation}\label{Eq:FourierIntermediateEquality}
	\FT\Big(\frac{1}{\|x\|^{d-\alpha}}\Big)(\xi)
	=
	\frac{\pi^{\frac{d}{2}-\alpha} \Gamma(\frac{\alpha}{2})}{\Gamma(\frac{d-\alpha}{2})}
	\frac{1}{\|\xi\|^\alpha}
\end{equation}
in $\SC(\R^d)'$ for any $\alpha\in (0,(d+1/2))$. Taking the inverse Fourier transform on both sides of \eqref{Eq:FourierIntermediateEquality} yields
\begin{equation}\label{Eq:FourierIntermediateEquality2}
	\FT\Big(\frac{1}{\|x\|^{\alpha}}\Big)(\xi)
	=
	\frac{\Gamma(\frac{d-\alpha}{2})}{\pi^{\frac{d}{2}-\alpha} \Gamma(\frac{\alpha}{2})} \frac{1}{\|\xi\|^{d-\alpha}}
\end{equation}
in $\SC(\R^d)'$ for any $\alpha\in (0,(d+1)/2)$. Now fix $\beta\in ((d-1)/2,d)$ and write $\beta=d-\alpha$ for some $\alpha\in (0, (d+1)/2)$. Then, \eqref{Eq:FourierIntermediateEquality} yields
\begin{equation}\label{Eq:FourierIntermediateEquality3}
	\FT\Big(\frac{1}{\|x\|^{\beta}}\Big)(\xi)
	=
	\frac{\pi^{\frac{d}{2}-\alpha} \Gamma(\frac{\alpha}{2})}{\Gamma(\frac{d-\alpha}{2})}
	\frac{1}{\|\xi\|^\alpha}
	=
	\frac{\Gamma(\frac{d-\beta}{2})}{\pi^{\frac{d}{2}-\beta} \Gamma(\frac{\beta}{2})} \frac{1}{\|\xi\|^{d-\beta}}
\end{equation}
in $\SC(\R^d)'$. Putting \eqref{Eq:FourierIntermediateEquality2} and \eqref{Eq:FourierIntermediateEquality3} together yields the conclusion for any 
$$
\alpha\in \Big(0,\frac{d+1}{2}\Big)\cup\Big(\frac{d-1}{2},d\Big)=(0,d)
.
$$
This concludes the proof.
%	and compute 
%	$$
%	\ps{\FT(g_{k,\alpha})}{\psi}
%	=
%	\ps{g_{k,\alpha}}{\wwidehat{\psi}}
%	=
%	\int_{\R^d} g_{k,\alpha}(x) \wwidehat{\psi}(x)\, dx 
%	.
%	$$
%	The sequence of functions $(g_{k,\alpha})$ converges almost everywhere to $g_{\infty,\alpha}(x):=\|x\|^{\alpha-n}\Ind{x\neq 0}$. We further have that 
%	$$
%	|g_{k,\alpha}(x)\psi(x)|\lesssim \frac{\wwidehat{\psi}(x)}{\|x\|^\alpha}
%	$$
%	uniformly in $k$, with $\|x\|^{-\alpha} \wwidehat{\psi}(x)\in L^1(\R^d)$ since $\alpha<n$ (recall that $\alpha<\frac{d+1}{2}$ and that $d\geq 2$) and $\wwidehat{\psi}\in\SC(\R^d)$.
\end{Proof}

\vspace{3mm}

We now proceed with the proof of Corollary \ref{CorolFTRiesz}.

\begin{Proof}{Corollary \ref{CorolFTRiesz}.}
Because $\alpha\in (0,d)$, Lemma \ref{LemFourierTransformRiesz} entails that 
$$
\FT\Big(\frac{1}{\|x\|^{\alpha}}\Big)(\xi)
=
\frac{\Gamma(\frac{d-\alpha}{2})}{\pi^{\frac{d}{2}-\alpha} \Gamma(\frac{\alpha}{2})} \frac{1}{\|\xi\|^{d-\alpha}}
.
$$
Since $\alpha\in (0,d)$, we thus have that $1/\|x\|^\alpha$ and its Fourier transform belong to $L^1_\loc(\R^d)\cap \SC_s(\R^d)'$. Consequently, Proposition \ref{PropLaplaceExplicit} (iii) entails that 
$$
(-\Delta)^s\Big(\frac{1}{\|x\|^\alpha}\Big)
=
(2\pi)^{2s} 
\frac{\Gamma(\frac{d-\alpha}{2})}{\pi^{\frac{d}{2}-\alpha} \Gamma(\frac{\alpha}{2})}
\FF^{-1}\Big(\|\xi\|^{2s} \frac{1}{\|\xi\|^{d-\alpha}}\Big)
.
$$ 
If $\alpha = d-2s$, then the last Fourier transform reduces to $\FF^{-1}(1)=\delta$, Dirac's distribution at $0$, which yields 
$$
(-\Delta)^s\Big(\frac{1}{\|x\|^{d-2s}}\Big) 
=
4^{s}\pi^{d/2} 
\frac{\Gamma(s)}{\Gamma(\frac{d-2s}{2})}\
\delta 
.
$$
If $\alpha\in (0,d-2s)$, Lemma \ref{LemFourierTransformRiesz} entails that 
\begin{equation}\label{eq:FTRiesz}
	(-\Delta)^s\Big(\frac{1}{\|x\|^\alpha}\Big)	
	=
	4^s\
	\frac{
		\Gamma(\frac{d-\alpha}{2}) \Gamma(\frac{\alpha+2s}{2})
	}{
		\Gamma(\frac{\alpha}{2})\Gamma(\frac{d-\alpha-2s}{2})
	} \frac{1}{\|x\|^{\alpha+2s}}
	.
\end{equation}
This concludes the proof.
\end{Proof}

\section{Complements for Section \ref{sec:PDEGeneral}} \label{Appendix:PDEGeneral}

\subsection{The operators $\LL_d$ and $\LL_d^*$} \label{Appendix:PDEGeneralOperators}

The operators $\LL_d$ and $\LL_d^*$ depend on $d$, and so do their domains $D(\LL_d)$ and $D(\LL_d^*)$. We let 
$$
D(\LL_d)
=
\begin{cases}
\SC^d(\R^d)' & \textrm{ if } d \textrm{ is odd,}\\
\SC_{1/2}^d(\R^d)' & \textrm{ if } d \textrm{ is even},\\
\end{cases}
\quad \text{and}\quad
D(\LL_d^*)
=
\begin{cases}
\SC(\R^d)' & \textrm{ if } d \textrm{ is odd,}\\
\SC_{1/2}(\R^d)' & \textrm{ if } d \textrm{ is even}.\\
\end{cases}
$$
Recall that $\SC(\R^d)'$ (resp. $\SC^d(\R^d)'$) denotes the space of $\C$ (resp. $\C^d$)--valued tempered distributions (see Appendix \ref{sec:Distributions}); see also Section \ref{sec:IntroLaplace} of the main paper for the definition of $\SC_{1/2}(\R^d)'$ and $\SC_{1/2}^d(\R^d)'$. Since $\SC^d_{1/2}(\R^d)'\subset \SC^d(\R^d)'$, we  have
$$
\SC^d_{1/2}(\R^d)'
\subset 
D(\LL_d)
\subset
\SC^d(\R^d)'
,\spa 
\forall\ d\in\N
.
$$ 
Also notice that since $(L^\infty(\R^d))^d\subset \SC_{1/2}^d(\R^d)'$ and $\Fpg\in (L^\infty(\R^d))^d$, then $\Fpg\in D(\LL_d)$ for any $d$.
%The domain $D(\LL_d)$ of $\LL_d$ is a space of distributions that depends on $d$, namely we let $D(\LL_d)=\SC^d(\R^d)'$, the space of $\C^d$-valued tempered distributions (see Appendix \ref{sec:Distributions}) if $d$ is odd, and $D(\LL_d)=\SC_{1/2}^d(\R^d)'$ (see Section \ref{sec:IntroLaplace}) if $d$ is even.

\begin{Defin}\label{DefinDiffOperator}
Let $d\in\N$ with $d\geq 1$. Define the operator $\LL_d:D(\LL_d)\to \SC(\R^d)'$ by letting
$$
\LL_d 
:=
\gamma_d
\begin{cases}
	(-\Delta)^{\frac{d-1}{2}} \nabla\ \cdot & \textrm{ if } d \textrm{ is odd,}\\
	(-\Delta)^{\frac{1}{2}} (-\Delta)^{\frac{d-2}{2}} \nabla\ \cdot & \textrm{ if } d \textrm{ is even,}\\
\end{cases}
$$
%	with domain
%	$$ 
%	D(\LL_d)
%	=
%	\begin{cases}
	%		\SC^d(\R^d)' & \textrm{ if } d \textrm{ is odd,}\\
	%		\SC_{1/2}^d(\R^d)' & \textrm{ if } d \textrm{ is even,}\\
	%	\end{cases}
%	$$
%(for the definition of $\SC^d(\R^d)'$ and $\SC_{1/2}^d(\R^d)'$, see Appendix \ref{sec:Distributions} and Section \ref{sec:IntroLaplace}, respectively), 
where $\nabla\ \cdot$ is the divergence operator, $(-\Delta)^k$ stands for the Laplacian operator $-\Delta$ taken $k$ times successively when $k\in\N$, and $(-\Delta)^{1/2}$ denotes the fractional Laplacian introduced in Section \ref{sec:IntroLaplace} of the main paper.
%on the domain $D(\LL_d)$, that also depends on $d$, made of tempered distributions and defined by
%$$
%D(\LL_d)
%=
%\begin{cases}
%	\SC^d(\R^d)' & \textrm{ if } d \textrm{ is odd,}\\
%	\SC_{1/2}^d(\R^d)' & \textrm{ if } d \textrm{ is even,}\\
%\end{cases}
%$$
%
%The domain $D(\LL_d)$ of $\LL_d$ is a space of tempered distributions that also depends on $d$ as we let 
%\end{Defin}
%\begin{Defin}[The adoint $\LL_d^*$]\label{DefinDiffOperatorAdjoint}
%	Let $n\in\N$ with $d\geq 1$. 
%	, and 
%	$$
%	D(\LL_d^*)
%	=
%	\begin{cases}
	%		\SC(\R^d)' & \textrm{ if } d \textrm{ is odd,}\\
	%		\SC_{1/2}(\R^d)' & \textrm{ if } d \textrm{ is even}.\\
	%	\end{cases}
%	$$
\end{Defin}

%\vspace{3mm}

We then define the formal adjoint $\LL_d^*:D(\LL_d^*)\subset \SC(\R^d)'\to \SC^d(\R^d)'$ of $\LL_d$ as
$$
\LL_d^*
:=
\gamma_d
\begin{cases}
\nabla (-\Delta)^{\frac{d-1}{2}} & \textrm{ if } d \textrm{ is odd,}\\
\nabla (-\Delta)^{\frac{1}{2}} (-\Delta)^{\frac{d-2}{2}} & \textrm{ if } d \textrm{ is even,}\\
\end{cases}
$$
where $\nabla$ stands for the gradient operator. We call $\LL_d^*$ the formal adjoint of $\LL_d$ because we have
$$
\ps{\LL_d \Lambda}{\psi}_{\SC',\SC}
=
\ps{\Lambda}{\LL_d^*\psi}_{(\SC^d)',\SC^d}
,\spa \forall\ \Lambda\in D(\LL_d),\ \forall\ \psi\in\SC(\R^d)
,
$$
and 
$$
\ps{\LL_d^* T}{\Psi}_{(\SC^d)',\SC^d}
=
\ps{T}{\LL_d \Psi}_{\SC',\SC}
,\spa \forall\ T\in D(\LL_d^*),\ \forall\ \Psi\in\SC(\R^d,\C^d)
.
$$
In particular, we have
$$
\int_{\R^d} (\LL_d \Psi)(x) \varphi(x)\, dx 
=
\int_{\R^d} \ps{\Psi(x)}{(\LL_d^* \varphi)(x)}\, dx
,\spa 
\forall\ \Psi\in\SC(\R^d,\C^d),\ \forall\ \varphi\in\SC(\R^d)
.
$$
%for any $\Psi\in\SC(\R^d,\C^d)$ and $\varphi\in\SC(\R^d)$.
%\begin{Defin}[The adjoint operator $\LL_d^*$]\label{DefinDiffOperatorAdjoint}
%	Let $n\in\N$ 
%\end{Defin}
Taking Fourier transforms shows that all differential operators involved in the definition of $\LL_d$ and $\LL_d^*$ commute over $D(\LL_d)$ and $D(\LL_d^*)$, respectively. This legitimates writing $\LL_d$ and $\LL_d^*$ in the more compact forms 
$$
\LL_d
=
\gamma_d\ (-\Delta)^{\frac{d-1}{2}}\nabla\ \cdot\ ,
\quad\textrm{and}\quad
\LL_d^* 
= 
\gamma_d\ \nabla (-\Delta)^{\frac{d-1}{2}},
$$
irrespective of $d\geq 1$.

\subsection{Proof of Theorem \ref{TheorEDPSchwartz}.}

The proof of Theorem \ref{TheorEDPSchwartz} requires Lemma \ref{LemDerivOK} which, in turn, relies on Lemma \ref{GreenFormula}, and can be found in Appendix C.2 of \cite{Evans1998}.

\begin{Lem}\label{GreenFormula}
Let $\Omega\subset\R^d$ be a regular and bounded open subset, and $\partial\Omega$ denote its boundary. Let $f,g\in \CC^1(\overline{\Omega})$. For any $i\in\{1,\ldots, n\}$, we have 
$$
\int_\Omega (\partial_i f)(x) g(x)\, dx 
=
\int_{\partial \Omega} f(x)f(x)\nu_i(x)\, d\sigma(x)
-
\int_\Omega f(x) (\partial_i g)(x)\, dx 
,
$$
where $\nu_i(x)$ is the $i$th component of the outer unit normal vector to $\Om$ at $x$, and $\sigma$ is the surface area measure on $\partial \Om$.
\end{Lem}

\vspace{3mm}

\begin{Lem}\label{LemDerivOK}
Let $d\geq 1$ and $\theta:\R^d\to\R$ be such that $\theta\in L^1(\BB^d)\cap L^\infty(\R^d\sm\BB^d)\cap \CC^1(\R^d\sm \{0\})$ and $\nabla \theta\in L^\infty(\R^d\sm\BB^d)$. Further assume that there exists a positive constant $C$ such that $\|x\|^{d-1}|\theta(x)|\leq C$ for all $x\in\BB^d$, and that the limit 
$$
\kappa(u)
:=
\lim_{\eta\downarrow 0} \eta^{d-1}\theta(\eta u)
%	,\spa 
%	\text{and} 
%	\spa
%	\kappa_{\infty}(u)
%	:=
%	\lim_{R\uparrow \infty} R^{d-1}\theta(R u)
$$
exists for all $u\in \SS^{d-1}$. Denote by $\sigma$ the surface area measure on $\SS^{d-1}$ and let
$$
\kappa_{i}
:=
\int_{\SS^{d-1}} \kappa(u)u_i\, d\sigma(u)
%	,\spa 
%	\text{and}
%	\spa 
%	\kappa_{\infty,i}
%	:=
%	\int_{\SS^{d-1}} \kappa_\infty(u)u_i\, d\sigma(u)
,
$$
for all $i=1,2,\ldots, d$, where $u_i$ denotes the $i$th component of $u$. Then the (distributional) derivative $\widetilde{\partial_i} \theta$ of $\theta$ is given by 
$$
\widetilde{\partial_i} \theta
=
\PV(\partial_i \theta)
+
\kappa_{i} \delta 
,
$$
where $\partial_i \theta$ stands for the pointwise $i$th partial derivative of $\theta$ defined on $\R^d\sm\{0\}$.
%	, i.e. 
%	$$
%	\ps{\tilde{\partial_i}\theta}{\psi}_{\SC',\SC}
%	=
%	\lim_{\eta\downarrow 0} \int_{\R^d\sm \BB^d_\eta} (\partial_i \theta)(z)\ \psi(z)\, dz 
%	+
%	\kappa_{0,i} \psi(0)
%	,\spa 
%	\forall\ \psi\in\SC(\R^d)
%	.
%	$$
%	We further have
%	$$
%	\ps{\tilde{\partial_i}\theta}{\psi}_{\SC',\SC}
%	=
%	\lim_{\substack{\eta\downarrow 0\\ R\uparrow\infty}} \int_{\BB^d_R\sm \BB^d_\eta} (\partial_i \theta)(z)\ (\psi(z)-\psi(0))\, dz 
%	+
%	\kappa_{\infty,i} \psi(0)
%	,\spa 
%	\forall\ \psi\in\SC(\R^d)
%	.
%	$$ 
%	(ii) If we further assume that $\|x\|^{d-1} |\theta(x)|\leq C$ holds for all $x\in\R^d$, and that the limit 
%	$$
%	\kappa_\infty(u)
%	:=
%	\lim_{R\uparrow\infty} R^{d-1} \theta(Ru)
%	$$
%	exists for all $u\in \SS^{d-1}$, then 
%	$$
%	\ps{\tilde{\partial_i} \theta}{\psi}
%	=
%	\kappa_{\infty,i}
\end{Lem}

%\textcolor{red}{Examples: $1/\|x\|^{d-1}$, $K_d$}

In particular, it follows from Lemma \ref{LemDerivOK} that
$$
\nabla\Big(\frac{1}{\|x\|^{d-1}}\Big) 
=
-(d-1)\PV\Big(\frac{x}{\|x\|^{d+1}}\Big) 
$$
holds in $\SC(\R^d)'$, and that for all $\alpha\in\N^d$ there exists a real constant $\kappa_\alpha$ such that
$$
\partial^\alpha K 
=
\PV(\partial^\alpha K) + \kappa_\alpha \delta
,
$$
holds in $\SC(\R^d)'$, where $\PV(\partial^\alpha K)$ denotes the Principal Value of the pointwise derivative of $K$ on $\R^d\sm\{0\}$---see the Notation subsection of Section \ref{sec:Introduction}.

\begin{Proof}{Lemma \ref{LemDerivOK}}
Fix $\psi\in\SC(\R^d)$ and $i\in\{1,2\ldots, d\}$. Because, the product $\theta\ \partial_i \psi$ is integrable on $\R^d$, we have 
$$
\ps{\widetilde{\partial_i}\theta}{\psi}_{\SC',\SC}
=
- \ps{\theta}{\partial_i\psi}_{\SC',\SC}
=
- \int_{\R^d} \theta(x) (\partial_i\psi)(x)\, dx
=
- \lim_{\substack{R\uparrow\infty\\ \eta\downarrow 0}} \int_{\BB^d_R\sm \BB^d_\eta} \theta(x) (\partial_i\psi)(x)\, dx
.
$$
Fix $R,\eta>0$. Denoting by $\partial(\BB^d_R\sm\BB^d_\eta)$ the boundary of $\BB^d_R\sm\BB^d_\eta$ and $\nu(x):=x/\|x\|$ the outer unit normal vector to that boundary at $x$, Lemma \ref{GreenFormula} entails that 
%	\begin{eqnarray*}
	%		\lefteqn{
		%			
		%		}
	%		\\[2mm]
	%		
	%	\end{eqnarray*}
$$
- \int_{\BB^d_R\sm \BB^d_\eta} \theta(x) (\partial_i\psi)(x)\, dx
=
- \int_{\partial(\BB^d_R\sm\BB^d_\eta)} \theta(x)\psi(x)\nu_i(x)\, d\sigma(x)
+
\int_{\BB^d_R\sm \BB^d_\eta} (\partial_i \theta)(x)\psi(x)\, dx
.
$$
We have 
$$
- \int_{\partial(\BB^d_R\sm\BB^d_\eta)} \theta(x)\psi(x)\nu_i(x)\, d\sigma(x)
=
\int_{\partial(\BB^d_\eta)} \theta(x)\psi(x)\nu_i(x)\, d\sigma(x)
- 
\int_{\partial(\BB^d_R)} \theta(x)\psi(x)\nu_i(x)\, d\sigma(x)
.
$$
Substituting $u=x/\eta$ and $u=x/R$ in the first and second integral of the r.h.s. of the last equality, respectively, yields 
\begin{eqnarray*} 
	\lefteqn{
		\hspace{-15mm}
		- \int_{\partial(\BB^d_R\sm\BB^d_\eta)} \theta(x)\psi(x)\nu_i(x)\, d\sigma(x)
		=
		\eta^{d-1} \int_{\SS^{d-1}} \theta(\eta u)\psi(\eta u)u_i\, d\sigma(u)
	}
	\\[2mm]
	&&
	\hspace{40mm}
	- 
	R^{d-1} \int_{\SS^{d-1}} \theta(R u)\psi(R u)u_i\, d\sigma(u)
	.
\end{eqnarray*}
Because $\eta^{d-1} \theta(\eta u)$ is bounded in $u$, uniformly in $\eta$, and converges pointwise to $\kappa(u)$ as $\eta\downarrow 0$, the first integral of the r.h.s. converges to 
$$
\int_{\SS^{d-1}} \kappa(u)\psi(0) u_i\, d\sigma(u)
=
\kappa_i \psi(0)
,
$$
and because $\theta$ is bounded on $\R^d\sm\BB^d$ and $\psi\in\SC(\R^d)$, the second integral of the r.h.s. converges to $0$ as $R\uparrow\infty$. It follows that
$$
\ps{\widetilde{\partial_i}\theta}{\psi}_{\SC',\SC}
=
\kappa_i \psi(0) 
+
\lim_{\eta\downarrow 0} \int_{\R^d \sm \BB^d_\eta} (\partial_i \theta)(x)\psi(x)\, dx
=
\ps{\PV(\partial_i \theta)+\kappa_i\ \delta}{\psi}_{\SC',\SC}
,
$$
which concludes the proof.
\end{Proof}

\vspace{3mm} 

We now proceed with the proof of Theorem \ref{TheorEDPSchwartz}.

\begin{Proof}{Theorem \ref{TheorEDPSchwartz}}
Recall that 
$
\Fpg(x)
%	=
%	(K*P)(x)
=
(K*f_P)(x)
%	=
%	\int_{\R^d} \frac{x-z}{\|x-z\|}f(z)\, dz
$
for all $x\in\R^d$, where $K$ is the kernel introduced at the beginning of Section \ref{sec:PDEGeneral} of the main paper. 
%We argued at the beginning of Section \ref{sec:PDEGeneral} that $\Fpg\in$. Because $\Fpg\in L^\infty(\R^d,\R^d)$, we have $\Fpg\in\SC_{1/2}^d(\R^d)' \subset\SC^d(\R^d)'$; see Section \ref{sec:IntroLaplace}. In particular, $\Fpg$ belongs to the domain $D(\LL_d)$ of $\LL_d$, irrespective of $d$. 
Because $K\in\SC^d(\R^d)'$ and $f_P\in\SC(\R^d)$, Proposition \ref{PropConvolution} entails that $\Fpg\in\CC^\infty(\R^d)$. When $d=1$ recall that $\Fpg = 2F_P-1$, where 
$$
F_P(x)
=
\int_{-\infty}^x f_P(t)\, dt
$$
is the cdf of $P$; see Section \ref{sec:Introduction} in the main paper. Because $f_P$ is continuous, the fundamental theorem of calculus yields
$$
(\LL_d \Fpg)(x)
=
\gamma_d\
(-\Delta)^{\frac{d-1}{2}} (\nabla\cdot \Fpg)(x)
=
\frac{1}{2} \frac{d \Fpg}{dx} (x)
=
\frac{1}{2} \frac{d(2F_P-1)}{dx}(x)
=
f_P(x)
.
$$
Therefore, the claim is proved when $d=1$.

Then, assume that $d\geq 2$. Because $K\in L^\infty(\R^d,\R^d)$, we have $K\in\SC^d(\R^d)'$. Therefore, Proposition \ref{PropConvolution} entails that the equality 
%	the comments above that $\Fpg\in\SC(\R^d)'$, since $f_P\in\SC(\R^d)$. By , we have
\begin{equation}
	\label{TFR}
	\FT (\Fpg) = \FT(K)\ \FT( f_P )
\end{equation}
holds in $\SC^d(\R^d)'$, because $f_P\in\SC(\R^d)$. Lemma \ref{LemFourierTransformRiesz} and the fact that $\Gamma(1/2)=\sqrt{\pi}$ yield
$$
\FT\Big(\frac{1}{\|x\|}\Big)(\xi)
=
\frac{\Gamma(\frac{d-1}{2})}{\pi^{\frac{d-1}{2}}}
\frac{1}{\|\xi\|^{d-1}} 
$$
in $\SC(\R^d)'$; recall that $1/\|x\|$ is a tempered distribution on $\R^d$ because $d\geq 2$. From the identities stated before Proposition \ref{PropConvolution}, we deduce that 
$$
\FT (K)
=
\FT\Big(\frac{x}{\|x\|}\Big)
=
-
\frac{1}{2i\pi} \nabla \FT\Big(\frac{1}{\|x\|}\Big)
=
-\frac{1}{2i\pi} \frac{\Gamma(\frac{d-1}{2})}{\pi^{\frac{d-1}{2}}}
\nabla \Big(\frac{1}{\|\xi\|^{d-1}}\Big)
$$
in $\SC^d(\R^d)'$. Recalling that $x\Gamma(x)=\Gamma(x+1)$ for all $x>0$, Lemma \ref{LemDerivOK} yields
\begin{equation}
	\label{TFN}
	(\FT K)(\xi)
	=
	\frac{\Gamma(\frac{d+1}{2})}{i \pi^{\frac{d+1}{2}}}\ \PV\Big(\frac{\xi}{\|\xi\|^{d+1}}\Big)
\end{equation}
in $\SC^d(\R^d)'$. Equation \eqref{TFR} then rewrites
$$
(\FT \Fpg)(\xi)
=
\frac{\Gamma(\frac{d+1}{2})}{i \pi^{\frac{d+1}{2}}}\ \PV\Big(\frac{\xi}{\|\xi\|^{d+1}}\Big) 
(\FT f_P )(\xi)
.
$$
%in $\SC^d(\R^d)'$. 
%It is easy to see that
On the one hand, we have 
$$
\ps{\xi}{ \PV\Big(\frac{\xi}{\|\xi\|^{d+1}}\Big)}
:=
\sum_{i=1}^d \xi_i\ \PV\Big( \frac{\xi_i}{\|\xi\|^{d+1}}\Big)
=
\frac{1}{\|\xi\|^{d-1}}
$$
in $\SC(\R^d)'$.
%	and $\xi_i \frac{\xi_i}{\|\xi\|^{d+1}}$ is to be understood, in accordance with point $5$, as the product of the smooth function $\xi_i$, every derivative of which has at most polynomial growth, with the distribution $\frac{\xi_i}{\|\xi\|^{d+1}}$ described above. 
%Further note that 
On the other hand, we have
$$
\ps{\xi}{(\FT \Fpg)(\xi)}
=
\sum_{i=1}^d \xi_i (\FT(\Fpg)_i)(\xi)
=
\frac{1}{2i\pi} \sum_{i=1}^d \FT(\partial_i (\Fpg)_i)(\xi)
=
\frac{1}{2i\pi}\big(\FT(\nabla \cdot \Fpg)\big)(\xi)
$$
in $\SC(\R^d)'$. It follows that 
\begin{equation}\label{EqEDPSchwartzIntermediateStep}
	\frac{1}{2i\pi} \FT(\nabla\cdot \Fpg)(\xi) 
	=
	\frac{\Gamma(\frac{d+1}{2})}{i \pi^{\frac{d+1}{2}}} \frac{1}{\|\xi\|^{d-1}} (\FT f_P )(\xi)
\end{equation}
holds in $\SC(\R^d)'$. Consider two cases. (A) Assume that $d\geq 3$ is odd. Then, $\frac{d-1}{2}\in \N$ and we have
$$
(\FT f_P )(\xi)
=
\frac{1}{2}\frac{\pi^{\frac{d-1}{2}}}{\Gamma(\frac{d+1}{2})} \|\xi\|^{d-1} \FT(\nabla\cdot \Fpg)(\xi) 
= 
\gamma_d\
\FT\big((-\Delta)^{\frac{d-1}{2}} \nabla\cdot \Fpg\big)(\xi)
$$
in $\SC(\R^d)'$, where
$
\gamma_d^{-1}
=
2^d\pi^{(d-1)/2}\Gamma((d+1)/2)
.
$
%	Since $d$ is odd then, according to point $3$, we have
%	$$
%	\|\xi\|^{d-1} \FT(\nabla\cdot \Fpg)(\xi) 
%	= 
%	\frac{1}{(2i\pi)^{d-1}} \FT(\Delta^{\frac{d-1}{2}} \nabla\cdot \Fpg)(\xi)
%%	=
%%	\frac{(-1)^{\frac{d-1}{2}}}{(2\pi)^{d-1}} \FT(\Delta^{\frac{d-1}{2}} \nabla\cdot \Fpg)(\xi)
%	, 
%	$$
In particular, the equality
$$
f_P
=
\gamma_d\
(-\Delta)^{\frac{d-1}{2}} (\nabla\cdot \Fpg)
=
\LL_d (\Fpg)
$$
holds in $\SC(\R^d)'$. The fact that $\Fpg$ belongs to $\CC^\infty(\R^d)$ ensures that the r.h.s. of the last equality is a continuous function. Because $f_P$ is also continuous, and equality holds in the sense of distributions, equality also holds pointwise. (B) Assume that $d\geq 2$ is even. Because $d-2$ is even, we deduce from (\ref{EqEDPSchwartzIntermediateStep}) that
$$
\frac{(\FT f_P )(\xi)}{\|\xi\|} 
=
\frac{1}{2}
\frac{\pi^{\frac{d-1}{2}}}{\Gamma(\frac{d+1}{2})}
\|\xi\|^{d-2} \FT(\nabla\cdot \Fpg)(\xi) 
=
\frac{1}{2}
\frac{\pi^{\frac{d-1}{2}}}{\Gamma(\frac{d+1}{2})}
\frac{1}{(2\pi)^{d-2}} 
\FT\big((-\Delta)^{\frac{d-2}{2}} (\nabla\cdot \Fpg)\big)(\xi) 
%	=
%	\frac{\Gamma(\frac{d+1}{2})}{i \pi^{\frac{d+1}{2}}} 
$$
holds in $\SC(\R^d)'$. Recalling that $\Fpg\in \SC_{1/2}^d(\R^d)'$ and that $\SC_{1/2}^d(\R^d)'$ is closed with respect to differentiation, we see that $u:=(-\Delta)^{(d-2)/2} (\nabla\cdot \Fpg)$ belongs to $\SC_{1/2}(\R^d)'$. Because $f_P\in\SC(\R^d)$, the previous equality yields $\|\xi\|\FT u(\xi)\in L^1(\R^d)$. Since $d\geq 2$, we also have $(\FT f_P )(\xi)/\|\xi\|\in L^1_\loc(\R^d)$, hence $\F u \in L^1_\loc(\R^d)$.
%		
%		Observe that $\wwidehat{f_P}(\xi)/\|\xi\|\in L^1_\loc(\R^d)$ (recall that $d\geq 2$), and that $\|\xi\|\FT u(\xi)\in L^1(\R^d)$ from the previous equality since $f_P\in\SC(\R^d)$
%		
%		
%		 It is clear that $\FT u\in L^1_\loc(\R^d)$ since $\wwidehat{f_P}(\xi)/\|\xi\|\in L^1_\loc(\R^d)$ (recall that $d\geq 2$) and $\|\xi\|\FT u(\xi)\in L^1(\R^d)$ from the previous equality since $f_P\in\SC(\R^d)$. 
Proposition \ref{PropLaplaceExplicit} (i) thus entails that $(-\Delta)^{1/2} (-\Delta)^{(d-2)/2} (\nabla\cdot \Fpg)\in\CC_0(\R^d)$ and that
$$
\FT (f_P)
=
\gamma_d\ \FT ( (-\Delta)^{\frac{1}{2}} (-\Delta)^{\frac{d-2}{2}} (\nabla\cdot \Fpg ))
$$
holds in $\SC(\R^d)'$, where $\gamma_d$ is the same constant as in (A). We deduce that 
$$
f_P
=
\gamma_d\ 
(-\Delta)^{\frac{1}{2}} (-\Delta)^{\frac{d-2}{2}} (\nabla\cdot \Fpg)
=
\LL_d (\Fpg)
$$
holds in $\SC(\R^d)'$. Since both sides of this last equality are continuous, equality also holds pointwise on $\R^d$, which concludes the proof.
\end{Proof}

\subsection{Proof of Theorem \ref{TheorEDPDistributions}}

The proof of Theorem \ref{TheorEDPDistributions} requires the following Lemma, which is Corollary 2.2.10 in \cite{BogachevMeasures}.

\begin{Lem}\label{LemConvDistributionNotContinuous}
Let $Q$ and $(Q_k)_{k\geq 1}$ be Borel probability measures on $\R^d$ such that $(Q_k)$ converges to $Q$ in distribution as $k\to\infty$. For any bounded and measurable map $g:\R^d\to\C$ such that $Q(D_g)=0$, where we let 
$$
D_g
:=
\{x\in\R^d : g \textrm{ is not continuous at } x\}
,
$$
we have $\int_{\R^d} g\, d Q_k \to \int_{\R^d} g\, dQ$ as $k\to\infty$.
\end{Lem}

\vspace{3mm}

We now proceed with the proof of Theorem \ref{TheorEDPDistributions}.

\begin{Proof}{Theorem \ref{TheorEDPDistributions}}
Assume first that there exists a sequence of Borel probability measures $(Q_k)$ on $\R^d$ such that $(Q_k)$ converges in law to $P$ as $k\to\infty$, and such that $Q_k$ admits a density $f_k\in\SC(\R^d)$ with respect to the Lebesgue measure for any $k$. Let $F_{Q_k}^\g$ denote the geometric cdf associated with the probability measure $Q_k$, for any $k$. Because $Q_k$ admits the density $f_k\in\SC(\R^d)$, Theorem \ref{TheorEDPSchwartz} entails that
\begin{equation}\label{eq:EDPSchwartzQk}
	f_k(x)
	=
	%	\gamma_d\  (-\Delta)^{\frac{d-1}{2}} (\nabla\cdot R_{Q_k})(x)
	(\LL_d F_{Q_k}^\g)(x)
	,\spa \forall\ k,\ \forall\ x\in\R^d
	.
\end{equation}
Letting $\psi\in\SC(\R^d)$, (\ref{eq:EDPSchwartzQk}) reads
\begin{equation}\label{Eq:TheorEDPDistributions1}
	\int_{\R^d} \psi(x) f_k(x)\, dx 
	=
	%	\gamma_d \int_{\R^d} \Big( R_{Q_k}(x), \nabla \big( (-\Delta)^{\frac{d-1}{2}} \psi \big) (x)  \Big)\, dx
	\int_{\R^d} \ps{F_{Q_k}^\g(x)}{ (\LL_d^* \psi)(x)}\, dx
	,
\end{equation}
for all $k$. We are going to show that the l.h.s. of (\ref{Eq:TheorEDPDistributions1}) converges to $\int_{\R^d} \psi(x)\, dP(x)$, and the r.h.s. of (\ref{Eq:TheorEDPDistributions1}) converges to 
$$
\int_{\R^d} \ps{\Fpg(x)}{(\LL_d^* \psi)(x)}\, dx
$$ 
as $k\to\infty$. Because $(Q_k)$ converges in law to $P$ as $k\to\infty$, and $\psi$ is continuous and bounded on $\R^d$, we have 
\begin{equation}\label{Eq:TheorEDPDistributions2}
	\int_{\R^d} \psi(x)f_k(x)\, dx 
	=
	\int_{\R^d} \psi(x)\, dQ_k(x)
	\to 
	\int_{\R^d} \psi(x)\, dP(x)
	%	= 
	%	\ps{P}{\psi}
\end{equation}
as $k\to\infty$. Let us now show that the r.h.s. of \eqref{Eq:TheorEDPDistributions1} converges. We first show that $F_{Q_k}^\g$ converges almost everywhere to $\Fpg$ as $k\to\infty$. For all $x\in\R^d$, define 
$$
g_x(z)
:=
\frac{x-z}{\|x-z\|}\Ind{z\neq x}
,\spa \forall\ z\in\R^d
.
$$
With the notations of Lemma \ref{LemConvDistributionNotContinuous}, we have $D_{g_x}=\{x\}$. Define 
$
A
:=
\{x\in\R^d : P[\{x\}]>0\}
.
$
Then, $A$ is at most countable and we have $P[D_{g_x}]=0$ for all $x\in\R^d\sm A$. Because $g_x$ is bounded and measurable for any $x\in\R^d$, Lemma \ref{LemConvDistributionNotContinuous} entails that 
$$
F_{Q_k}^\g(x)
=
\int_{\R^d} g_x(z)\, dQ_k(z)
\to 
\int_{\R^d} g_x(z)\, dP(z)
=
\Fpg(x)
,\spa \forall\ x\in\R^d\setminus A
,
$$
as $k\to\infty$. Because $A$ is at most countable, $F_{Q_k}^\g$ converges to $\Fpg$ almost everywhere as $k\to\infty$. To apply Lebesgue's dominated convergence theorem to the r.h.s. of \eqref{Eq:TheorEDPDistributions1}, observe that $\LL_d^*(\psi)\in L^1(\R^d)$.
%	$\nabla \big( (-\Delta)^{\frac{d-1}{2}} \psi \big)\in L^1(\R^d)$. 
Indeed, if $d$ is even we have $(-\Delta)^{(d-2)/2}\psi\in\SC(\R^d)$ since $\psi\in\SC(\R^d)$. It follows that
$$(-\Delta)^{\frac{1}{2}}((-\Delta)^{\frac{d-2}{2}}\psi)\in \SC_{1/2}(\R^d),$$ whence 
$$
\LL_d^*(\psi)
=
\gamma_d\
\nabla \Big((-\Delta)^{\frac{1}{2}}(-\Delta)^{\frac{d-2}{2}}\psi\Big)
\in \SC_{1/2}(\R^d,\C^d)\subset L^1(\R^d)
.
$$  
If $d$ is odd, $\nabla \big( (-\Delta)^{(d-1)/2} \psi \big)$ obviously belongs to $\SC(\R^d,\C^d)$, which is a subset of $L^1(\R^d)$. Because the sequence of functions $(F_{Q_k}^\g)_k$ is uniformly norm-bounded by $1$, and $F_{Q_k}^\g$ converges to $\Fpg$ almost everywhere as $k\to\infty$, Lebesgue's dominated convergence theorem entails that 
\begin{equation}\label{Eq:TheorEDPDistributions3}
	\int_{\R^d} \ps{F_{Q_k}^\g(x)}{(\LL_d^*\psi)(x)}\, dx 
	\to 
	\int_{\R^d} \ps{\Fpg(x)}{(\LL_d^*\psi)(x)}\, dx 
\end{equation}
as $k\to\infty$. Putting \eqref{Eq:TheorEDPDistributions1}, \eqref{Eq:TheorEDPDistributions2}, and \eqref{Eq:TheorEDPDistributions3} together yields
$$
\int_{\R^d} \psi(x)\, dP(x)
=
\int_{\R^d} \ps{\Fpg(x)}{(\LL_d^* \psi)(x)}\, dx 
,\spa \forall\ \psi\in\SC(\R^d)
.
$$
It follows that
$
P
=
\LL_d(\Fpg)
$
in $\SC(\R^d)'$.

It remains to show that there exists a sequence of Borel probability measures $(Q_k)$ on $\R^d$ converging in law to $P$ as $k\to\infty$, such that $Q_k$ admits a density $f_k\in\SC(\R^d)$ with respect to the Lebesgue measure for any $k$. Let $Z$ be a random $d$-vector with law $P$. Let $\rho\in\CC^\infty_c(\R^d)$ be such that $0\leq \rho\leq 1$ and $\int_{\R^d} \rho(x)\, dx =1$. In particular, $\rho$ is a probability density on $\R^d$. Let $Y$ be a random $d$-vector with density $\rho$ and independent from $Z$. For any $k$, define 
$
Z_k 
:= 
Z + Y/k
.
$
Because $(Z_k)$ converges to $Z$ in probability as $k\to\infty$, $(Z_k)$ converges in law to $Z$. Furthermore, observe that $Z_k$ admits the density $s_k:=\rho_k*P$ with respect to the Lebesgue measure, where $\rho_k(x):= k^d \rho(kx)$ for all $k$ and $x\in\R^d$. In particular, $s_k\in\CC^\infty(\R^d)$ since $\rho_k\in\CC^\infty_c(\R^d)$. Because $(Z_k)$ converges in law to $Z$, we have 
$$
\int_{\R^d} g(x) s_k(x)\, dx 
\to 
\int_{\R^d} g(x)\, dP(x)
$$
for any continuous and bounded map $g:\R^d\to\C$ as $k\to\infty$. For all $k$ there exists $r_k>0$ such that 
$
\int_{\R^d\sm \BB^d_{r_k}} s_k(x)\, dx 
< 
1/k
.
$
For all $k$, let $\chi_k\in\CC^\infty_c(\R^d)$ be such that $0\leq \chi_k\leq 1$ on $\R^d$, $\chi_k = 1$ over $\BB^d_{r_k}$, and $\chi_k=0$ over $\R^d\sm \BB^d_{1+r_k}$. Then, define $f_k(x):= \chi_k(x) s_k(x)$ for any $k$ and $x\in\R^d$. Because $(s_k)\subset\CC^\infty(\R^d)$, we have $(f_k)\subset\CC^\infty_c(\R^d)$. In particular, $(f_k)\subset\SC(\R^d)$. Let $g:\R^d\to\C$ be a continuous and bounded map. We have
\begin{eqnarray*}
	\lefteqn{
		%			\hspace{-10mm}
		\Big| 
		\int_{\R^d} g(x) f_k(x)\, dx 
		-
		\int_{\R^d} g(x)s_k(x)\, dx
		\Big|	
		\leq 
		\int_{\R^d} |g(x)| (\chi_k(x)-1) s_k(x)\, dx 
	}	
	\\[2mm]
	&&
	\leq 
	\int_{\R^d\sm \BB^d_{r_k}} |g(x)| s_k(x)\, dx 
	\leq 
	\|g\|_{L^\infty(\R^d)}\int_{\R^d\sm \BB^d_{r_k}} s_k(x)\, dx 
	\leq 
	\frac{1}{k} \|g\|_{L^\infty(\R^d)}
	.
\end{eqnarray*}
Because $\int_{\R^d} g(x)s_k(x)\, dx\to \int_{\R^d} g(x)\, dP(x)$ as $k\to\infty$, we deduce that 
$$
\int_{\R^d} g(x)f_k(x)\, dx
\to 
\int_{\R^d} g(x)\, dP(x)
$$
for any continuous and bounded map $g:\R^d\to\C$ as $k\to\infty$. Letting $Q_k$ be the probability measure with density $f_k\in\SC(\R^d)$ for any $k$ yields the conclusion.
\end{Proof}

\subsection{Proof of Proposition \ref{PropRankIntermediateRegularity}}

\begin{Proof}{Proposition \ref{PropRankIntermediateRegularity}}
As a preliminary remark, let us observe that since $P$ is non-atomic over $\Om$, a straightforward application of Lebesgue's dominated convergence theorem yields that $\Fpg$ is continuous over $\Om$. Since $\|\Fpg(x)\|\leq 1$ for all $x\in\R^d$, we deduce that $\Fpg$ is continuous over $\Om$ and bounded on $\R^d$. 

We now prove the result by induction. Since $\Fpg\in\CC^0(\Om)$, let $0\leq k\leq \ell-1$ be such that $\Fpg\in\CC^{k}(\Om)$ with 
\begin{equation}\label{eq:PropRankIntermediate1}
	(\partial^\alpha \Fpg)(x)
	=
	\E[(\partial^\alpha K)(x-Z)]
\end{equation}
for all $x\in\Om$ and $\alpha\in\N^d$ with $|\alpha|\leq k$. Let us show that $\Fpg\in\CC^{k+1}(\Om)$ and that (\ref{eq:PropRankIntermediate1}) holds for any $\alpha\in\N^d$ with $|\alpha|\leq k+1$.
%	Let $\alpha\in\N^d$ with $0\leq |\alpha|\leq \ell-1$ and assume that $\Fpg\in\CC^{|\alpha|}(\Om)$ with $\partial^\beta \Fpg(x)=\E[(\partial^\beta K)(x-X)]$ for any $x\in\Om$ and $\beta\in\N^d$ with $|\beta|\leq |\alpha|$. Let us start by showing that $\Fpg\in\CC^{|\alpha|+1}(\Om)$. 
For this purpose, let $\alpha\in\N^d$ with $|\alpha|=k$, $x\in\Om$, and $r>0$ be such that $\overline{\BB^d_{2r}(x)}\subset\Om$. Let $j\in\{1,\ldots, n\}$ and $e_j$ be the $j$th vector of the canonical basis of $\R^d$. We are going to show that 
$$
\frac{(\partial^\alpha \Fpg)(x+h e_j)-(\partial^\alpha \Fpg)(x)}{h}
\to 
\E[(\partial_j\partial^\alpha K)(x-Z)]
$$
as $h\to 0$ and that the limit is continuous over $\Om$. Notice that for all $h$ such that $|h|<r$, we have $x+he_j\in \BB^d_{r}(x)\subset\Om$. For all $h$, define 
$$
S_h:=\{x+she_j : s\in [0,1]\}
$$ 
be the segment with endpoints $x$ and $x+he_j$. Since $S_h\subset\Om$ and $P$ has a density on $\Om$, we have 
$$
\frac{\partial^\alpha \Fpg(x+h e_j)-\partial^\alpha \Fpg(x)}{h}
=
\E\bigg[
\frac{(\partial^\alpha K)(x+he_j-Z)-(\partial^\alpha K)(x-Z)}{h}\Ind{Z\in \R^d\sm S_h}
\bigg]
$$
for any $h$. To take the limit as $h\to 0$ under the above expectation, we will show that the integrand is a uniformly $P$-integrable family indexed by $h$ and converges $P$-almost surely as $h\to 0$. Since $K\in\CC^\infty(\R^d\sm\{0\})$, observe that
$$
\frac{(\partial^\alpha K)(x+he_j-z)-(\partial^\alpha K)(x-z)}{h}
=
\int_0^1 (\partial_j \partial^\alpha K)(x+she_j-z)\, ds 
$$
for any $z\in\R^d\sm S_h$. The latter obviously converges to $(\partial_j \partial^\alpha K)(x-z)$ as $h\to 0$, for any $z\in\R^d\sm S_h$. Let us now show that the family of random vectors
$$
\bigg( \int_0^1 (\partial_j \partial^\alpha K)(x+she_j-Z)\Ind{Z\in\R^d\sm S_h}\, ds \bigg)_{|h|<r}
$$
is uniformly $P$-integrable. It is enough to show that there exists $\delta>0$ such that
$$
\sup_{|h| < r} 
\E\bigg[
\Big| \int_0^1 (\partial_j \partial^\alpha K)(x+she_j-Z)\, ds \Big|^{1+\delta}
\Ind{Z\in\R^d\sm S_h}
\bigg]
<
\infty 
.
$$
Fix $\delta>0$ and let us specify its value later. Observe that for all $\beta\in\N^d$ there exists a positive constant~$C_\beta$ such that $|\partial^\beta K(x)|\leq C_{\beta} \|x\|^{-|\beta|}$ for all $x\in\R^d\sm\{0\}$. Therefore, there exists $C>0$ such that
\begin{align*}
	\Big| 
	\int_0^1 (\partial_j \partial^\alpha K)(x+she_j-z)\, ds 
	\Big|^{1+\delta}
	&\leq
	\int_0^1 |(\partial_j \partial^\alpha K)(x+she_j-z)|^{1+\delta}\, ds 
	\\[2mm]
	&\leq 
	C\int_0^1 \frac{1}{\|x+she_j-z\|^{(1+k)(1+\delta)}}\, ds 
	,
\end{align*}
%\begin{eqnarray*}
%	\lefteqn{
	%		\hspace{-10mm}
	%		\Big| 
	%		\int_0^1 (\partial_j \partial^\alpha K)(x+she_j-z)\, ds 
	%		\Big|^{1+\delta}
	%	}
%	\\
%	&&
%	\hspace{-5mm}
%	\leq
%	\int_0^1 |(\partial_j \partial^\alpha K)(x+she_j-z)|^{1+\delta}\, ds 
%	\\
%	&&
%	\hspace{-5mm}
%	\leq 
%	C\int_0^1 \frac{1}{\|x+she_j-z\|^{(1+k)(1+\delta)}}\, ds 
%\end{eqnarray*}
for any $z\in\R^d\sm S_h$ and $h$, by Jensen's inequality. It follows from Fubini's theorem that 
\begin{eqnarray*}
	\lefteqn{
		\sup_{|h| < r} 
		\E\bigg[
		\Big| \int_0^1 (\partial_j \partial^\alpha K)(x+she_j-Z)\, ds \Big|^{1+\delta}
		\Ind{Z\in\R^d\sm S_h}
		\bigg]
	}
	\\[2mm]
	&&
	\leq 
	C \sup_{|h| < r} 
	\int_0^1
	\E\bigg[
	\frac{1}{\|x+she_j-Z\|^{(1+k)(1+\delta)}}
	\Ind{Z\in\R^d\sm S_h}
	\bigg]\, ds 
	\\[2mm]
	&&
	\leq 
	C \sup_{|h| < r} 
	\sup_{s\in [0,1]}
	\E\bigg[
	\frac{1}{\|x+she_j-Z\|^{(1+k)(1+\delta)}}
	\Ind{Z\in\R^d\sm S_h}
	\bigg]
	%		\\[2mm]
	%		&&
	%		\leq 
	%		C \sup_{|h| < \kappa} 
	%		\E\bigg[
	%		\frac{1}{|x+he_j-Z|^{1+k+\delta}}
	%		\Ind{Z\in\R^d\sm S_h}
	%		\bigg] 
	.
\end{eqnarray*}
Fix $s\in [0,1]$ and $h$ with $|h|<r$. We have 
\begin{eqnarray*}
	\lefteqn{
		\E\bigg[
		\frac{1}{\|x+she_j-Z\|^{(1+k)(1+\delta)}}
		\Ind{Z\in\R^d\sm S_h}
		\bigg] 
	}
	\\[2mm]
	&&
	\hspace{5mm}
	\leq 
	\frac{1}{r^{(1+k)(1+\delta)}}
	+
	\E\bigg[
	\frac{1}{\|x+she_j-Z\|^{(1+k)(1+\delta)}}
	\Ind{Z\in \BB^d_r(x+she_j)\sm S_h}
	\bigg] 
	.
\end{eqnarray*}
Since $P$ admits the density $f_\Om$ on $\Om$, Hölder's inequality yields
\begin{eqnarray*}
	\lefteqn{
		\hspace{-5mm}
		\E\bigg[
		\frac{1}{\|x+she_j-Z\|^{(1+k)(1+\delta)}}
		\Ind{Z\in \BB^d_r(x+she_j)\sm S_h}
		\bigg] 
	}
	\\[2mm]
	&&
	=
	\int_{\BB^d_r(x+she_j)\sm S_h} \frac{1}{\|x+she_j-z\|^{(1+k)(1+\delta)}} f_\Om(z)\, dz
	\\[2mm]
	&&
	\leq 
	\bigg(\int_{\BB^d_r} \frac{1}{\|z\|^{q(1+k)(1+\delta)}}\, dz\bigg)^{1/q}
	\bigg(\int_{\BB^d_r(x+she_j)} |f_\Om(z)|^p\, dz\bigg)^{1/p}
	,
\end{eqnarray*}
where $p$ is such that $f_\Om\in L^p_\loc(\Om)$ and $q=p/(p-1)$ is the conjugate exponent of $p$. The fact that $k\leq \ell-1$ and $p>d/(d-\ell)$ implies that 
$$
p>\frac{d}{d-(1+k)}
$$ 
and $q<d/(1+k)$. Let us therefore choose $\delta>0$ small enough such that $$
q<\frac{d}{(1+k)(1+\delta)}
.
$$
In particular, we have 
$$
\int_{\BB^d_r} \frac{1}{\|z\|^{q(1+k)(1+\delta)}}\, dz 
< 
\infty 
.
$$
Since $\BB^d_r(x+she_j)\subset \BB^d_{2r}(x)$ with $\overline{\BB^d_{2r}(x)}\subset\Om$, and $f_\Om\in L^p_\loc(\Om)$, we also have 
$$
\sup_{s\in [0,1]}
\sup_{|h|<r}
\int_{\BB^d_r(x+she_j)} |f_\Om(z)|^p\, dz
\leq 
\int_{\BB^d_{2r}(x)} |f_\Om(z)|^p\, dz
< 
\infty 
.
$$
We deduce that 
$$
\sup_{|h| < r} 
\E\bigg[
\Big| \int_0^1 (\partial_j \partial^\alpha K)(x+she_j-Z)\, ds \Big|^{1+\delta}
\Ind{Z\in\R^d\sm S_h}
\bigg]
<
\infty
.
$$
Consequently, the family of random vectors 
$$
\bigg( \int_0^1 (\partial_j \partial^\alpha K)(x+she_j-Z)\Ind{Z\in\R^d\sm S_h}\, ds \bigg)_{|h|<r}
$$
is uniformly $P$-integrable. It follows from Lebesgue-Vitali's theorem that 
$$
\frac{(\partial^\alpha \Fpg)(x+h e_j)-(\partial^\alpha \Fpg)(x)}{h}
\to 
\E[(\partial_j\partial^\alpha K)(x-Z)]
$$
for any $x\in\Om$ as $h\to 0$. Let us show that $x\mapsto \E[(\partial_j\partial^\alpha K)(x-Z)]$ is continuous over $\Om$. Fix $x\in\Om$ and a sequence $(x_m)\subset\Om$ converging to $x$ as $m\to\infty$. The familiy of random vectors $((\partial_j\partial^\alpha K)(x_m-Z))_{m\in\N}$ converges $P$-almost surely to $(\partial_j\partial^\alpha K)(x-Z)$ as $m\to\infty$ since $P$ is non-atomic, and is uniformly $P$-integrable since 
$$
\sup_{m\in\N} 
\E[|(\partial_j\partial^\alpha K)(x_m-Z)|^{1+\eta}]
\lesssim 
\sup_{m\in\N} \E\bigg[\frac{1}{\|x_m-Z\|^{(1+k)(1+\eta)}}\bigg]
<
\infty 
$$
for $\eta$ small enough, by the previous computations. It follows that $\partial^\alpha \Fpg\in \CC^1(\Om)$ and that 
$$
(\partial_j \partial^\alpha \Fpg)(x)
=
\E[(\partial_j\partial^\alpha K)(x-Z)]
$$
for all $x\in\R^d$.
%	, any $j\in\{1,\ldots, n\}$ and $\alpha\in\N^d$ with $0\leq |\alpha|=k$. 
Since $\alpha\in\N^d$ with $|\alpha|=k$ was arbitrary, we deduce that $\Fpg\in \CC^{k+1}(\Om)$  and that \eqref{eq:PropRankIntermediate1} holds for any $\alpha\in\N^d$ with $|\alpha|\leq k+1$. We conclude by induction.\\

Assume that $\Om=\R^d$, and let $f_P$ stand for the density $f_\Om\in L^p(\R^d)$. We showed that $\Fpg\in\CC^{\ell}(\R^d)$. Then, fix $\alpha\in\N^d$ with $1\leq |\alpha|\leq \ell$ and let us show that $\partial^\alpha \Fpg$ converges to $0$ at infinity. For the sake of convenience, write $k:=|\alpha|$. We have already noticed that
$$
|\partial^\alpha \Fpg(x)|
\leq 
C\ 
\E\Big[\frac{1}{\|x-Z\|^k}\Ind{Z\neq x}\Big]
=:
C\ g(x)
$$
for all $x\in\R^d$ and some positive constant $C$. Therefore, it is enough to show that $g$ converges to $0$ at infinity. Let $(x_m)\subset\R^d$ be such that $|x_m|\to\infty$ as $m\to\infty$. A standard application of Lebesgue's dominated convergence theorem entails that
$$
\E\Big[\frac{1}{\|x_m-Z\|^k}\Ind{Z\in\R^d\sm \BB^d(x_m)}\Big]
\to
0
$$
as $m\to\infty$.
Next observe that 
\begin{align*}
	\Big|\E\Big[\frac{1}{\|x_m-Z\|^k}\Ind{Z\in\BB^d(x_m)}\Big]\Big|
	&=
	\int_{\BB^d(x_m)} \frac{1}{\|z-x_m\|^k} f_P(z)\, dz 
	\\[2mm]
	&\leq 
	C\bigg(\int_{\BB^d(x_m)} |f_P(z)|^p\, dz\bigg)^{1/p} 
	\bigg(\int_{\BB^d} \frac{1}{\|z\|^{qk}}\, dz \bigg)^{1/q},
\end{align*}
%\begin{eqnarray*}
%	\lefteqn{
	%		\hspace{-3mm}
	%		\Big|\E\Big[\frac{1}{\|x_m-Z\|^k}\Ind{\BB^d(x)}\Big]\Big|
	%	}
%	\\[2mm]
%	&&
%	\hspace{3mm}
%	=
%	\int_{\BB^d(x_m)} \frac{1}{\|z-x_m\|^k} f_P(z)\, dz 
%	\\[2mm]
%	&&
%	\hspace{3mm}
%	\leq 
%	C\bigg(\int_{\BB^d(x_m)} |f_P(z)|^p\, dz\bigg)^{1/p} 
%	\bigg(\int_{\BB^d} \frac{1}{\|z\|^{qk}}\, dz \bigg)^{1/q},
%\end{eqnarray*}
where $q=p/(p-1)$ is the conjugate exponent of $p$. Since $p>d/(d-\ell)$ and $k\leq \ell$, we have~$p>d/(d-k)$ whence $q<d/k$. In particular, $qk<d$. It follows that 
$$
\int_{\BB^d} \frac{1}{\|z\|^{qk}}\, dz
<
\infty 
.
$$
It remains to show that $\int_{\BB^d(x_m)} |f_P(z)|^p\, dz \to 0$ as $m\to\infty$. Since $f_P\in L^p(\R^d)$, let $\nu$ be the non-negative finite measure defined by 
$$
\nu(B)
:=
\int_B |f_P(z)|^p\, dz 
$$
for any Borel subset $B\subset\R^d$. We then have 
$$
\int_{\BB^d(x_m)} |f_P(z)|^p\, dz 
=
\nu(\BB^d(x_m))
$$
for any $m$. Furthermore, we have $\nu(\R^d\sm \BB^d_{\|x_m\|-1})
\to 
0$ as $m\to\infty$ since $\nu$ is finite and $\|x_m\|\to \infty$ as $m\to\infty$. It follows that
$$
\E\Big[\frac{1}{\|x_m-Z\|^k}\Ind{Z\in\BB^d(x_m)}\Big]
\to 
0
$$
as $m\to\infty$. We deduce that $\partial^\alpha \Fpg$ converges to $0$ at infinity for any $\alpha\in\N^d$ with $1\leq |\alpha|\leq \ell$, which concludes the proof.
%	 Let us write 
%	 $$
%	 h(x)
%	 =
%	 (u_1*f_P)(x)+(u_2*f_P)(x)
%	 $$
%	 for any $x\in\R^d$, where $u_1(x)=\frac{1}{\|x\|^{k+1}}\Ind{0<\|x\|<1}$ and $u_2(x)=\frac{1}{\|x\|^{k+1}}\Ind{\|x\|>1}$. Since $u_2\in L^\infty(\R^d)$ and $f_P\in L^1(\R^d)$, Hausdorff-Young's inequality yields $u_2*f_P\in L^\infty(\R^d)$. We have already observed that $u_1\in L^q(\R^d)$, where $q=\frac{p}{p-1}$ is the conjugate exponent of $p$. Since $f_P\in L^p(\R^d)$ and $u_1\in L^q(\R^d)$ with $\frac1p + \frac1q = 1$, Hausdorff-Young's inequality entails that $u_1*f_P\in L^\infty(\R^d)$. It follows that $h\in L^\infty(\R^d)$, whence $\partial^\alpha \Fpg\in\CC^1_b(\R^d)$ provided $f_P\in L^p(\R^d)$.	 
%	$$
%	\partial^\beta \Fpg(x) 
%	= 
%	\E[(\partial^\beta K)(x-Z)]
%	$$ 
%	for any $x\in\Om$ and $\beta\in\N^d$ with $|\beta|\leq 1+\alpha$. 
%	Since $|\alpha|\leq \ell-1$ was arbitrary, we conclude by induction that $\Fpg\in\CC^\ell(\Om)$ with $\partial^\alpha \Fpg(x)=\E[(\partial^\alpha K)(x-Z)]$ for any $x\in\Om$ and $\alpha\in\N^d$ with $|\alpha|\leq \ell$.
\end{Proof}

\subsection{Proof of Proposition \ref{PropContinuousCounterExample}}

The proof of Proposition \ref{PropContinuousCounterExample} requires the following Lemma. 

\begin{Lem}\label{LemContinuityTransfer}
Let $X,Y,Z$ be Banach spaces, and $A:X\to Y$ be a (not necessarily continuous) linear map. Let $j:Y\to Z$ be a (not necessarily linear) continuous and injective map. If $j\circ A:X\to Z$ is continuous, then $A:X\to Y$ is a bounded operator.
\end{Lem}

\begin{Proof}{Lemma \ref{LemContinuityTransfer}}
Let $z_n=(x_n, Ax_n)\subset X\times Y$ be a sequence in the graph of $A$ such that $(z_n)$ converges in $X\times Y$ for the product topology on $X\times Y$, \ie $(x_n)\to x\in X$ and $(Ax_n) \to y\in Y$. On the one hand, the continuity of $j$ entails that $j(Ax_n)\to j(y)$ in $Z$. On the other hand, we have $j(Ax_n)=(j\circ A)(x_n)\to (j\circ A)(x)=j(Ax)$ by continuity of $j\circ A$. By uniqueness of limits in $Z$, we deduce that $j(Ax)=j(y)$. Since $j$ is injective, we then have $y=Ax$. It follows that the graph of $A$ is closed which, by the Closed Graph Theorem, entails that $A$ is continuous. Together with the fact that $A$ is a linear map, this entails that $A$ is a bounded operator.
\end{Proof}

\vspace{3mm} 

We now proceed with the proof of Proposition \ref{PropContinuousCounterExample}. 

\begin{Proof}{Proposition \ref{PropContinuousCounterExample}}
Let us assume, \emph{ad absurdum}, that the distributional derivative $\partial^\alpha \Fpg$ belongs to $L^\infty(\Om)$ for any Borel probability measure $P$ with a continuous and compactly supported density on $\R^d$, i.e. $\partial^\alpha(K*f)\in L^\infty(\Om)$ for any density $f\in\CC^0_c(\R^d)$. By linearity, we thus have $\partial^\alpha (K*f)\in L^\infty(\Om)$ for any $f\in \CC^0_c(\R^d)$, not necessarily a density. For any $f\in\CC^0_c(\R^d)$, Lemma \ref{PropConvolGeneral} entails that the equality $\partial^\alpha (K*f)=(\partial^\alpha K)*f$ holds in $\D(\R^d)'$. Lemma \ref{LemDerivOK}, and the comments thereafter, provides 
$$
\partial^\alpha K 
=
\PV(\partial^\alpha K) + \kappa_\alpha \delta
,
$$
where $\PV(\partial^\alpha K)$ denotes the Principal Value of the pointwise derivative of $K$ on $\R^d\sm\{0\}$---see the Notation subsection of Section \ref{sec:Introduction}---and $\kappa_\alpha$ is a real constant. We then have
\begin{equation}\label{eq:DerConv}
	\partial^\alpha (K*f)
	=
	\PV(\partial^\alpha K)*f + \kappa_\alpha f
	%	\\[2mm]
	=: T(f)+\kappa_\alpha f
	,\spa 
	\forall\ f\in \CC^0_c(\R^d)
	.
\end{equation}
%	 Now let us assume, \emph{ad absurdum}, that the distribution $\partial^\alpha \Fpg$ is a bounded function on $\Om$ for any Borel probability measure $P$ with a continuous and compactly supported density on $\R^d$.
%	 ; in the sequel, we will write $P\in\CC^0_c(\R^d)$ in short. 
Since $\partial^\alpha(K*f)\in L^\infty(\Om)$ for any $f\in \CC^0_c(\R^d)$, then $T(f)\in L^\infty(\Om)$ since $f\in L^\infty(\Om)$ as well. 
%	We then have a linear map $T:\CC^0_c(\R^d)\to L^\infty(\Om)$.
%	In particular, for any density $f\in\CC^0_c(\R^d)$ the distributional equality (\ref{eq:DerConv}) is in fact an equality between bounded functions on $\Om$, so that
%	%	 the convolution $(\partial^\alpha K_d)*f_P$ 
%	\begin{equation}\label{eq:SingInt}
	%		\big((\partial^\alpha K)*f\big)(x)
	%		=
	%		\lim_{\eta\downarrow 0} 
	%		\int_{\R^d\sm \BB^d_\eta(x)}
	%		%	\E\Big[(\partial^\alpha K)(x-Z)\ \I\big[\|Z-x\|>\eta\big]\Big]
	%		(\partial^\alpha K)(x-z)\  f(z)\, dz
	%		+
	%		\kappa_\alpha f(x)
	%		=:
	%		(Tf)(x)
	%		+
	%		\kappa_a f(x)
	%		,
	%	\end{equation}
%	holds for almost every $x\in\Om$. By linearity, (\ref{eq:SingInt}) holds for any $f\in \CC^0_c(\R^d)$, not necessarily a density. It follows that $Tf\in L^\infty(\Om)$ for any $f\in\CC^0_c(\R^d)$.
Fix an open and bounded subset $U\subset \Om$ such that $\overline{U}\subset \Om$. Let $V:=\overline{U}$ and define $C$ as the closure of the Minkowski sum $V + \BB^d_\ve$ with $\ve>0$ small enough such that $C\subset \Om$, \ie $x\in C$ if and only if $x$ is at distance at most $\ve$ of $V$. In particular, $C$ is compact and we have $V\subset C\subset \Om$.  For $f\in \CC^0_c(V)$, decompose $T(f)$ as
\begin{align}\label{eq:SingInt2}
	T(f)
	=
	\big( (\partial^\alpha K)\ \I[\R^d\sm\BB^d_\ve]\big)*f
	+
	\PV\Big( (\partial^\alpha K)\ \I[\BB^d_\ve]\Big)*f 
	%		\\[2mm] 
	%		\int_{\R^d\sm \BB^d_\ve(x)} (\partial^\alpha K)(x-z)\  f(z)\, dz
	%		+
	%		\lim_{\eta\downarrow 0} \int_{\BB^d_\ve(x)\sm \BB^d_\eta(x)} (\partial^\alpha K)(x-z)\  f(z)\, dz
	=:
	H(f)
	+
	G(f)
	,
\end{align}
where the Principal Value is no longer needed for $H$ since $z\mapsto (\partial^\alpha K)(z) \I[\|z\|>\ve]$ is a bounded map on $\R^d$ so that its convolution with any $f\in L^1(\R^d)$ is a pointwsie well-defined and bounded map on $\R^d$. In particular, this implies that $G(f)\in L^\infty(\Om)$ for any $f\in \CC^0_c(V)$. 

We observed after Proposition \ref{PropRankIntermediateRegularity} of the main paper that $f\mapsto \partial^\alpha(K*f)$ is a bounded operator on $L^2(\R^d)$; it thus follows from (\ref{eq:DerConv}) that $T:L^2(\R^d)\to L^2(\R^d)$ is also a bounded operator. Since $|(\partial^\alpha K)(z)|\lesssim \|z\|^{-d}$, the map $z\mapsto (\partial^\alpha K)(z) \I[\|z\|>\ve]$ is in $L^2(\R^d)$; Young's inequality then implies that $H:L^2(\R^d)\to L^2(\R^d)$ is a bounded operator. It follows that $G$ is a bounded operator on $L^2(\R^d)$ as well, with operator norm denoted by $\|G\|_{2}$. Letting $\supp$ denote the support, we have 
$$
\supp(G)
\subset 
\overline{ \supp\Big((\partial^\alpha K) \I[\BB^d_\ve] \Big) + \supp(f) }
=
\overline{\BB^d_\ve + \supp(f)}
.
$$
Consequently, for any $f\in \CC^0_c(V)$ we have $\supp(G (f))\subset C$. It follows that
$$
\|G(f)\|_{L^2(C)}
=
\|G(f)\|_{L^2(\R^d)}
\leq 
\|G\|_{2} \|f\|_{L^2(\R^d)}
=
\|G\|_{2} \|f\|_{L^2(V)}
\leq 
\|G\|_{2} \lambda(V)^{1/2} \|f\|_{L^\infty(V)}
$$
for all $f\in \CC^0_c(V)$, where $\lambda$ denotes the Lebesgue measure in $\R^d$. Therefore, $G$ is a bounded operator between the Banach spaces $(\CC^0_c(V),\|\cdot\|_{L^\infty(V)})$ and $(L^2(C),\|\cdot\|_{L^2(C)})$. Since $G(f)\in L^\infty(\Om)$ for all $f\in \CC^0_c(V)$, and $L^\infty(\Om)$ continuously embeds into $L^2(C)$ since $C\subset \Om$, Lemma \ref{LemContinuityTransfer} entails that $G$ is a bounded operator between $\CC^0_c(V)$ and $L^\infty(\Om)$, with operator norm denoted by $\|G\|_{\infty}$. In other words, we have 
$$
\|G(f)\|_{L^\infty(\Om)}
\leq 
\|G\|_{\infty} \|f\|_{L^\infty(V)}
,\spa 
\forall\ f\in \CC^0_c(V)
.
$$
Since $\CC^\infty_c(V)$ is dense in $\CC^0_c(V)$, we have 
%	It follows that
$$
%	B
%	:=
\sup_{\substack{f\in\CC^\infty_c(V)\\ \|f\|_{L^\infty(V)}\leq 1}} \|G(f)\|_{L^\infty(\Om)}
= 
\|G\|_{\infty}
<
\infty 
.
$$
Arguing along the lines described after Proposition \ref{PropRankIntermediateRegularity} of the main paper, a simple application of Lebesgue's dominated convergence theorem entails that $G(f)$ is a continuous map on $\Om$ for all $f\in\CC^\infty_c(V)$ since any such $f$ is Lipschitz on $\R^d$. This implies that
$$
\|G(f)\|_{L^\infty(\Om)}
=
\sup_{x\in\Om} |(Gf)(x)|
,\spa 
\forall\ f\in\CC^\infty_c(V)
.
$$
Consequently, we have
$$
%	\sup_{\substack{f\in\CC^\infty_c(V)\\ \|f\|_{L^\infty(V)}\leq 1}} \|T_2 f\|_{L^\infty(\Om)}
%	=
\sup_{x\in\Om} 
\sup_{\substack{f\in\CC^\infty_c(V)\\ \|f\|_{L^\infty(V)}\leq 1}} |(G f)(x)|
=
\|G\|_{\infty}
%	=
%	\sup_{\substack{f\in\CC^\infty_c(V)\\ \|f\|_{L^\infty(V)}\leq 1}} \|T_2 f\|_{L^\infty(\Om)}
%	=
%	\|T_2\|_{\infty,\infty}
%	<
%	\infty 
.
$$
Fix $x\in\Om$ such that $\BB^d_\ve(x)\subset V$---such an $x$ always exists by taking $\ve$ small enough in the definition of $C$ since $V$ is the closure of an open set. For all $\delta\in (0,\ve)$, define $V_\delta:=\BB^d_\ve(x)\sm \BB^d_\delta(x)$. Then 
$$
\sup_{\substack{f\in\CC^\infty_c(V_\delta)\\ \|f\|_{L^\infty(V_\delta)}\leq 1}} |(G f)(x)|
\leq 
\sup_{\substack{f\in\CC^\infty_c(V)\\ \|f\|_{L^\infty(V)}\leq 1}} |(G f)(x)|
=
\|G\|_{\infty}
.
$$
Note that the limit as $\eta\downarrow 0$ in the definition of $(G f)(x)$ disappears when $f$ is supported in $V_\delta$; indeed, we have
$$
(G f)(x)
=
\int_{\BB^d_\ve(x)\sm \BB^d_\delta(x)} (\partial^\alpha K)(x-z) f(z)\, dz 
,\spa 
\forall\ f\in \CC^\infty_c(V_\delta)
.
$$ 
In particular, these integrals are not singular since $\|x-z\|\geq \delta$; the map $z\mapsto (\partial^\alpha K)(x-z)$ is thus integrable on $V_\delta$. For all $\delta\in (0,\ve)$, this yields 
$$
\|G\|_{\infty}
\geq 
\sup_{\substack{f\in\CC^\infty_c(V_\delta)\\ \|f\|_{L^\infty(V_\delta)}\leq 1}} |(G f)(x)|
=
\int_{\BB^d_\ve(x)\sm\BB^d_\delta(x)} |(\partial^\alpha K)(x-z)|\, dz 
=
\int_{\BB^d_\ve\sm \BB^d_\delta} |(\partial^\alpha K)(z)|\, dz 
.
$$
We already observed after Proposition \ref{PropRankIntermediateRegularity} of the main paper that $(\partial^\alpha K)(z)=\|z\|^{-d} p_\alpha(z/\|z\|)$ for some polynomial $p_\alpha:\R^d\to\R$. Letting $\sigma$ denote the surface area measure on $\SS^{d-1}$, the spherical change of variables in $\R^d$ provides 
$$
\int_{\BB^d_\ve\sm \BB^d_\delta} |(\partial^\alpha K)(z)|\, dz 
=
\int_{\delta}^\ve r^{d-1} \bigg( \int_{\SS^{d-1}} \frac{1}{r^d} |p_\alpha(u)|\, d\sigma(u)\bigg)\, dr 
=
c_\alpha \log(\ve/\delta)
,
$$
with $0<c_\alpha = \int_{\SS^{d-1}} |p_\alpha(u)|\, d\sigma(u) < \infty$ since $p_\alpha$ is non-vanishing and continuous on $\R^d$. In particular, we have 
$$
0 
<
c_\alpha \log(\ve/\delta)
\leq 
\|G\|_{\infty}
,\spa 
\forall\ \delta\in (0,\ve)
.
$$
Taking $\delta\downarrow 0$ yields $\|G\|_{\infty}=\infty$, which is a contradiction.
\end{Proof}

\subsection{Proof of Theorem \ref{TheorRankRegularityOdd}}

To prove Theorem \ref{TheorRankRegularityOdd}, we will need the following propositions, known in the literature as \emph{elliptic regularity} results, for the usual Laplacian $(-\Delta)$ and the fractional Laplacian $(-\Delta)^{1/2}$.\\

%	\begin{Defin}
%		Fix $d\geq 1$ and $\Om\subset\R^d$ be an open and bounded set. Let $g\in L^2(\Om)$ and $u\in H^1(\Om)$. We say that $u$ satisfies $-\Delta u = f$ in the weak sense in $\Om$ if 
%		$$
%		\int_{\Om} \ps{\nabla u(x)}{\nabla \varphi(x)}\, dx 
%		=
%		\int_{\Om} f(x)\varphi(x)\, dx
%		,\spa \forall\ \varphi\in\CC^\infty_c(\Om)
%		.
%		$$
%		%	for any $$.
%	\end{Defin}

%The following proposition is Theorem 2 of \S 6.3.1 in \cite{Evans1998}, and will play a crucial role in our proofs.
%% can be found in \cite{Evans1998}, p. 314.
%
%\begin{Prop}[Elliptic regularity, I]\label{Prop:EllipticRegHk}
%	Let $d\geq 1$, and $\Om\subset\R^d$ be an open and bounded set. Let $g\in H^k(\Om)$ for some $k\in\N$, and $u\in H^1(\Om)$ be such that $-\Delta u = g$ in the weak sense in $\Om$. Then $u\in H^{k+2}_\loc(\Om)$, i.e $u\in H^{k+2}(V)$ for any open subset $V\subset\Om$ such that $\overline{V}\subset \Om$. In particular, $u$ admits weak derivatives of order $k+2$ in $\Om$, and we have $-\Delta u=f$ almost everywhere in $\Om$, where $\Delta u=\sum_{i=1}^d \partial_i^2 u$ and $\partial_1^2 u, \ldots, \partial_d^2 u$ are weak derivatives of $u$.
%\end{Prop}

The following proposition is Corollary $2.17$ in \cite{RosOtonRegulairty}. 

\begin{Prop}\label{PropSchauder}
Let $d\geq 1$. Fix $\beta\in (0,1)$ and $k\in\N$, and let $f\in \CC^{k,\beta}(\BB^d)$. If $u\in H^1(\BB^d)\cap L^\infty(\BB^d)$ satisfies $-\Delta u = f$ in $\D(\BB^d)'$, then $u\in \CC^{k+2,\beta}_\loc(\BB^d)$.
\end{Prop}

\vspace{3mm}

In the proof of Theorem \ref{TheorRankRegularityOdd} we will use a straightforward generalization of Proposition \ref{PropSchauder}, stated in the next corollary	.

\begin{Corol}\label{CorolSchauderGeneral}
Let $d\geq 1$ and $\Omega\subset\R^d$ an open subset. Fix $\beta\in (0,1)$ and $k\in\N$, and let $f\in \CC^{k,\beta}_\loc(\Om)$. If $u\in H^1_\loc(\Om)\cap L^\infty_\loc(\Om)$ satifies $-\Delta u = f$ in $\D(\Om)'$, then $u\in \CC^{k+2,\beta}_\loc(\Om)$.
\end{Corol}

\begin{Proof}{Corollary \ref{CorolSchauderGeneral}}
Let $x_0\in\Om$ and $r>0$ be such that $\overline{\BB^d_r(x_0 )}\subset \Om$. For any $x\in \BB^d$, let $\tilde{u}(x):=u((x-x_0)/r)$ and $\tilde{f}(x):=r^{-2} f((x-x_0)/r)$. Since $\Delta u = f$ in the weak sense in $\BB^d_r(x_0)$, a direct computation entails that $\Delta \tilde{u} = \tilde{f}$ in the weak sense in $\BB^d$. Since $u\in H^1(\BB^d_r(x_0))\cap L^\infty(\BB^d_r(x_0))$ and $f\in \CC^{k,\beta}(\BB^d_r(x_0))$, we have $\tilde{u}\in H^1(\BB^d)\cap L^\infty(\BB^d)$, and $\tilde{f}\in \CC^{k,\beta}(\BB^d)$. It follows from Proposition \ref{PropSchauder} that $\tilde{u}\in \CC^{k+2,\beta}(\BB^d)$. This implies that $u\in \CC^{k+2,\beta}(\BB^d_r(x_0))$. Now let $V\subset \Om$ be an open subset such that $\overline{V}\subset\Om$. Then $V$ can be covered by a finite number of balls of the form $\BB^d_r(x_0)$ with $\overline{\BB^d_r(x_0)}\subset\Om$. Since $u$ is of class $\CC^{k+2,\beta}$ on each one of these balls, we have $u\in \CC^{k+2,\beta}(V)$. We conclude that $u\in \CC^{k+2,\beta}_\textrm{ loc}(\Om)$.
\end{Proof}

\vspace{3mm}

We now prove an analogue of Corollary \ref{CorolSchauderGeneral} for the fractional Laplacian $(-\Delta)^{1/2}$.

\begin{Prop}\label{PropSchauderFrac}
Let $d\geq 2$ be even and $\Omega\subset\R^d$ an open subset. Fix $\beta\in (0,1)$ and $k\in\N$, and let $f\in\CC^{k,\beta}_\loc(\Om)$. If $u\in L^\infty_\loc(\Om)$ satifies $(-\Delta)^{1/2} u = f$ in $\D(\Om)'$, then $u\in \CC^{k+1,\beta}_\loc(\Om)$.
\end{Prop}

\begin{Proof}{Proposition \ref{PropSchauderFrac}.}
Let $U_1\subset U_2$ be bounded subsets of $\Om$ such that $\overline{U_1}\subset U_2$ and $\overline{U_2}\subset\Om$. Let $\eta\in\CC^\infty_c(\Om)$ be such that $0\leq \eta\leq 1$ and $\eta=1$ on $U_2$. Since $f$ is of class $\CC^{k,\beta}$ on the support of $\eta$, then $\eta f\in \CC^{k,\beta}(\R^d)$. Let $u_0\in L^\infty(\R^d)$ be such that $(-\Delta)^{1/2} u_0 = \eta f$ in $\R^d$. To see that such a $u_0$ always exists, it is enough to consider the (compactly supported) probability measures $Q^+$ and $Q^-$ with density $p_{Q^+} = c_+\ (\eta f)^+$ and $p_{Q^-} = c_-\ (\eta f)^-$, for normalizations constants constant $c_+$ and $c_-$, where $(\eta f)^+$ and $(\eta f)^-$ denote the positive and negative part of $\eta f$ respectively (so that $\eta f=(\eta f)^+ - (\eta f)^-$). Letting
$$
u_0
:=
c_+^{-1} \gamma_d\ (-\Delta)^{\frac{d-2}{2}}(\nabla\cdot F_{Q^+}^\g) - c_-^{-1} \gamma_d\ (-\Delta)^{\frac{d-2}{2}}(\nabla\cdot F_{Q^-}^\g) 
$$
yields 
$$
(-\Delta)^{\frac12} u_0 
= 
c_+^{-1} \LL_d(F_{Q^+}^\g) - c_-^{-1} \LL_d(F_{Q^-}^\g) 
= 
c_+^{-1} p_{Q^+} - c_-^{-1} p_{Q^-}
=\eta f
$$ 
in $\D(\R^d)'$, by Theorem \ref{TheorEDPDistributions} of the main paper. It remains to observe that $u_0$ is bounded over $\R^d$, which follows from Proposition \ref{PropRankIntermediateRegularity} and the fact that $\eta f$ is bounded on $\R^d$. Consequently, we have
$$
(-\Delta)^{\frac12}u_0 = \eta f = f = (-\Delta)^{\frac12} u
$$
on $U_2$. Proposition 2.22 in \cite{Silvestre2007}---and the comments thereafter---yields $u_0-u\in \CC^\infty(U_2)$. In particular, we have $u_0-u\in\CC^{k+1,\beta}(U_1)$. Assume that $u_0\in \CC^{k+1,\beta}(U_1)$ as well. We then have $u\in \CC^{k+1,\beta}(U_1)$ which, since $U_1$ was arbitrary in the first place, yields $u\in\CC^{k+1,\beta}_\loc(\Om)$ as was to be shown. Therefore, it remains to prove that $u_0\in\CC^{k+1,\beta}(U_1)$.
%		; when $d\geq 2$ is even, this follows from Proposition \ref{PropRankIntermediateRegularity} and the fact that $\eta f$ is bounded on $\R^d$, whereas, when $d\geq 3$ is odd, this follows by writing 
%		$$
%		(-\Delta)^{\frac{d-2}{2}} (\nabla\cdot F_{Q^\pm}^\g)
%		=
%		(-\Delta)^{1/2} (-\Delta)^{\frac{d-3}{2}} (\nabla \cdot F_{Q^\pm}^\g)
%		,
%		$$
%		noticing that $(-\Delta)^{\frac{d-3}{2}} (\nabla \cdot F_{Q^\pm}^\g)$ 
First assume that $k=0$. Since $\eta f\in\CC^{0,\beta}(\R^d)$ and $u_0\in L^\infty(\R^d)$ with $(-\Delta)^{1/2} u_0 = \eta f$, Proposition 2.8 of \cite{Silvestre2007} yields $u_0\in \CC^{1,\beta}(\R^d)$. Now assume that $k\geq 1$. Since $(-\Delta)^{1/2}u_0 = \eta f$ in $\D(\R^d)'$, we have $-\Delta u_0 = (-\Delta)^{1/2}(\eta f)$ in $\D(\R^d)'$. Because $\eta f\in \CC^{k,\beta}(\R^d)$, Proposition 2.7 of \cite{Silvestre2007} entails that $(-\Delta)^{1/2}(\eta f)\in \CC^{k-1,\beta}(\R^d)$. Also notice that we already proved that $u_0\in \CC^{1,\beta}(\R^d)$, so that $u_0\in H^1_\loc(\R^d)\cap L^\infty_\loc(\R^d)$. Consequently, Corollary \ref{CorolSchauderGeneral} entails that $u_0\in \CC^{k+1,\beta}_\loc(\R^d)$. In particular, we have $u_0\in \CC^{k+1,\beta}(U_1)$, which concludes the proof.
\end{Proof}

\vspace{3mm} 

We now turn to the proof of Theorem \ref{TheorRankRegularityOdd}.

\begin{Proof}{Theorem \ref{TheorRankRegularityOdd}} 

Because $P$ is non-atomic on $\Om$ recall that $\nabla h_P = \Fpg$ over $\Om$; see Definition \ref{DefinObjectiveFunction} and the comments thereafter. Since $\nabla\cdot \Fpg = \Delta h_P$, Theorem \ref{TheorEDPDistributions} entails that the equality $f_\Om=-\gamma_d\ (-\Delta)^{(d+1)/2} h_P
$ holds in $\SC(\R^d)'$. The fact that $f_\Om\in\CC^{k,\beta}_\loc(\Om)$ further implies that $f_\Om\in L^p_\loc(\Om)$ for any $p>n$, so that  $\Fpg\in \CC^{d-1}(\Om)$ by Proposition \ref{PropRankIntermediateRegularity}. In particular, we have $h_P\in\CC^d(\Om)$. Because of the different nature of $\LL_d$ when $d$ is odd or even, we split the rest of the proof and consider the case where (A) $d$ is odd, or (B) $d$ is even separately.

(A) Let us assume that $d\geq 3$ is odd. Since $h_P\in \CC^d(\Om)$ then in particular $(-\Delta)^{(d-1)/2}h_P$ belongs to $H^1_\loc(\Om)\cap L^\infty_\loc(\Om)$. Because $f_\Om\in \CC^{k,\beta}_\loc(\Om)$ and
$
f_\Om 
=
-\gamma_d\ (-\Delta) \big( (-\Delta)^{(d-1)/2} h_P\big)
$ holds in $\D(\R^d)'$, hence also in $\D(\Om)'$,
Corollary \ref{CorolSchauderGeneral} entails that $(-\Delta)^{(d-1)/2} h_P \in \CC^{k+2,\beta}_\loc(\Om)$. Since $(-\Delta)^{(d-1-2j)/2} h_P\in H^1_\loc(\Om)\cap L^\infty_\loc(\Om)$ for all $j=0,1,\ldots, (d-1)/2$, repeating the argument $(d+1)/2$ times yields $h_P\in \CC^{k+d+1,\beta}_\loc(\Om)$ and implies that $\Fpg\in\CC^{k+d,\beta}_\loc(\Om)$. In particular, we have $\Fpg \in \CC^{d}(\Om)$ which implies that $\LL_d(\Fpg)$ is continuous over $\Om$. Because $f_\Om$ is also continuous over $\Om$ and $f_\Om = \LL_d (\Fpg)$ holds in $\D(\Om)'$, the equality holds pointwise in $\Om$.

%Define $h_0=-f_\Om$ and
%$
%h_j:=\gamma_d\ (-\Delta)^{\frac{d+1}{2}-j} h_P
%$
%for all $j\in\{1, \ldots, (d+1)/2\}$.
%Let us show recursively that $h_j\in \CC^{k + 2j,\beta}_\loc(\Om)$ for all $j\in\{0,\ldots, (d+1)/2\}$. For $j=0$, this follows from the fact that $f_\Om$ belongs to $\CC^{k,\beta}_\loc(\Om)$. Fix $j\in\{0,\ldots, (d-1)/2\}$ and assume that $h_j\in \CC^{k+2j,\beta}_\loc(\Om)$, by induction. We need to show that $h_{j+1}$ belongs to $\CC^{k+2j+2,\beta}_\loc(\Om)$. To that end, observe that $h_{j+1}$ belongs to $H^1_\loc(\Om)\cap L^\infty_\loc(\Om)$. Indeed, the fact that $h_P\in\CC^d(\Om)$ entails that 
%$$
%h_{j+1}\in\CC^{2j+1}(\Om)
%\subset 
%\CC^1(\Om)
%\subset 
%H^1_\loc(\Om)\cap L^\infty_\loc(\Om)
%.
%$$
%Because $-\Delta h_{j+1} = h_j$ in $\SC(\R^d)'$ with $h_j\in\CC^{k+2j,\beta}_\loc(\Om)$ and $h_{j+1}\in H^1_\loc(\Om)\cap L^\infty_\loc(\Om)$, Corollary \ref{CorolSchauderGeneral} yields that $h_j\in\CC^{k+2j+2,\beta}_\loc(\Om)$. It follows by induction that $h_P=h_{(d+1)/2}$ belongs to $\CC^{d+k+1,\beta}_\loc(\Om)$. Consequently, $\Fpg=\nabla h_P$ belongs to $\CC^{d+k,\beta}_\loc(\Om)$. Because $f_\Om$ is continuous over $\Om$ and $\Fpg$ is of class $\CC^d$ on $\Om$, the fact that $f_\Om = \LL_d (\Fpg)$ in $\SC(\R^d)'$ entails that equality actually holds pointwise over $\Om$.
%\end{Proof}
%
%\begin{Proof}{Theorem \ref{TheorRankRegularityEven}}
(B) Assume that $d\geq 2$ is even. Since $h_P\in \CC^d(\Om)$, we have $(-\Delta)^{d/2} h_P\in L^\infty_\loc(\Om)$. As a consequence, since $f_\Om\in \CC^{k,\beta}_\loc(\Om)$ and $f_\Om = -\gamma_d\ (-\Delta)^{1/2}\big((-\Delta)^{d/2} h_P\big)$, Proposition \ref{PropSchauderFrac} entails that $(-\Delta)^{d/2} h_P\in\CC^{k+1,\beta}_\loc(\Om)$. The fact that $h_P\in\CC^d(\Om)$ further entails that $(-\Delta)^{\frac{d-2j}{2}} h_P\in H^1_\loc(\Om)\cap L^\infty_\loc(\Om)$ for all $j=0,1,\ldots, d/2$. Applying the same inductive reasoning as in Part (A) thus provides $h_P\in \CC^{k+1+d,\beta}_\loc(\Om)$ and implies that $\Fpg\in \CC^{k+d,\beta}_\loc(\Om)$. In particular $\Fpg\in \CC^{d,\beta}_\loc(\Om)$, so that $(-\Delta)^{(d-2)/2} (\nabla\cdot \Fpg)\in \CC^{1,\beta}_\loc(\Om)$. Proposition \ref{PropLaplaceExplicit} (i) of the main paper then entails that $\LL_d(\Fpg)=\gamma_d\ (-\Delta)^{1/2} \big((-\Delta)^{(d-2)/2}(\nabla\cdot \Fpg)\big)$ is continuous over $\Om$. Because $f_\Om$ is also continuous over $\Om$ and the equality $f_\Om=\LL_d(\Fpg)$ holds in the sense of distributions, it also holds pointwise in $\Om$.
\end{Proof}

\section{Proofs for Section \ref{sec:Localisation}}\label{Appendix:Localisation}

\begin{Proof}{Proposition \ref{PropLocalOdd}}
Theorem \ref{TheorEDPDistributions} entails that
$$
\int_{\R^d} \psi(x)\, dP(x)
=
\int_{\R^d} \ps{\Fpg(x)}{ (\LL_d^*\psi)(x)}\, dx 
,
$$
and 
$$
\int_{\R^d} \psi(x)\, dQ(x)
=
\int_{\R^d} \ps{F_Q^\g(x)}{ (\LL_d^* \psi)(x)}\, dx 
$$
for any $\psi\in\SC(\R^d)$. In particular, these equalities hold for any $\psi\in\CC^\infty_c(\Om)$. Because $(d-1)/2$ is an integer, $\LL_d^*=\gamma_d\ \nabla (-\Delta^{(d-1)/2})$ is a (non-fractional) differential operator. In particular, $\LL_d^*\psi$ is also supported in $\Om$. Because $\Fpg=R_Q$ over $\Om$, we have 
$$
\int_{\Om} \psi(x)\, dP(x)
=
\int_{\Om} \psi(x)\, dQ(x)
,\spa \forall\ \psi\in\CC^\infty_c(\Om)
.
$$
It follows that $P(E)=Q(E)$ for any Borel subset $E\subset\Om$.
\end{Proof}

\vspace{3mm} 

\begin{Proof}{Proposition \ref{PropLocalEven}} 
Theorem \ref{TheorEDPDistributions} entails that 
$$
P
=
\gamma_d\ (-\Delta)^{\frac{d-1}{2}}(\nabla\cdot \Fpg)
,\spa 
\text{and} 
\spa 
Q
=
\gamma_d\ (-\Delta)^{\frac{d-1}{2}}(\nabla\cdot F_Q^\g)
$$
in $\SC(\R^d)'$. Because any probability measure belongs to the domain $\SC_{1/2}(\R^d)'$ of $(-\Delta)^{1/2}$---see Section~\ref{sec:IntroLaplace}---we have 
$$
(-\Delta)^{\frac12}P
=
\gamma_d\ (-\Delta)^{\frac{d}{2}}(\nabla\cdot \Fpg)
,\spa 
\text{and} 
\spa 
(-\Delta)^{\frac12} Q
=
\gamma_d\ (-\Delta)^{\frac{d}{2}}(\nabla\cdot F_Q^\g)
$$
in $\SC(\R^d)'$. Since $d$ is even, then $(-\Delta)^{d/2}$ only involves the usual Laplacian $(-\Delta)$ taken $d/2$ times. Because $\Fpg=F_Q^\g$ on $\Om$, we have 
$$
(-\Delta)^{\frac{d}{2}}(\nabla\cdot F_P^\g)
=
(-\Delta)^{\frac{d}{2}}(\nabla\cdot F_Q^\g)
\spa 
\text{on}\quad \Om 
$$
in $\SC(\R^d)'$. It follows that $P=Q$ on $\Om$ and $(-\Delta)^{1/2}P = (-\Delta)^{1/2}Q$. As a consequence, recalling that any probability measure belongs to the negative Sobolev space $H^{-\beta}(\R^d)$ as soon as $\beta>d/2$ (see Section~\ref{sec:Sobolev}), Theorem 1.2 in \cite{GhoshEtAl2020} yields $P=Q$ on $\R^d$.
\end{Proof}

\vspace{3mm}

\begin{Proof}{Proposition \ref{ExplicitSpherOdd}}
Since $P$ is spherically symmetric, Proposition 2.2 (i) in \cite{GirStu2017} entails that there exists a map $g:[0,\infty)\to [0,\infty)$ such that $\Fpg(x)=g(\|x\|) x/\|x\|$ for all~$x\neq 0$. Theorem \ref{TheorRankRegularityOdd} implies that $g\in\CC^d((0,\infty))$ and 
$$
f_P(x) 
=
\gamma_d\ (-\Delta)^{\frac{d-1}{2}}\nabla\cdot\Big(g(\|x\|)\frac{x}{\|x\|}\Big) 
,\spa 
\forall\ x\neq 0 
.
$$
Define the operators $D,T:\CC^1\to\CC^0$ by 
$$
(D\varphi)(s) 
:=
\varphi'(s)
\spa 
\text{and}
\spa 
(T\varphi)(s)
:= 
\varphi'(s) + (d-1)\frac{\varphi(s)}{s}	
.
%	(T_2\varphi)(s) 
%	:=
%	\varphi''(s) + (d-1)\frac{\varphi'(s)}{s}
$$
A straightforward computation provides
$$
\nabla\cdot \Big(g(\|x\|)\frac{x}{\|x\|}\Big) 
=
(Tg)(\|x\|)
,\spa 
\forall\ x\in\R^d\sm\{0\} 
.
$$
This implies that, for any differentiable map $x\mapsto \psi(\|x\|)$, we have
$$
\Delta\big(\psi(\|x\|)\big) 
=
\nabla\cdot\big(\nabla(\psi(\|x\|))\big)
=
\nabla\cdot\Big( (D\psi)(\|x\|)\frac{x}{\|x\|}\Big)
=
\big((T\circ D)\psi\big)(\|x\|)
,\spa 
\forall\ x\in\R^d\sm\{0\}
.
$$
For the sake of simplicity, we will abbreviate $T\circ D$ by $TD$. Letting $f:[0,\infty)\to[0,\infty)$ be such that $f_P(x)=f(\|x\|)$ for all $x\in\R^d$, and $\bar{f}:=\gamma_d^{-1} (-1)^{(d-1)/2} f$, the equality $f_P = \gamma_d\ (-\Delta)^{(d-1)/2}(\nabla\cdot \Fpg)$ thus reads
$$
\big((TD)^{\frac{d-1}{2}} T g\big)(s) =  \bar{f}(s)
,\spa 
\forall\ s>0 
.
$$
For all $k=0,\ldots, (d-1)/2$, define $u_k:=(TD)^k Tg$. We will obtain the desired expression for $g$ by solving the equation $Tg = u_0$ in $g$, while $u_0$ will be obtained by iteratively solving $(TD)u_k = u_{k+1}$ in $u_k$, with initial step $u_{(d-1)/2}=\bar{f}$. Multiplying the equality $(T g)(s)=u_0(s)$ by $s^{d-1}$ entails that $(s^{d-1} g)'(s) = s^{d-1}u_0(s)$ for all $s>0$. Since $g\in\CC^d((0,\infty))$, then $u_0$ is continous. In particular, the map $r\mapsto r^{d-1} u_0(r)$ is integrable on each interval, which implies that for some constant $c$, we have
$$
g(s)
=
\frac{1}{s^{d-1}} \int_0^s t^{d-1} u_0(t)\, dt + \frac{c}{s^{d-1}} 
,\spa\forall\ s>0 
.
$$
Since $g$ and $u_0$ are bounded near $s=0$, we have $c=0$. Let us show by induction that for all $k\in\{0,\ldots, (d-1)/2\}$ there exists a map $q_k: (0,1]\to\R$ that does not depend on $P$ such that 
\begin{equation}\label{eq:InductG}
	g(s)
	=
	\sum_{i=0}^{k-1} \Big( \E\Big[\frac{1}{\|Z\|^{2i+1}}\Big] w_i s^{2i+1} \Big)   + \int_0^s t^{2k} q_k\Big(\frac{t}{s}\Big) u_k(t)\, dt
	,\spa\forall\ s>0 
	,
\end{equation}
where, for all $i\in\{0,\ldots, (d-1)/2\}$, 
$$
w_i 
=
(-1)^i 2^{2i+1} \frac{\Gamma(\frac{d+1}{2}) \Gamma(\frac{2i+1}{2}) }{ \sqrt{\pi}\ \Gamma(\frac{d-2i-1}{2}) } 
\int_0^1 t^{2i} q_i(t)\,dt 
.
$$
In addition, we have $q_0(y) = y^{d-1}$ for all $y\in (0,1]$, and for all $k\in\{1,\ldots, (d-1)/2\}$ we have
$$
q_k(y)
=
%	\alpha_{0}^{(k)} y^{d-1}
\sum_{i=0}^{k} \alpha_{i}^{(k)} y^{d-1-2i} 
+
\sum_{i=1}^k \beta_{i}^{(k)} \frac{1}{y^{2i-1}}
,\spa\forall\ y\in (0,1]
,
$$
for constants $\alpha_{0}^{(k)}> 0$ and $\beta_i^{(k)}>0$, $i=1,\ldots, k$. We already showed that this holds for $k=0$ with $q_0(y)=y^{d-1}$. Let us assume this holds for some $k\in\{0,\ldots, (d-3)/2\}$ and show that it holds for $k+1$. Notice that  $(TD)u_k = u_{k+1}$. To solve this equation in $u_k$, let us decompose the problem into $T\theta_k = u_{k+1}$ and $\theta_k=Du_k$. Reasoning as before yields
$$
\theta_k(r)
=
\frac{1}{r^{d-1}} \int_0^r y^{d-1} u_{k+1}(y)\, dy
,\spa\forall\ r>0 
,
$$
while the equality $Du_k = \theta_k$ gives
$$
u_k(t) 
=
u_k(0) + \int_0^t \theta_k(r)\, dr
,\spa\forall\ t>0 
.
$$
Combining the expressions established for $\theta_k$ and $u_k$, Fubini's theorem provides
$$
u_k(t)
=
u_k(0)
+
\frac{1}{d-2} \int_0^t r\big(1-(r/t)^{d-2}\big)u_{k+1}(r)\, dr 
,\spa \forall\ t>0 
.
$$
Plugging the last expression of $u_k$ in the expression (\ref{eq:InductG}) for $k$, which holds by our inductive step, we find
\begin{eqnarray*}
	\lefteqn{
		g(s)
		=
		\sum_{i=0}^{k-1} \Big( \E\Big[\frac{1}{\|Z\|^{2i+1}}\Big] w_i s^{2i+1} \Big)  
		+
		u_k(0) \int_0^s t^{2k} q_k\Big(\frac{t}{s}\Big)\, dt 
	}
	\\[2mm]
	&&
	\hspace{15mm}
	+
	\frac{1}{d-2} \int_0^s \Big( t \int_t^s (1-(t/r)^{d-2}) r^{2k} q_k\Big(\frac{r}{s}\Big)\,dr\Big)  u_{k+1}(t)\, dt
	.
\end{eqnarray*}
Observe that $(-1)^k u_k(0)=(-\Delta)^k(\nabla\cdot \Fpg)(0)$, and recall that $\nabla\cdot(x/\|x\|)=(d-1)\|x\|^{-1}$. Therefore, letting $Z$ be a random $d$-vector with law $P$, Proposition \ref{PropRankIntermediateRegularity} and Corollary \ref{CorolFTRiesz} of the main paper yield 
$$
u_k(0)
=
(d-1)\ (-1)^k (-\Delta)^k \Big( \E\Big[\frac{1}{\|Z-x\|}\Big]\Big)_{|x=0}
=
(-1)^k 2^{2k+1} \frac{\Gamma(\frac{d+1}{2}) \Gamma(\frac{2k+1}{2}) }{ \sqrt{\pi}\ \Gamma(\frac{d-2k-1}{2}) } 
\E\Big[\frac{1}{\|Z\|^{2k+1}}\Big]
.
$$
Changing the variable $t':=t/s$ yields 
$$
\sum_{i=0}^{k-1} \Big( \E\Big[\frac{1}{\|Z\|^{2i+1}}\Big] w_i s^{2i+1} \Big)
+
u_k(0) \int_0^s t^{2k} q_k\Big(\frac{t}{s}\Big)\, dt 
=
\sum_{i=0}^{k} \E\Big[\frac{1}{\|Z\|^{2i+1}}\Big] w_i s^{2i+1}
.
$$
After changing the variable $r':=r/s$, we find
$$
t \int_t^s (1-(t/r)^{d-2}) r^{2k} q_k\Big(\frac{r}{s}\Big)\,dr
=
t^{2k+2} \Big(\frac{s}{t}\Big)^{2k+1} \Big( \int_{\frac{t}{s}}^1 r^{2k}q_k(r)\, dr - \Big(\frac{t}{s}\Big)^{d-2} 
\int_{\frac{t}{s}}^1 \frac{ r^{2k} q_k(r) }{ r^{d-2} }\, dr  \Big)
.
$$
We deduce that 
$$
g(s)
=
\sum_{i=0}^{k} \Big( \E\Big[\frac{1}{\|Z\|^{2i+1}}\Big] w_i s^{2i+1} \Big)
+
\int_0^s t^{2(k+1)} q_{k+1}\Big(\frac{t}{s}\Big) u_{k+1}(t)\, dt
,\spa \forall\ s>0 
,
$$
with 
%	$$
%	m_{k+1}(s)
%	=
%	m_k(s)
%	+
%	c_k \int_0^s t^{2k} q_k\Big(\frac{t}{s}\Big)\, dt
%	=
%	m_k(s)
%	+
%	c_k s^{2k+1} \int_0^1 t^{2k} q_k(t)\, dt
%	$$
%	and 
$$
q_{k+1}(y)
=
\frac{1}{d-2}\frac{1}{y^{2k+1}} \Big( \int_y^1 r^{2k}q_k(r)\, dr - y^{d-2}\int_y^1 \frac{r^{2k} q_k(r)}{r^{d-2}}\, dr \Big)
,\spa \forall\ y\in (0,1]
.
$$
Since $q_k$ has no term of order $-(2k+1)$ or $d-3-2k$, then $r^{2k}q_k(r)$ and $r^{2k-(d-2)}q_k(r)$ are Laurent polynomials with no term of order $-1$. It follows that $q_{k+1}$ is a Laurent polynomial as well. In addition, any term $r^\ell$ in $q_k(r)$ gives rise to the following terms for $q_{k+1}(r)$:
$$
\frac{1}{d-2}\bigg\{  \Big(\frac{1}{2k-(d-3)+\ell}-\frac{1}{2k+1+\ell}\Big)y^\ell - \frac{y^{d-3-2k}}{2k-(d-3)+\ell} + \frac{1}{(2k+1+\ell)y^{2k+1}}\bigg\}
.
$$
Since $\alpha_0^{(k)}>0$ by our inductive step, we deduce that
$$
\alpha_0^{(k+1)}
=
\alpha_0^{(k)} \frac{1}{d-2} \Big(\frac{1}{2k-(d-3)+(d-1)}-\frac{1}{2k+1+(d-1)}\Big)
=
%	\frac{d-2}{(2k+2)(2k+d)} \alpha_0^{(k)}
%	=
\frac{\alpha_0^{(k)}}{(2k+2)(2k+d)} 
>
0
.
$$
In addition, it is clear that $\beta_i^{(k)}>0$ for all $i=1,\ldots, k$. This concludes the proof by induction. In particular, we have
\begin{align*}
	g(s)
	&=
	\sum_{i=0}^{\frac{d-3}{2}} \Big( \E\Big[\frac{1}{\|Z\|^{2i+1}}\Big] w_i s^{2i+1} \Big) 
	+ 
	\int_0^s t^{d-1} q_{(d-1)/2}\Big(\frac{t}{s}\Big) u_{(d-1)/2}(t)\,dt
	\\[2mm]
	&=
	\sum_{i=0}^{\frac{d-3}{2}} \Big( \E\Big[\frac{1}{\|Z\|^{2i+1}}\Big] w_i s^{2i+1} \Big) 
	+ 
	(-1)^{\frac{d-1}{2}}\gamma_d^{-1} \int_0^s t^{d-1} q_{(d-1)/2}\Big(\frac{t}{s}\Big) f(t)\,dt
	.
\end{align*}
Denoting by $\sigma(\SS^{d-1})$ the surface area of $\SS^{d-1}$ and letting 
$$
p(y)
:=
\frac{(-1)^{\frac{d-1}{2}}}{\gamma_d \sigma(\SS^{d-1})} 
q_{(d-1)/2}(y)
,\spa \forall\ y\in (0,1]
,
$$
the spherical change of variable yields 
$$
g(s)
=
\sum_{i=0}^{\frac{d-3}{2}} \Big( \E\Big[\frac{1}{\|Z\|^{2i+1}}\Big] w_i s^{2i+1} \Big) 
+ 
\E\Big[ p\Big(\frac{\|Z\|}{s}\Big) \I\big[\|Z\|\leq s\big]\Big] 
,\spa \forall\ s>0
.
$$
This concludes the proof
%	In particular, the term in $r^{d-1}$ in $p_{k+1}$ is obtained from $p_k$, through the previous expression relating $p_{k+1}$ and $p_k$, solely from the term in $r^{d-1}$ in $p_k$. We thus see from the previous computation that $\alpha coefficient is strictly positive
\end{Proof}

%\subsection{Proof of Proposition \ref{PropEquivLocalMethods}} 

\vspace{3mm}

\begin{Proof}{Proposition \ref{PropEquivLocalMethods}}
(i) Because $P^*$ admits the null density on $\R^{d+1}_+$, Theorem \ref{TheorRankRegularityOdd} entails that $F_{P^*}^\g\in\CC^\infty(\R^{d+1}_+)$, so that $U\in\CC^\infty(\R^{d+1}_+)$ as well, and that $-\Delta U(x^*)=0$ for all $x^*\in\R^{d+1}_+$.\\

(ii) Let $Z$ be a random $d$-vector with law $P$ and $Z^*$ a random $(d+1)$-vector with law $P^*$. Since $P^*$ admits the null density on $\R^{d+1}_+$, Proposition \ref{PropRankIntermediateRegularity} yields
$$
U(x^*)
=
2\gamma_{d+1} \E\big[(-\Delta)^{\frac{d-2}{2}}(\nabla\cdot K_{d+1})(x^*-Z^*)\big]
,\spa \forall\ x^*\in\R^{d+1}_+
.
$$
Since $f_P\in L^p_\loc(\Om)$ for some $p>d$, then Proposition \ref{PropRankIntermediateRegularity} yields
$$
u(x)
=
\gamma_d\ \E\big[(-\Delta)^{\frac{d-2}{2}}(\nabla\cdot K_d)(x-Z)\big]
,\spa 
\forall\ x\in\Om
.
$$
Because $(\nabla\cdot K_d)(z) = (d-1)\|z\|^{-1}$ for all $z\in\R^d\sm\{0\}$ and $(\nabla\cdot K_{d+1})(z)=d\|z\|^{-1}$ for all $z\in\R^{d+1}\sm\{0\}$, Corollary \ref{CorolFTRiesz} provides 
\begin{equation}\label{eq:ExpU}
	U(x^*)
	=
	2d\gamma_{d+1} \frac{2^{d-2}\Gamma(\frac{d}{2}) \Gamma(\frac{d-1}{2})}{\Gamma(\frac12)}
	\E\Big[\frac{1}{\|x^*-Z^*\|^{d-1}}\I[Z^*\neq x^*]\Big]
	,\spa 
	\forall\ x^*\in\R^{d+1}_+
	,
\end{equation}
and
$$
u(x)
=
(d-1)\gamma_d \frac{2^{d-2} \Gamma(\frac{d-1}{2})^2}{\Gamma(\frac12)^2}
\E\Big[\frac{1}{\|x-Z\|^{d-1}}\I[Z\neq x]\Big]
,\spa 
\forall\ x\in\Om
.
$$
Writing $Z=(Z_1,\ldots, Z_d)$, the definition of $Z^*$ entails that we can take $Z^*=(Z_1,\ldots, Z_d, 0)$. Consequently, we have for all $x^*=(x,s)\in\Om\times (0,\infty)$
$$
\E\Bigg[\frac{1}{\|x^*-Z^*\|^{d-1}}\I[Z^*\neq x^*]\Bigg]
=
\E\Bigg[\frac{1}{\big(\|x-Z\|^2 + s^2\big)^{\frac{d-1}{2}}}\I[Z\neq x]\I[s\neq 0]\Bigg]
.
$$
Since the expression under the last expectation is bounded, uniformly in $s\in [0,\infty)$, by the $P$-integrable random variable~$1/\|x-Z\|^{d-1}$, Lebesgue's dominated convergence theorem entails that $U$ extends continuously on $\R^d\times [0,\infty)=\overline{\R^{d+1}_+}$. In addition, we have $U(x,0)=u(x)$ for all $x\in\Om$ provided 
$$
2d\gamma_{d+1} \frac{2^{d-2}\Gamma(\frac{d}{2}) \Gamma(\frac{d-1}{2})}{\Gamma(\frac12)}
=
(d-1)\gamma_d \frac{2^{d-2} \Gamma(\frac{d-1}{2})^2}{\Gamma(\frac12)^2}
,
$$
\ie provided
$
2 \gamma_{d+1} d\ \Gamma(d/2) \Gamma(1/2)
=
\gamma_d (d-1)\Gamma((d-1)/2)
.
$
Recalling that $\Gamma(1/2)=\sqrt{\pi}$ and $x\Gamma(x)=\Gamma(x+1)$ for all $x>0$, the definition of $\gamma_d$ and $\gamma_{d+1}$---see Section \ref{sec:PDEGeneral} of the main paper---entails that the last identity reduces to $2=2$. We deduce that $U(x,0)=u(x)$ for all $x\in\Om$.\\

(iii) Since $P$ admits the null density on $\R^{d+1}_+$, Proposition \ref{PropRankIntermediateRegularity} implies that $(\partial_{d+1} U)(x^*)$ can be computed by applying $\partial_{d+1}$ under the expectation defining $U$. Letting $x^*=(x_1^*,\ldots, x_{d+1}^*)$ and $Z^*=(Z_1^*,\ldots,Z_{d+1}^*)$, and recalling (\ref{eq:ExpU}) and the definition of $\gamma_{d+1}$ yields
\begin{align*}
	-(\partial_{d+1} U)(x^*)
	&=
	-\frac{\Gamma(\frac{d-1}{2})}{2\pi^{\frac{d+1}{2}}}  \E\Big[\partial_{d+1}\Big(\frac{1}{\|x^*-Z^*\|^{d-1}}\Big) \I[Z^*\neq x^*]\Big]
	\\[2mm]
	&=
	\frac{\Gamma(\frac{d+1}{2})}{\pi^{\frac{d+1}{2}}}  \E\Big[\frac{x^*_{d+1}-Z_{d+1}^*}{\|x^*-Z^*\|^{d+1}}  \I[Z^*\neq x^*]\Big]
	,\spa \forall\ x^*\in\R^{d+1}_+
	.
\end{align*}
Writing $Z^*=(Z_1,\ldots, Z_d, 0)$ as in the proof of (ii) above, we have for all $x^*=(x,\ve)\in\R^{d+1}_+$
$$
-(\partial_{d+1} U)(x,\ve)
=
\frac{\Gamma(\frac{d+1}{2})}{\pi^{\frac{d+1}{2}}}  \E\Big[\frac{\ve}{\|x^*-Z^*\|^{d+1}}  \I[Z^*\neq x^*]\Big]
=
\frac{\Gamma(\frac{d+1}{2})}{\pi^{\frac{d+1}{2}}} 
\E\Bigg[\frac{\ve}{\big(\|x-Z\|^2 + \ve^2\big)^{\frac{d+1}{2}}}\I[Z\neq x]\Bigg]
.
$$
%Writing $Z^*=(Z_1^*, \ldots, Z_{d+1}^*)$, we then have 
%$$
%-(\partial_{d+1} G)(\tilde{x},x_{d+1})
%=
%2\gamma_{d+1} \Gamma(d+1)
%\E\bigg[\frac{x_{d+1}-Z_{d+1}^*}{\|x-Z^*\|^{d+1}}\Ind{Z^*\neq x}\bigg]
%,\spa \forall\ (\tilde{x},x_{d+1})\in\R^d\times (0,\infty)
%.
%$$
%Let us show that 
%$$
%\lim_{x_{d+1} \downarrow 0}
%\E\bigg[\frac{x_{d+1}-Z_{d+1}^*}{\|x-Z^*\|^{d+1}}\Ind{Z\neq x}\bigg]
%=
%\frac{\pi^{\frac{d+1}{2}}}{\Gamma(\frac{d+1}{2})}
%f_P(\tilde{x})
%.
%$$
%For any $x\in\R^{d+1}_+$, we have 
Define the map $\rho:\R^d\to [0,1]$ by letting 
$$
\rho(y)
=
\frac{\Gamma(\frac{d+1}{2})}{\pi^{\frac{d+1}{2}}}  
\frac{1}{(1+\|y\|^2)^{\frac{d+1}{2}}}
,\spa 
\forall\ y\in\R^d
.
$$
We have $\rho\geq 0$, and one can further show that $\int_{\R^d} \rho(y)\, dy = 1$. Then, letting $\rho_\ve(y):=\ve^{-d} \rho(y/\ve)$ for all $y\in\R^d$ and $\ve>0$, we have
$$
-(\partial_{d+1} U)(x,\ve)
=
\E\big[\rho_\ve(x-Z)\big] 
=
\int_{\R^d} \rho_\ve(x-z)\, dP(z)
%=
%(\rho_\ve * P)(x)
,\spa 
\forall\ (x,\ve)\in\R^{d+1}_+
.
$$
Since $P$ has a bounded density $f_P$ on $\R^d$, and $f_P$ is continuous on $\Om$, a change of variables combined with Lebesgue's dominated convergence theorem yield
$$
f_P(x)
=
\lim_{\ve\downarrow 0} -(\partial_{d+1} U)(x,\ve)
,\spa 
\forall\ x\in\Om 
.
$$
Because $f\in \CC^{0,\beta}_\loc(\R^d)$, Theorem \ref{TheorRankRegularityOdd} ensures that $f_P(x)=((-\Delta)^{1/2} u)(x)$ for all $x\in\Om$, which concludes the proof.
%$$
%f_P(x)
%=
%\gamma_d\ (-\Delta)^{1/2} (-\Delta)^{\frac{d-2}{2}}(\nabla\cdot \Fpg)(x)
%=
%((-\Delta)^{1/2} g)(\tilde{x})
%,\spa \forall\ \tilde{x}\in\R^d
%,
%$$
%which concludes the proof.
\end{Proof}

\vspace{3mm}

%\section{Proofs for Section \ref{sec:ProbaOpenSets}} 

\begin{Proof}{Proposition \ref{PropLowerScaleEven}} 
Fix $\nu\in\R^d$ orthogonal to $H$. Since $P$ admits the null density on $\R^d\sm H$, Theorem \ref{TheorRankRegularityOdd} entails that $\Fpg\in \CC^\infty(\R^d\sm H)$. In addition, Proposition \ref{PropRankIntermediateRegularity} yields 
$$
\big(\partial_{\nu} (-\Delta)^{\frac{\kappa}{2}-1}(\nabla\cdot F_P^\g)\big)(x+\ve \nu)
=
\int_{\R^d} \big(\partial_{\nu} (-\Delta)^{\frac{\kappa}{2}-1}(\nabla\cdot K_d)\big)(x+\ve \nu-z)\, dP(z)
$$
for all $x\in H$ and $\ve>0$.
Recalling that $(\nabla\cdot K_d)(z)=(d-1)\|z\|^{-1}$, then Corollary \ref{CorolFTRiesz} yields for all $x,z\in H$
\begin{align*}
	\big(\partial_{\nu} (-\Delta)^{\frac{\kappa}{2}-1}(\nabla\cdot K_d)\big)(x+\ve \nu-z)
	&=
	2^{\kappa-2}(d-1)(\kappa-1) \frac{\Gamma(\frac{d-1}{2}) \Gamma(\frac{\kappa-1}{2})}{\sqrt{\pi}\ \Gamma(\frac{d-\kappa+1}{2})} 
	\ps{\frac{x+\ve \nu - z}{\|x+\ve\nu -z\|^{\kappa+1}}}{\nu}
	\\[2mm]
	&= 
	\frac{2^{\kappa}\ \Gamma(\frac{d+1}{2}) \Gamma(\frac{\kappa+1}{2})}{\sqrt{\pi}\ \Gamma(\frac{d-\kappa+1}{2})} 
	\frac{\ve}{ (\|x-z\|^2 + \ve^2)^{\frac{\kappa+1}{2}} }
	\\[2mm]
	&= 
	\frac{2^{\kappa}\ \Gamma(\frac{d+1}{2}) \Gamma(\frac{\kappa+1}{2})}{\sqrt{\pi}\ \Gamma(\frac{d-\kappa+1}{2})}	\frac{\ve^{-\kappa}}{ (\|x-z\|^2 + 1)^{\frac{\kappa+1}{2}} }
	,
\end{align*}
where we used that $\nu$ is orthogonal to $H$.
Define $\rho:\R^\kappa\to[0,\infty)$ by letting 
$$
\rho(y)
=
\frac{\Gamma(\frac{\kappa+1}{2})}{\pi^{\frac{\kappa+1}{2}}}  
\frac{1}{(1+\|y\|^2)^{\frac{\kappa+1}{2}}}
,\spa 
\forall\ y\in\R^\kappa
.
$$
We have $\rho\geq 0$, and one can further show that $\int_{\R^\kappa} \rho(y)\, dy = 1$. Then, letting $\rho_\ve(y):=\ve^{-\kappa} \rho(y/\ve)$ for all $y\in\R^\kappa$ and $\ve>0$, we have
$$
-\frac{\Gamma(\frac{d-\kappa+1}{2})}{(4\pi)^{\frac{\kappa}{2}} \Gamma(\frac{d+1}{2})}\big(\partial_{\nu} (-\Delta)^{\frac{\kappa}{2}-1}(\nabla\cdot K_d)\big)(x+\ve \nu-z)
=
\int_{\R^\kappa} \rho_\ve(x-z) f_P(z)\, dz
%=
%(\rho_\ve * P)(x)
,
$$
where, by an obvious abuse of notation, we identified $H$ with $\R^\kappa$ and $x,z$ with their representative in $\R^\kappa$.
Since $P$ has a bounded and continuous density $f_P$ on $\R^d$, a change of variables combined with Lebesgue's dominated convergence theorem yields
$$
f_P(x)
=
\lim_{\ve\downarrow 0} \int_{\R^\kappa} \rho_\ve(x-z) f_P(z)\, dz
,\spa\forall\ x\in H
,
$$
which concludes the proof.
\end{Proof}

\section{Proofs for Section \ref{sec:DepthRegions}}

The proof of Proposition~\ref{PropReIndexedRanksConv} requires the following Lemma.

\begin{Lem} \label{LemDifferenceVectUnit}
For all $v,w\in\R^d\setminus\{0\}$, we have
$$
\bigg\|\frac{v}{\|v\|}-\frac{w}{\|w\|}\bigg\|\leq 2 \min\bigg\{1,\frac{\|v-w\|}{\|v\|}, \frac{\|v-w\|}{\|w\|}\bigg\}
.
$$
\end{Lem}

\vspace{3mm}

We now proceed with the proof of Proposition \ref{PropReIndexedRanksConv}.

\begin{Proof}{Proposition \ref{PropReIndexedRanksConv}}
Fix $x\in\R^d$ and let us decompose 
\begin{eqnarray*}
	\lefteqn{
		\hspace{-5mm}
		\widetilde{F}_{P_N}^\g(x)-\widetilde{F}_P^\g(x)
		=
		\Big( w_{P_N}\big(\|F_{P_N}^\g(x)\|\big)-w_P\big(\|\Fpg(x)\|\big)\Big) \frac{F_{P_N}^\g(x)}{\|F_{P_N}^\g(x)\|}
	}
	\\[2mm]
	&&
	+
	\Big( w_P\big(\|F_{P_N}^\g(x)\|\big) - w_P\big(\|F_P^\g(x)\|\big) \Big) \frac{F_{P_N}^\g(x)}{\|F_{P_N}^\g(x)\|}
	+
	w_P(\|F_P^\g(x)\|)\bigg( \frac{F_{P_N}^\g(x)}{\|F_{P_N}^\g(x)\|} - \frac{F_{P}^\g(x)}{\|F_{P}^\g(x)\|} \bigg)
	\\[2mm]
	&&
	=: (I) + (II) + (III)
	.
\end{eqnarray*}
Let us denote by $f_P$ the density of $P$. Since $f_P$ is bounded, Lemma 2 in \cite{Mottonen97} entails that 
$$\sup_{x\in\R^d} \|F_{P_N}^\g(x)-\Fpg(x)\|\to 0$$ 
$P$-almost surely as $N\to\infty$. In particular, $\sup_{x\in\R^d}\|(II)\| \to 0$ $P$-almost surely as well provided $w_P$ is uniformly continuous on $[0,1)$. To see that $w_P$ is uniformly continuous, recall that $F_P^\g$ is invertible with inverse $Q_P^\g$ since $P$ is non-atomic and not supported on a line (see Section \ref{sec:Relabeling} of the main paper). Then, denoting the Jacobian matrix of $Q_P^\g$ at $x$ by $J_x Q_P^\g$, we have
$$
w_P(\beta)
=
P\big[Q_P^\g(\BB^d_\beta)\big]
=
\int_{\BB^d_\beta} f_P(Q_P^\g (x))\ |\det(J_x Q_P^\g)|\, dx
=
\int_0^\beta g(r)\, dr 
,\spa \forall\ \beta\in (0,1)
,
$$
where, denoting by $\sigma$ the surface area measure on $\SS^{d-1}$, we let
$$
g(r)
:=
r^{d-1} \int_{\SS^{d-1}} f_P(Q_P^\g (ru))\ |\det(J_{ru} Q_P^\g)|\, d\sigma(u)
,\spa \forall r>0 
.
$$	
Since $w_P$ is a bounded map, then $g$ is integrable on $(0,\infty)$. It follows that for all $\ve>0$, there exists $\delta>0$ such that for any Borel measurable set $A\subset (0,\infty)$ with Lebesgue measure smaller than $\delta$ we have $\int_A g(x)\, dx < \ve$. This implies that $w_P$ is uniformly continuous, hence that $\sup_{x\in\R^d} \|(II)\| \to 0$ $P$-almost surely as $N\to\infty$.

We obviously have $\sup_{x\in\R^d} \|(I)\| \leq \sup_{t\in [0,1]} |w_{P_N}(t)-w_P(t)|$.
If $P_N$ is the empirical cdf associated with a sample $Z_1,\ldots, Z_N$, define the empirical cdf $\widetilde{w}_N$ associated with the sample $\|\Fpg(Z_1)\|,\ldots, \|\Fpg(Z_N)\|$, i.e 
$$
\widetilde{w}_N(t)
=
\frac{1}{N}\sum_{j=1}^N \Ind{\|\Fpg(Z_i)\|\leq t}
,\spa\forall\ t\in [0,1]
.
$$
Letting $\Delta_N:=\sup_{x\in\R^d} \|F_{P_N}^\g(x)-\Fpg(x)\|$, we have for all $t\in [0,1]$
\begin{align*}
	w_{P_N}(t)-w_P(t)
	&=
	\frac{1}{N}\sum_{j=1}^N \Ind{\|F_{P_N}^\g(Z_i)\|\leq t} - w_P(t)
	\\[2mm]
	&\leq 
	\frac{1}{N}\sum_{j=1}^N \Ind{\|\Fpg(Z_i)\|\leq t + \Delta_N} - w_P(t)
	\\[2mm]
	&= 
	\widetilde{w}_N(t+\Delta_N)-w_P(t)
	\\[2mm]
	&=
	\big(\widetilde{w}_N(t+\Delta_N)-w_P(t+\Delta_N)\big) 
	+ 
	\big( w_P(t+\Delta_N)-w_P(t) \big)
	.
\end{align*}
The Glivenko-Cantelli Theorem entails that $\sup_{s\in [0,1]} |\widetilde{w}_N(s)-w_P(s)|\to 0$ $P$-almost surely as $N\to\infty$, and the uniform continuity of $w_P$ and the fact that $\Delta_N\to0$ $P$-almost surely as $N\to\infty$ entail that $\sup_{t\in [0,1]} |w_P(t+\Delta_N)-w_P(t)|\to 0$ $P$-almost surely as $N\to\infty$. Since a similar reasoning yields 
$$
w_{P_N}(t)-w_P(t)
\geq 
\big(\widetilde{w}_N(t-\Delta_N)-w_P(t-\Delta_N)\big) 
+ 
\big( w_P(t-\Delta_N)-w_P(t) \big) 
,
$$
for all $t\in [0,1]$, we deduce that $\sup_{x\in\R^d} \|(I)\|\to 0$ $P$-almost surely as $N\to\infty$. 

By \ref{LemDifferenceVectUnit}, and using the fact that $w_P$ is non-decreasing and takes its values in $[0,1]$, we have
\begin{align*}
	\|(III)\|
	&\leq 
	2 w_P(\|\Fpg(x)\|) \min\bigg\{1,\frac{\|F_{P_N}^\g(x)-\Fpg(x)\|}{\|\Fpg(x)\|}\bigg\}
	\\[2mm]
	&\leq 
	2 w_P(\ve)\Ind{\|\Fpg(x)\|\leq \ve}
	+
	\frac{w_P(\ve)}{\ve} \|F_{P_N}^\g(x)-\Fpg(x)\| \Ind{\|\Fpg(x)\|> \ve}
	.
\end{align*}
Since $\sup_{x\in\R^d} \|F_{P_N}^\g(x)-\Fpg(x)\|\to 0$ $P$-almost surely as $N\to\infty$, we have 
$$
\limsup_{N\to\infty}\sup_{x\in\R^d} \|\widetilde{F}_{P_N}^\g(x)-\widetilde{F}_P^\g(x)\|
\leq 
2 w_P(\ve)
.
$$
Because $\ve>0$ was arbitrary and $w_P$ is right-continuous (recall it is a distribution function), we have
$$
\limsup_{N\to\infty}
\sup_{x\in\R^d} \|\widetilde{F}_{P_N}^\g(x)-\widetilde{F}_P^\g(x)\|
\leq 
2w_P(0)
.
$$
But we have $w_P(0)=P[\|\Fpg(Z_1)\|=0]=P[\Fpg(Z_1)=0]$. Because the equation $\Fpg(z)=0$ has a unique solution, given by the geometric median of $P$---see the comments after Definition \ref{DefinQuantiles}---and $P$ is non-atomic, we have $w_P(0)=0$. It follows that 
$$
\sup_{x\in\R^d} \|\widetilde{F}_{P_N}^\g(x)-\widetilde{F}_P^\g(x)\|
\to 
0
$$
$P$-almost surely as $N\to\infty$.
\end{Proof}

\vspace{3mm}

\begin{Proof}{Proposition \ref{PropRegQuantiles}}
In Part (i), Proposition \ref{PropRankIntermediateRegularity} entails that $\Fpg$ belongs to $\CC^\ell(\R^d)$, while in Part (ii) Theorem \ref{TheorRankRegularityOdd} ensures that $\Fpg\in\CC^{d+k}(\R^d)$. In both cases, let us write $\Fpg\in \CC^j(\R^d)$ with $j\in\{\ell,d+k\}$. We also have
$$
(\partial^\alpha\Fpg)(x)
=
\E[(\partial^\alpha K)(x-Z)]
,\spa \forall\ x\in\R^d
$$
for any $\alpha\in\N^d$ with $|\alpha|\leq j$, where $Z$ is a random $d$-vector with law $P$. Because $P$ has a density with respect to the Lebesgue measure, Theorem 6.2 in \cite{KonPai1} entails that $\Fpg:\R^d\to\BB^d$ is a homeomorphism with inverse $(\Fpg)^{-1}=Q_P^\g$. In particular, $\Fpg$ is injective on $\R^d$. Assume for now that the Jacobian matrix $J_x\Fpg$ of $\Fpg$ at $x\in\R^d$ is invertible for all $z\in\R^d$. Then, the Inverse Function Theorem entails that $Q_P^\g=(\Fpg)^{-1}$ is a diffeomorphism of class $\CC^j(\BB^d)$. It remains to show that $J_x\Fpg$ is invertible at any $x\in\R^d$. For this purpose, fix $x\in\R^d$ and recall that $\partial_i \Fpg(x)=\E[(\partial_i K)(x-Z)]$ for all $i=1,2,\ldots, d$, with
$$
J_x K=\frac{1}{\|x\|}\Big(I_d-\frac{xx'}{\|x\|^2}\Big)
,\spa \forall\ x\in\R^d\sm\{0\},
$$
where $I_d$ stands for the $d\times d$ identity matrix and $x'$ denotes the transpose of $x$. Consequently, we have
$$
J_x\Fpg
=
\E\bigg[\frac{1}{\|x-Z\|}\Big(I_d-\frac{(x-Z)(x-Z)'}{\|x-Z\|^2}\Big)\Ind{Z\neq x}\bigg]
.
$$
The matrix $J_x\Fpg$ is obviously symmetric and non-negative definite. It is thus enough to show that it is positive definite. Assume, ad absurdum, that there exists $v\in\SS^{d-1}$ such that $\ps{v}{(J_x\Fpg) v} = 0$, \ie 
$$
\E\bigg[\frac{1}{\|x-Z\|}\Big(1-\ps{v}{ \frac{x-Z}{\|x-Z\|}}^2\Big)\Ind{Z\neq x}\bigg]
=
0
.
$$
We then have 
$$
\frac{1}{\|x-Z\|}\Big(1-\ps{v}{ \frac{x-Z}{\|x-Z\|}}^2\Big)\Ind{Z\neq x}
=
0
$$
$P$-almost surely. Because $P$ admits a density, we have $\|x-Z\|^{-1}\Ind{Z\neq x}\neq 0$ with $P$-probability $1$. Consequently, we have 
$$
\Big|\ps{v}{ \frac{x-Z}{\|x-Z\|}}\Big|
=
1
$$
with $P$-probability $1$. This implies that $P$ is supported on the line through $x$ with direction $v$, a contradiction. We deduce that $J_x\Fpg$ is positive definite, hence invertible. This concludes the proof.		
\end{Proof}

\section{Examples of reconstruction} \label{sec:Examples}

In this section we compute $\LL_d(\Fpg)$ when $P$ is a standard normal and a standard Cauchy distribution in $\R^2$ and $\R^3$. We show that the result coincides with the density of $P$, as established in Section \ref{sec:PDEGeneral}. Because the nature of the operator $\LL_d$ depends on whether $d$ is odd or even, computing $\LL_2$ and $\LL_3$ is requires different approaches.\\

Recall that the standard normal and standard Cauchy distributions on $\R^d$ have density 
$$
f(x)
=
\frac{1}{(2\pi)^{d/2}} e^{-\|x\|^2/2}
,
\quad \textrm{and} \quad 
f(x)
=
\frac{\Gamma(\frac{d+1}{2})}{\pi^{\frac{d+1}{2}}}
\frac{1}{(1+\|x\|^2)^{\frac{d+1}{2}}}
,
$$
respectively. Because these distributions are spherically symmetric, Proposition 2.2 (i) in \cite{GirStu2017} entails that there exists a function $g:[0,\infty)\to [0,\infty)$ such that their geometric cdf $\Fpg$ writes 
$$
\Fpg(x)=g(\|x\|)\frac{x}{\|x\|}
,\spa \forall\ x\in\R^d
.
$$
When $P$ is the standard normal distribution on $\R^d$, it was shown in Section 4.3 of \cite{Oja2010} that 
\begin{equation}\label{eq:RankGaussian}
g(r)
=
\frac{1}{\sqrt{2}}\frac{\Gamma(\frac{d+1}{2})}{\Gamma(\frac{d+2}{2})}\
re^{-r^2/2}\ 
%\exp(-r^2/2)\
{}_1F_1\Big( \frac{d+1}{2}; \frac{d+2}{2}; \frac{r^2}{2} \Big)
,\spa \forall\ r>0
,
\end{equation}
where~$_1F_1$ is the confluent hypergeometric function; see, e.g., (13.2.2) in~\cite{Olvetal2010}. When $P$ is the standard Cauchy distribution on $\R^d$, it was shown in Section 4.3 of \cite{Oja2010} that 
\begin{equation}\label{eq:RankCauchy}
g(r)
=
\frac{1}{\sqrt{\pi}}\frac{\Gamma(\frac{d+1}{2})}{\Gamma(\frac{d+2}{2})} \frac{r}{1+r^2} \
{}_2F_1\Big( \frac{d+1}{2},1; \frac{d+2}{2}; \frac{r^2}{1+r^2} \Big)
,\spa \forall\ r>0,
\end{equation}
where~$_2F_1$ is the Gaussian hypergeometric function; see, e.g., (15.2.1) in~\cite{Olvetal2010}.\\

Both when $d=2$ and $d=3$, we need to compute the divergence $\nabla\cdot \Fpg=\sum_{i=1}^d \partial_i (\Fpg)_i$ of $\Fpg$, where $\partial_i$ denotes the partial derivative with respect to the $i$th coordinate and $(\Fpg)_i$ stand for the $i$th component of $\Fpg$. A straightforward computation gives 
\begin{equation}\label{eq:DivRankSpherique}
(\nabla\cdot \Fpg)(x) = g'(\|x\|)+(d-1)\frac{g(\|x\|)}{\|x\|}
,\spa \forall\ x\in\R^d
.
\end{equation}

%By an obvious abuse of notations, we will write 
%$$
%\forall\ r>0,\spa
%f(r)
%=
%\frac{1}{(2\pi)^{d/2}} \exp(-r^2/2)
%,
%\quad \textrm{and} \quad 
%f(r)
%=
%\frac{\Gamma(\frac{d+1}{2})}{\pi^{\frac{d+1}{2}}}
%\frac{1}{(1+r^2)^{\frac{d+1}{2}}}
%.
%$$

\subsection{Dimension $3$}

Letting $h(\|x\|):=(\nabla\cdot \Fpg)(x)$ for any $x\in\R^3$, \eqref{eq:DivRankSpherique} yields 
$
h(r)
=
g'(r)+2g(r)/r$ for any $r>0$. Recalling that $\gamma_3 = (8\pi)^{-1}$ (see the beginning of Section \ref{sec:PDEGeneral}), a straightforward computation yields 
$$
(\LL_3 \Fpg)(x)
=
\gamma_3\ (-\Delta) \big(\nabla\cdot \Fpg\big)(x)
=
-\frac{1}{8\pi} \Big(h''(\|x\|)+2 \frac{h'(\|x\|)}{\|x\|}\Big)
,\spa \forall\ x\in\R^3
.
$$
Letting $e\in\SS^2$ be arbitrary, we then need to show that 
$$
-\frac{1}{8\pi} \Big(h''(r)+2 \frac{h'(r)}{r}\Big)
=
f(re)
,\spa \forall\ r>0.
$$
%where $h(r)=g'(r)+2 g(r)/r$.

\subsubsection{Trivariate Gaussian distribution}
\label{secGauss3}

Fix $r>0$. Taking $d=3$ in \eqref{eq:RankGaussian} yields
\begin{align*}
g(r)
&=
\frac{2\sqrt{2}}{3\sqrt{\pi}}\ r
e^{-r^2/2}\
{}_1F_1\Big( 2; \frac{5}{2}; \frac{r^2}{2} \Big)
\\[2mm]
&=
\frac{2\sqrt{2}}{3\sqrt{\pi}}\ r\  
{}_1F_1\Big( \frac{1}{2}; \frac{5}{2}; -\frac{r^2}{2} \Big)
,
\end{align*}
where the last equality follows from~(13.2.39) in~\cite{Olvetal2010}. Using~(13.3.2) in~\cite{Olvetal2010} with~$a=\frac{1}{2}$ and~$b=\frac{3}{2}$ gives
\begin{eqnarray*}
g(r)
&= & 
\frac{4\sqrt{2}}{3\sqrt{\pi}r}  
\bigg\{
\frac{r^2}{2}\, {}_1F_1\Big( \frac{1}{2}; \frac{5}{2}; -\frac{r^2}{2} \Big)
\bigg\}
\\[2mm]
&= & 
\frac{4\sqrt{2}}{3\sqrt{\pi}r}  
\bigg\{
\frac{3}{4}\, {}_1F_1\Big( \frac{1}{2}; \frac{1}{2}; -\frac{r^2}{2} \Big)
+
\frac{3}{2}\Big(\frac{r^2}{2}-\frac{1}{2} \Big)\, {}_1F_1\Big( \frac{1}{2}; \frac{3}{2}; -\frac{r^2}{2} \Big)
\bigg\}
\\[2mm]
&= & 
\frac{4\sqrt{2}}{3\sqrt{\pi}r}  
\bigg\{
\frac{3}{4}\exp(-r^2/2)
+
\frac{3}{4}(r^2-1)\, \frac{\sqrt{\pi}}{\sqrt{2}r} \textrm{ erf}(r/\sqrt{2})
\bigg\}
,
\end{eqnarray*}
where the last equality follows from~(13.6.1) and~(13.6.7) in~\cite{Olvetal2010}; here,
$$
\textrm{ erf}(r) 
= 
\frac{2}{\sqrt{\pi}}\int_0^r e^{-t^2/2}\,dt
=
2\Phi(\sqrt{2} r)-1
$$
is the error function (throughout, $\Phi$ and~$\phi$ stand for the cumulative distribution function and probability density function of the standard normal distribution, respectively). Thus,
\begin{eqnarray*}
g(r)
&= & 
\frac{\sqrt{2}}{\sqrt{\pi}r}  
e^{-r^2/2}
+
\frac{r^2-1}{r^2} \textrm{ erf}(r/\sqrt{2})
%&= & 
%\frac{\sqrt{2}}{\sqrt{\pi}s}  
%\exp(-s^2/2)
%+
%\frac{s^2-1}{s^2} (2\Phi(s)-1)
%\\[2mm]
\\[2mm]
&= & 
\frac{2}{r} \phi(r)  
+
\frac{r^2-1}{r^2} (2\Phi(r)-1)
%	\\[2mm]
%	&=&
%	2\Phi(r)-1 - \frac{2\Phi(r)-1-2r\phi(r)}{r^2}
\cdot
\end{eqnarray*}
Using the fact that $\phi'(r)=-r\phi(r)$, straightforward computations give
$$
g'(r)
=
2\ \frac{2\Phi(r)-1-2r\phi(r)}{r^3}
.
$$
%Noting that
%$$
%(2\Phi(r)-1-2r\phi(r))'
%=
%2\phi(r)-2\phi(r)-2r\phi'(r)
%=
%2r^2\phi(r)
%,
%$$
%for any $r>0$ this provides
%$$
%g'(r)
%=
%2\phi(r) - \frac{2r^4\phi(r)-2r(2\Phi(r)-1-2r\phi(r))}{r^4}
%=
%2 \frac{2\Phi(r)-1-2r\phi(r)}{r^3}
%,
%$$
%for any $r>0$. 
This yields
$$
h(r)
=  
g'(r) + 2 \frac{g(r)}{r}
=
\frac{2(2\Phi(r)-1)}{r} 
.
$$
%\begin{eqnarray*}
%	\forall\ r>0,\spa
%	h(r)
%	&= & 
%	g'(r) + 2 \frac{g(r)}{r}
%	\\[2mm]
%	&= & 
%	2 \frac{2\Phi(r)-1-2r\phi(r)}{r^3}
%	+
%	\frac{2}{r} 
%	\bigg(
%	2\Phi(r)-1 - \frac{2\Phi(r)-1-2r\phi(r)}{r^2}
%	\bigg)
%	\\[2mm]
%	&= & 
%	\frac{2(2\Phi(r)-1)}{r} 
%	\cdot
%\end{eqnarray*}
It follows that 
$
h'(r)
%=
%\frac{4r\phi(r)-2(2\Phi(r)-1)}{r^2} 
=
(4r\phi(r)-4\Phi(r)+2)/r^2
,
$
and
$
h''(r)
=
-4\phi(r)
-2\frac{4r\phi(r)-4\Phi(r)+2}{r^3} 
$
.
%\begin{eqnarray*}
%	h''(r)
%	&= & 
%	\frac{r^2(4\phi(r)+4r\phi'(r)-4\phi(r))-2r(4r\phi(r)-4\Phi(r)+2)}{r^4} 
%	\\[2mm]
%	&= & 
%	-4\phi(r)
%	-2\frac{4r\phi(r)-4\Phi(r)+2}{r^3} 
%	,
%\end{eqnarray*}
%for any $r>0$. 
Hence,
$$
-\frac{1}{8\pi} 
\bigg(
h''(r)+\frac{2}{r} h'(r)
\bigg)
=
\frac{\phi(r)}{2\pi}
=
\frac{1}{(2\pi)^{3/2}} 
e^{-r^2/2} 
=
f(re)
,
$$
which, as expected, coincides with the probability density function of the trivariate standard normal distribution.

\subsubsection{Trivariate Cauchy distribution}
\label{secCauchy3}

Fix $r>0$. Taking $d=3$ in \eqref{eq:RankCauchy} yields
\begin{align*}
g(r)
&=
\frac{4r}{3\pi(1+r^2)} 
\,
{}_2F_1\Big( 2,1; \frac{5}{2}; \frac{r^2}{1+r^2} \Big)
\\[2mm]
&=
\frac{4r}{3\pi\sqrt{1+r^2}} 
\,
{}_2F_1\Big( \frac{1}{2},\frac{3}{2}; \frac{5}{2}; \frac{r^2}{1+r^2} \Big)
,
\end{align*}
where we used~(15.8.1) in~\cite{Olvetal2010}. Applying Identity~92 on page~473 of~\cite{Prud90} then provides
\begin{align*}
g(r)
&=
\frac{2}{\pi r^2} 
\Big( (1+r^2) \arcsin\Big( \frac{r}{\sqrt{1+r^2}}\Big)  - r \Big) 
\\[2mm]
&=
\frac{2}{\pi r^2} 
\Big( (1+r^2) \arctan(r) - r \Big) 
.
\end{align*}
Direct computations then yield
$$
h(r)
=
g'(r) + 2 \frac{g(r)}{r}
=
\frac{4\arctan(r)}{\pi r}  
,
$$
hence
$$
-\frac{1}{8\pi} 
\bigg(
h''(r)+\frac{2}{r} h'(r)
\bigg)
=
\frac{1}{\pi^2(1+r^2)}
=
f(re)
,
$$
which coincides with the probability density function of the trivariate standard Cauchy distribution.

\subsection{Dimension $2$}

Recall that $(-\Delta)^{1/2}u$ is defined through
$$
((-\Delta)^{\frac12}u)(x)
=
2\pi\ \mathcal{F}^{-1}\big(\|\xi\|\mathcal{F}u(\xi)\big)(x)
,\spa \forall x\in\R^2;
.
$$
Because $\Gamma(3/2)=\sqrt{\pi}/2$, we have $\gamma_2=(2\pi)^{-1}$. Letting $u=\nabla\cdot \Fpg$, it follows that
$$
(\LL_2 \Fpg)(x)
=
\gamma_2\ (-\Delta)^{\frac12}(\nabla\cdot \Fpg)(x)
=
\mathcal{F}^{-1}\big(\|\xi\|\mathcal{F}u(\xi)\big)(x)
$$
for all $x\in\R^2$. Writing $u(x)=h(\|x\|)$, (\ref{eq:DivRankSpherique}) yields $h(r)=g'(r)+g(r)/r$ for any $r>0$. A straightforward computation gives
$$ 
(\F u)(\xi)
=
\int_0^\infty 
h(r)
\Bigg(
\int_0^{2\pi}
e^{-i (2\pi r\|\xi\|)\cos \theta}
\, d\theta 
\Bigg)
r
\,
dr
=
2\pi
\int_0^\infty 
h(r)
J_0(2\pi r\|\xi\|)
r
\,
dr
,
$$
for all $\xi\in\R^2$, where
$$
J_0(z)
=
(2\pi)^{-1}
\int_0^{2\pi}
e^{-i z \cos \theta}
\, d\theta 
$$
is the Bessel function of the first kind with order zero. Writing~$\tilde{h}$ for the function defined through~$\tilde{h}(r)=\sqrt{r}\ h(r)$, we thus have 
\begin{align*}
(\F u)(\xi)
&=
\sqrt{\frac{2\pi}{\|\xi\|}}
\int_0^\infty 
\sqrt{r}\ h(r)
J_0\big(r (2\pi\|\xi\|)\big)
\sqrt{r(2\pi\|\xi\|)}
\,
dr
\\[2mm]
&=
\sqrt{\frac{2\pi}{\|\xi\|}}
(\mathcal{H}_0 \tilde{h})(2\pi\|\xi\|)
,
\end{align*}
where  
$$
(\mathcal{H}_0 \phi)(r)
:=
\int_0^\infty 
\phi(s)
J_0(sr)
\sqrt{sr}
\,
ds
$$
is the Hankel transform of~$\phi$ with order zero; see, e.g., page~1 in~\cite{Ober72}. Finally, to compute $\F^{-1}(\|\xi\| \F u(\xi))(x)$ we will use the fact that the restriction of $\F$ and $\F^{-1}$ to spherically symmetric functions coincide.

%%%%%%%%%%%%%%%%%%%%%%%%%%%%%%%%%%

\subsubsection{Bivariate Gaussian distribution}
\label{secGauss2}

Fix $r>0$. Taking $d=2$ in \eqref{eq:RankGaussian} yields
\begin{align*}
g(r)
&=
\frac{\sqrt{\pi}}{2\sqrt{2}}\ 
r e^{-r^2/2}\
{}_1F_1\Big( \frac{3}{2}; 2; \frac{r^2}{2} \Big)
\\[2mm]
&=
\frac{\sqrt{\pi}}{\sqrt{2}r}
\bigg\{
\frac{r^2}{2}
e^{-r^2/2}\
{}_1F_1\Big( \frac{3}{2}; 2; \frac{r^2}{2} \Big)
\bigg\}
.
\end{align*}
Hence, applying (13.3.21) in~\cite{Olvetal2010}
%, then (10.25.2),  
provides
\begin{align*}
g'(r)
&=
-\frac{\sqrt{\pi}}{2\sqrt{2}}
e^{-r^2/2}\
{}_1F_1\Big( \frac{3}{2}; 2; \frac{r^2}{2} \Big)
+
\frac{\sqrt{\pi}}{\sqrt{2}r}
\bigg\{
e^{-r^2/2}\
{}_1F_1\Big( \frac{1}{2};1; \frac{r^2}{2} \Big)
\bigg\}
r
\\[2mm]
&=
-
\frac{g(r)}{r}
+\frac{\sqrt{\pi}}{\sqrt{2}}
e^{-r^2/2}\
{}_1F_1\Big( \frac{1}{2};1; \frac{r^2}{2} \Big)
.
\end{align*}
Therefore, (13.6.9) in~\cite{Olvetal2010} yields
\begin{align*}
h(r)
&=
g'(r) + \frac{g(r)}{r}
=
\frac{\sqrt{\pi}}{\sqrt{2}}
e^{-r^2/2}\
{}_1F_1\Big( \frac{1}{2};1; \frac{r^2}{2} \Big)
\\[2mm]
&=
\frac{\sqrt{\pi}}{\sqrt{2}}
e^{-r^2/4}
\mathcal{I}_0\Big(\frac{r^2}{4}\Big)
,
\end{align*}
where $\mathcal{I}_0$ is the modified Bessel function of order $0$. Using~(2.126) in~\cite{Ober72} with~$a=1/4$, we obtain that~$u(x)=h(\|x\|)$ satisfies
\begin{align*}
(\F u)(\xi)
&=
\sqrt{\frac{2\pi}{\|\xi\|}}
(\mathcal{H}_0 \tilde{h})(2\pi\|\xi\|)
\\[2mm]
&=
\sqrt{\frac{2\pi}{\|\xi\|}}
\frac{\sqrt{\pi}}{\sqrt{2}}
\frac{1}{\sqrt{(\pi/2)2\pi\|\xi\|}}
e^{-(2\pi\|\xi\|)^2/2}
\\[2mm]
&=
\frac{1}{\|\xi\|}
e^{-2\pi^2\|\xi\|^2}
,
\end{align*}
so that 
$
\|\xi\| (\F u)(\xi)
=
e^{-2\pi^2\|\xi\|^2}
$. 
Using~(2.23)  in~\cite{Ober72} with~$a=2\pi^2$, we then have 
\begin{align*}
\mathcal{F}^{-1}\big(\|\xi\|\mathcal{F}u(\xi)\big)(x)
&=
\sqrt{\frac{2\pi}{\|x\|}}
\bigg(
\frac{1}{4\pi^2}
\sqrt{2\pi\|x\|}
e^{-4\pi^2\|x\|^2/(8\pi^2)}
\bigg)
\\[2mm]
&=
\frac{1}{2\pi}
e^{-\|x\|^2/2}
.
\end{align*}
Thus,
$$
(\LL_2 \Fpg)(x)
=
\mathcal{F}^{-1}(\|\xi\|\mathcal{F}u(\xi))(x)
=
\frac{1}{2\pi}
e^{-\|x\|^2/2}
,\spa \forall x\in\R^2
,
$$
which is indeed the probability density function of the bivariate standard normal distribution.

%%%%%%%%%%%%%%%%%%%%%%%%%%%%%%%%%%

\subsubsection{Bivariate Cauchy distribution}
\label{secCauchy2}

Fix $r>0$. Taking $d=2$ in \eqref{eq:RankCauchy} yields
$$
g(r)
=
\frac{r}{2(1+r^2)} \,
{}_2F_1\Big( \frac{3}{2},1; 2; \frac{r^2}{1+r^2} \Big)
.
$$
Hence, applying~(15.8.1) in~\cite{Olvetal2010}, then Identity~84 on page~473 of~\cite{Prud90}, provides
$$
g(r)
=
\frac{r}{2\sqrt{1+r^2}} \,
{}_2F_1\Big( \frac{1}{2},1; 2; \frac{r^2}{1+r^2} \Big)
=
\frac{r}{1+\sqrt{1+r^2}}
\cdot
$$
Therefore, direct computation yields
$$
h(r)
=
g'(r) + \frac{g(r)}{r}
=
\frac{1}{\sqrt{1+r^2}}
.
$$
Now, using~(2.19) in~\cite{Ober72} with~$a=1$, we obtain that~$u(x)=h(\|x\|)$ satisfies
$$
(\F u)(\xi)
=
\sqrt{\frac{2\pi}{\|\xi\|}}
(\mathcal{H}_0 \tilde{h})(2\pi\|\xi\|)
=
\sqrt{\frac{2\pi}{\|\xi\|}}
\frac{1}{\sqrt{2\pi\|\xi\|}}
e^{-2\pi\|\xi\|}
,
$$
for all $\xi\in\R^2$, so that 
$
\|\xi\| \mathcal{F}u(\xi)
=
e^{-2\pi\|\xi\|}
$. 
Thus, using~(2.23)  in~\cite{Ober72} with~$a=2\pi$, we have 
\begin{align*}
\mathcal{F}^{-1}\big(\|\xi\|\mathcal{F}u(\xi)\big)(x)
&=
\sqrt{\frac{2\pi}{\|x\|}}
\bigg(
\frac{2\pi\sqrt{2\pi\|x\|}}{(4\pi^2+(2\pi\|x\|)^2)^{3/2}}
\bigg)
\\[2mm]
&=
\frac{1}{2\pi(1+\|x\|^2)^{3/2}}
.
\end{align*}
It follows that
$$ 
(\LL_2 \Fpg)(x)
=
\mathcal{F}^{-1}(\|\xi\|\mathcal{F}u(\xi))(x)
=
\frac{1}{2\pi(1+\|x\|^2)^{3/2}}
,\spa \forall x\in\R^2
,
$$
which is the probability density function of the bivariate standard Cauchy distribution.

\end{document}